\documentclass[article,11pt]{amsart}
\usepackage{amstext,amsfonts,amsmath,tabu,graphicx,tikz,stmaryrd, amssymb,MnSymbol}
\usetikzlibrary[calc,decorations.pathreplacing]
\usepackage[utf8]{inputenc}
\usepackage[english]{babel}
\usepackage[T1]{fontenc}
\usepackage[a4paper]{geometry}
\usepackage{hyperref,url,calc,comment,xspace,ifthen}
\usepackage{enumitem}
\usepackage{multirow}
\usepackage[width=\linewidth]{caption}
\usepackage{overpic}
\usepackage{subcaption}
\usepackage{pict2e} 
\usetikzlibrary{arrows}
\usepackage[backend=bibtex, natbib=true, language=english,
style=authoryear,
style=numeric,
maxcitenames=1 ,maxbibnames=10
]{biblatex}
\bibliography{paper1.bib} 

\DeclareSourcemap{
	\maps[datatype=bibtex]{
		\map{
			\step[fieldsource=note, final]
			\step[fieldset=addendum, origfieldval, final]
			\step[fieldset=note, null]
		}
	}
}

\definecolor{uuuuuu}{rgb}{0.26666666666666666,0.26666666666666666,0.26666666666666666}
\definecolor{qqqqff}{rgb}{0.,0.,1.}
\definecolor{ffqqqq}{rgb}{1.,0.,0.}
\definecolor{ududff}{rgb}{0.30196078431372547,0.30196078431372547,1.}
\definecolor{zzttqq}{rgb}{0.6,0.2,0.}
\definecolor{geogebraBlue}{rgb}{0.301961,0.301961,1.0}
\definecolor{geogebraBrown}{rgb}{0.6,0.2,0}
\definecolor{geogebraGreen}{rgb}{0, 0.39216, 0}
\definecolor{geogebraPurple}{rgb}{0.498039,0,1}

\newcommand{\R}{\mbox{$\mathbb{R}$}}
\renewcommand{\S}{\mbox{$\mathbb{S}$}}

\newcommand{\Z}{\mbox{$\mathbb{Z}$}}
\newcommand{\Pro}{\mbox{$\mathbb{P}$}}
\newcommand{\Mink}{\mbox{$\mathbb{P}\mathbb{R}^{4,1}$}}
\newcommand{\LightCone}{\mbox{$\mathcal{L}$}}

\newcommand{\tr}{\ensuremath{\operatorname{trace}}}
\newcommand{\Mat}[1]{\ensuremath{\mathbf{#1}}}

\newcommand{\area}{\ensuremath{\operatorname{area}}}
\newcommand{\DV}{\ensuremath{\operatorname{cr}}}

\newcommand{\innerprod}[2]{\ensuremath{\langle #1 , #2 \rangle}}

\def\fe{\mbox{\boldmath$e$}}

\def\fw{\mbox{\boldmath$w$}}

\newcommand{\quadf}{\ensuremath{\mathcal{Q}}} 
\newcommand{\net}{\ensuremath{f}} 
\newcommand{\nesh}{\ensuremath{g}} 
\newcommand{\dgrad}[1][]{\ensuremath{\ifthenelse{\equal{#1}{}}{\mathcal{D}}{\mathcal{D}#1}}} 
\newcommand{\discdifu}[1]{\ensuremath{\delta_u #1}}
\newcommand{\discdifv}[1]{\ensuremath{\delta_v #1}}
\newcommand{\discdif}[2]{\ensuremath{\delta_{#1} #2}}

\newcommand{\cbp}[1][]{\ensuremath{\ifthenelse{\equal{#1}{}}{c}{c_{#1}}}}
\newcommand{\whiteface}[1][]{\ensuremath{\ifthenelse{\equal{#1}{}}{\mathcal{W}}{\mathcal{W}_{#1}}}}
\newcommand{\blackface}[1][]{\ensuremath{\ifthenelse{\equal{#1}{}}{\mathcal{B}}{\mathcal{B}_{#1}}}}
\newcommand{\gaussimage}[1][]{\ensuremath{\ifthenelse{\equal{#1}{}}{n}{n_{#1}}}}
\newcommand{\faceGaussimage}{\ensuremath{N}}
\newcommand{\shape}{\ensuremath{S}} 
\newcommand{\shapeMat}{\ensuremath{\Sigma}}
\newcommand{\moeb}[1][]{\ensuremath{\ifthenelse{\equal{#1}{}}{s}{s_{#1}}}}

\newcommand{\suppl}[2]{\ensuremath{\ell(#1,#2)}}
\newcommand{\firstFundi}{\ensuremath{\operatorname{I}}}
\newcommand{\secondFundi}{\ensuremath{\operatorname{II}}}

\theoremstyle{plain}
\newtheorem{satz}{Theorem}

\newtheorem{lem}[satz]{Lemma}
\newtheorem{kor}[satz]{Corollary}
\theoremstyle{definition}
\newtheorem{defi}[satz]{Definition}
\newtheorem{defiLemma}[satz]{Lemma and Definition}
\theoremstyle{remark}
\newtheorem{bem}[satz]{Remark}

\newenvironment{bew}[1][]{
	\ifthenelse{\equal{#1}{}}{\begin{proof}}{\begin{proof}[#1]}
}
{
	\end{proof}
}

\newcounter{count}

\title{Discrete Isothermic Nets Based on Checkerboard Patterns}
\author{Felix Dellinger\\
	{\tiny 
	Institute of Geometry, TU Graz. Kopernikusgasse 24, 8010 Graz, Austria\\
	Institut für diskrete Mathematik und Geometrie, TU Wien, Wiedner Hauptstrasse 8-10, A-1040 Vienna, Austria}}

\begin{document}
	
\maketitle

\section{Abstract}

This paper studies the discrete differential geometry of the checkerboard pattern inscribed in a quadrilateral net by connecting edge midpoints. It turns out to be a versatile tool which allows us to consistently define principal nets, Koenigs nets and eventually isothermic nets as a combination of both.

Principal nets are based on the notions of orthogonality and conjugacy and can be identified with sphere congruences that are entities of Möbius geometry. Discrete Koenigs nets are defined via the existence of the so called conic of Koenigs. We find several interesting properties of Koenigs nets, including their being dualizable and having equal Laplace invariants. 

Koenigs nets that are also principal are defined as isothermic nets. We prove that the class of isothermic nets is invariant under both dualization and Möbius transformations. Among other things, this allows a natural construction of discrete minimal surfaces and their Goursat transformations.

\section{Introduction}

Discretizing principal curvature nets is of great interest not only from a differential geometric point of view, but also in geometry processing, computer graphics and even freeform architecture \cite{pottmann-2007-pm, liu-Pottmann-conical}. The most prominent versions of discrete principal nets are circular nets and conical nets \cite{BobenkoBuch, liu-Pottmann-conical}. A new discretization, suggested by \citep{CheckerboardIntro}, is based on the checkerboard pattern inscribed in a quadrilateral net constructed by connecting edge midpoints. This approach has already proven to be useful in various applications \cite{Jiang2020Apl1, weingarten-aag-2020, UsingIsometries}. Its effectiveness suggests that there is more to the concept than just the good numerical approximation qualities already hinted at by \citep{CheckerboardIntro}. Indeed, we find that a rich discrete theory can be built upon these checkerboard patterns.

A checkerboard pattern is a quadrilateral net where every second face is a parallelogram. The edges of these parallelograms can be seen as discrete derivatives. If the faces in between the parallelograms are all planar we speak of a conjugate checkerboard pattern. If additionally the parallelograms are all rectangles we speak of a principal checkerboard pattern. As the concept of checkerboard patterns is Euclidean in nature, it is surprising that we can show principal nets to be Möbius invariant if they are seen as sphere congruences \cite{techter2021discrete}. Lifting these sphere congruences to the projective model of Möbius geometry preserves principality and offers the appropriate environment to efficiently study these geometric objects.

For a net with planar faces the supporting lines of neighboring edges intersect. Every face can be associated with six such intersection points. In \cite{DoliwasKoenigsNets} discrete Koenigs nets have been characterized by the property that these six points lie on a common conic section, the so called conic of Koenigs \cite{KoenigsOriginalPaper}. We apply this definition to a checkerboard pattern. The resulting discrete Koenigs nets enjoy several interesting properties such as projective invariance and the existence of dual nets similar to the approach in \cite{KoenigsNetsBobenkoSuris}. Usually, Koenigs nets have been known as \emph{nets with equal Laplace invariants}. While this property has been lost with previous discretizations of Koenigs nets, we manage to retain it in a natural way.

We define discrete isothermic nets as discrete Koenigs nets that are also principal. Analogous to the classical smooth theory, the class of discrete isothermic nets is invariant under both dualization and Möbius transformations. This is not only interesting from a theoretical point of view, but also offers a practical way to define and construct discrete minimal surfaces as surfaces that are dual to their own Gauß image. Consequently the dual of any isothermic net on the unit sphere can be seen as a minimal surface. All of these steps can now be easily discretized with our approach.


\section{Checkerboard patterns}\label{Section:CheckerboardPattern}

\subsection{Preliminaries.}

\begin{figure}
	\begin{subfigure}[t]{0.45\linewidth}
		\begin{overpic}[width=\linewidth,tics=10]
			{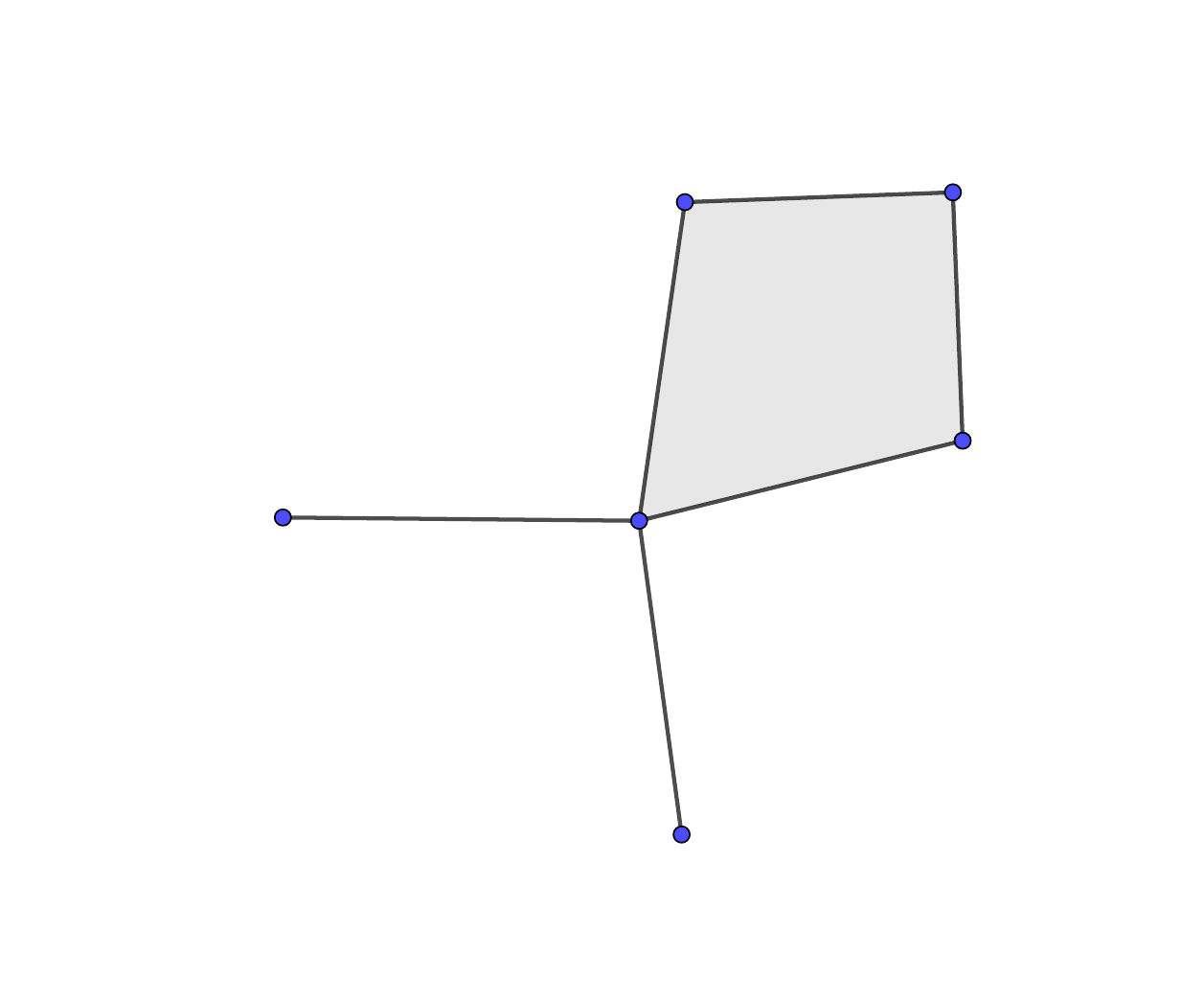}
			\put(54,42){\textcolor{geogebraBlue}{$f$}}
			\put(81,47){\textcolor{geogebraBlue}{$f_1$}}
			\put(58,14){\textcolor{geogebraBlue}{$f_{\bar{2}}$}}
			\put(23,43){\textcolor{geogebraBlue}{$f_{\bar{1}}$}}
			\put(54,69){\textcolor{geogebraBlue}{$f_{2}$}}
			\put(80,69){\textcolor{geogebraBlue}{$f_{12}$}}
			\put(63,52){$\quadf_f$}		
		\end{overpic}
	\end{subfigure}
	\begin{subfigure}[t]{0.45\linewidth}
		\includegraphics[width=\linewidth]{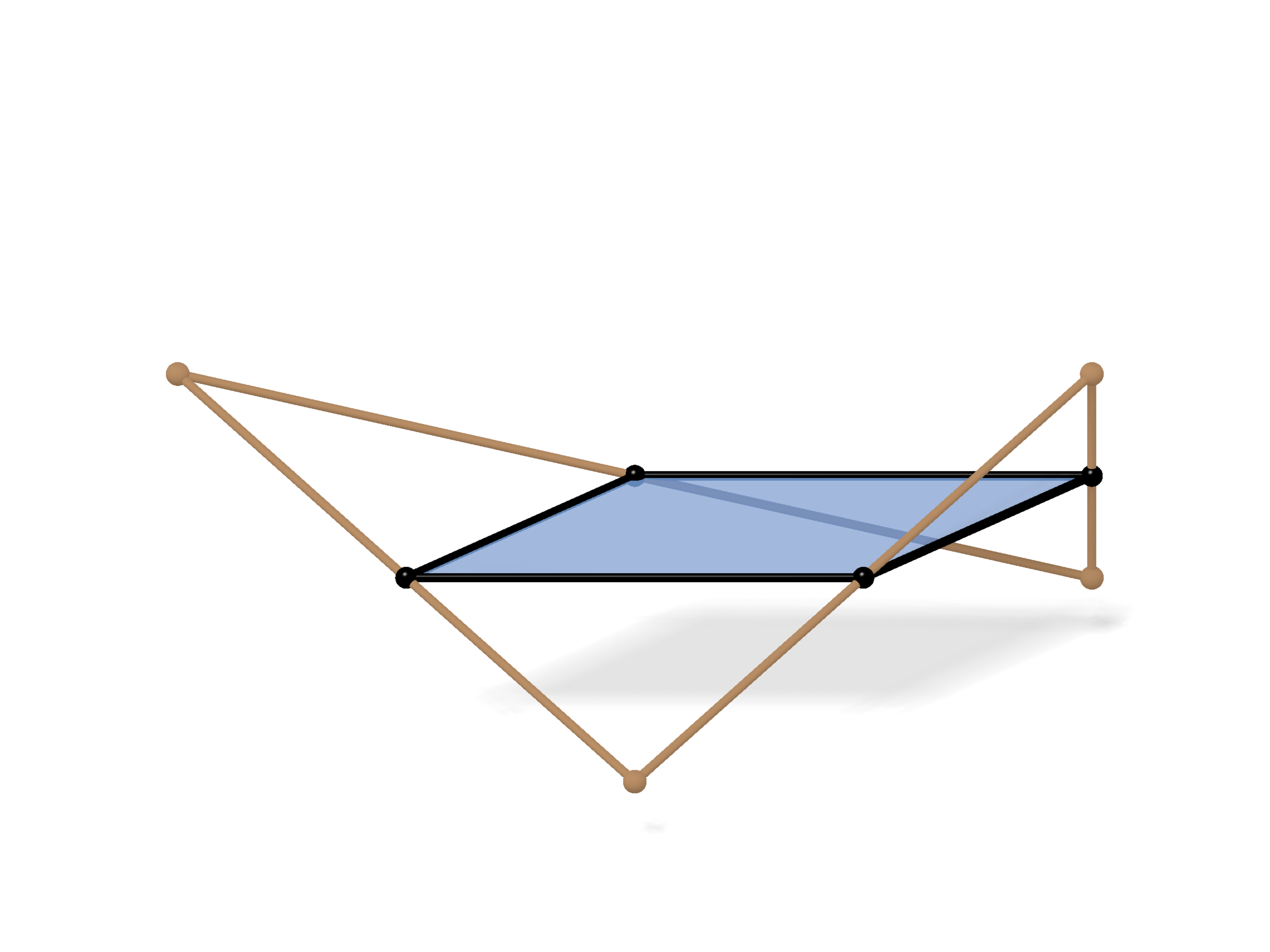}
	\end{subfigure}
	\vskip -10mm
	\caption{\textit{Left:} Notation for vertices and faces. \textit{Right:} An inscribed first order face, which is always a parallelogram.}
	\label{Fig:PointNotationAndInscribedFace}
\end{figure}
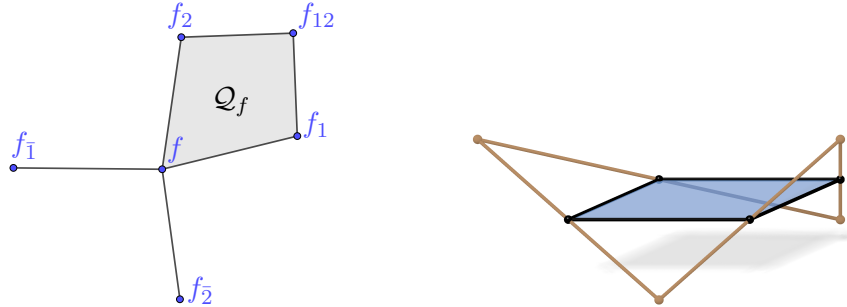

In this paper we study two-dimensional nets $f : D \to \R^3$. All our constructions are local which is why we can always assume $D=\Z^2$. To denote the one-ring or two-ring neighborhood of a vertex $f(k,l)$ we use the \emph{shift notation} as can be seen in Figure \ref{Fig:PointNotationAndInscribedFace}, left. 
The index $i$ resp. $\bar{i}$ indicates that the $i$-th coordinate is increased resp. decreased by one with $i \in \{1,2\}$. For instance, 
\begin{align*}
f_{1 \phantom{2}}(k,l) &= f(k+1,l), \quad f_{\bar{2}\phantom{2}}(k,l) = f(k,l-1), \\
f_{\bar{2}\bar{2}}(k,l) &= f(k,l - 2), \quad f_{12}(k,l) = f(k + 1, l + 1).
\end{align*}

We call the images of $f$ the vertices and the pairs $(f,f_1)$ or $(f,f_2)$ the edges of the net. Further we denote by $\quadf_f$ the face $(f, f_1, f_{12}, f_2)$. If no confusion can arise, we drop the index and just write $\quadf$.

\begin{defi}
	A \emph{checkerboard pattern} is a regular quad net where every second face is a parallelogram: $\quadf_{f}(k,l)$ is a parallelogram if $k+l \equiv 0\ (\text{mod}\ 2)$.
\end{defi}

Even if at first glance the definition of checkerboard patterns seems quite restrictive, they are actually very natural objects. From any given net $\net$ we can easily construct a checkerboard pattern $\cbp[\net]$ by midpoint subdivision as described in \cite{CheckerboardIntro}: The vertices of $\cbp[\net]$ are the edge midpoints of $\net$. There are then two kinds of faces in $\cbp[\net]$. The first type of face is formed by the midpoints of edges of each face $\quadf$ of $\net$ (compare Figure \ref{Fig:PointNotationAndInscribedFace}, right). It is elementary that these faces are parallelograms whose edges are parallel to the two diagonals of $\quadf$. We will refer to them as first order faces, as their edges can be interpreted as discrete first order derivatives. We denote the first order face associated to the quadrilateral $\quadf_f(k,l)$ by $\blackface[f](k,l)$.

\begin{figure}[h]
	\begin{subfigure}[t]{0.49\linewidth}
		\begin{overpic}[width=\linewidth]{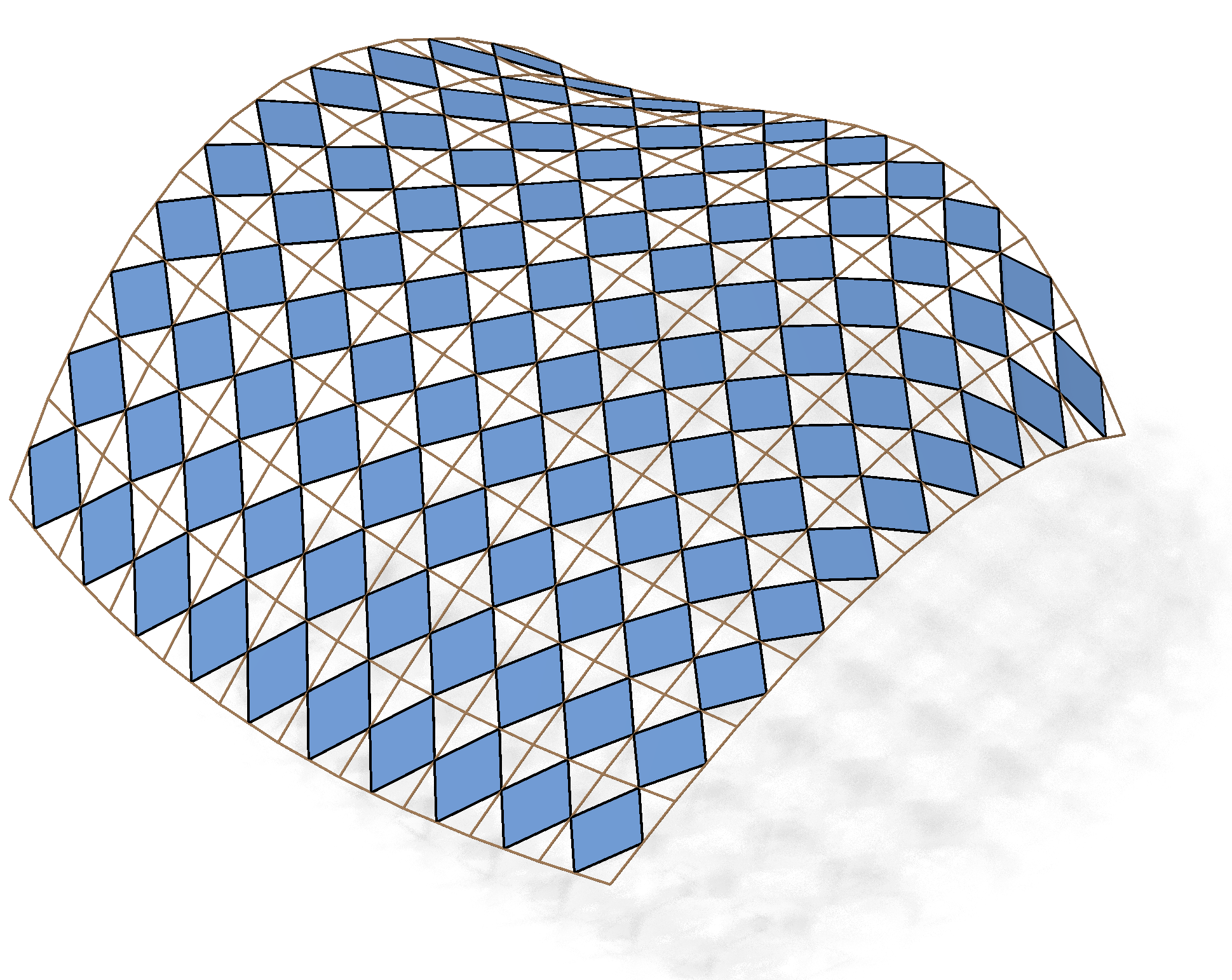}
		\end{overpic}
	\end{subfigure}
	\begin{subfigure}[t]{0.49\linewidth}
		\begin{overpic}[width=\linewidth]{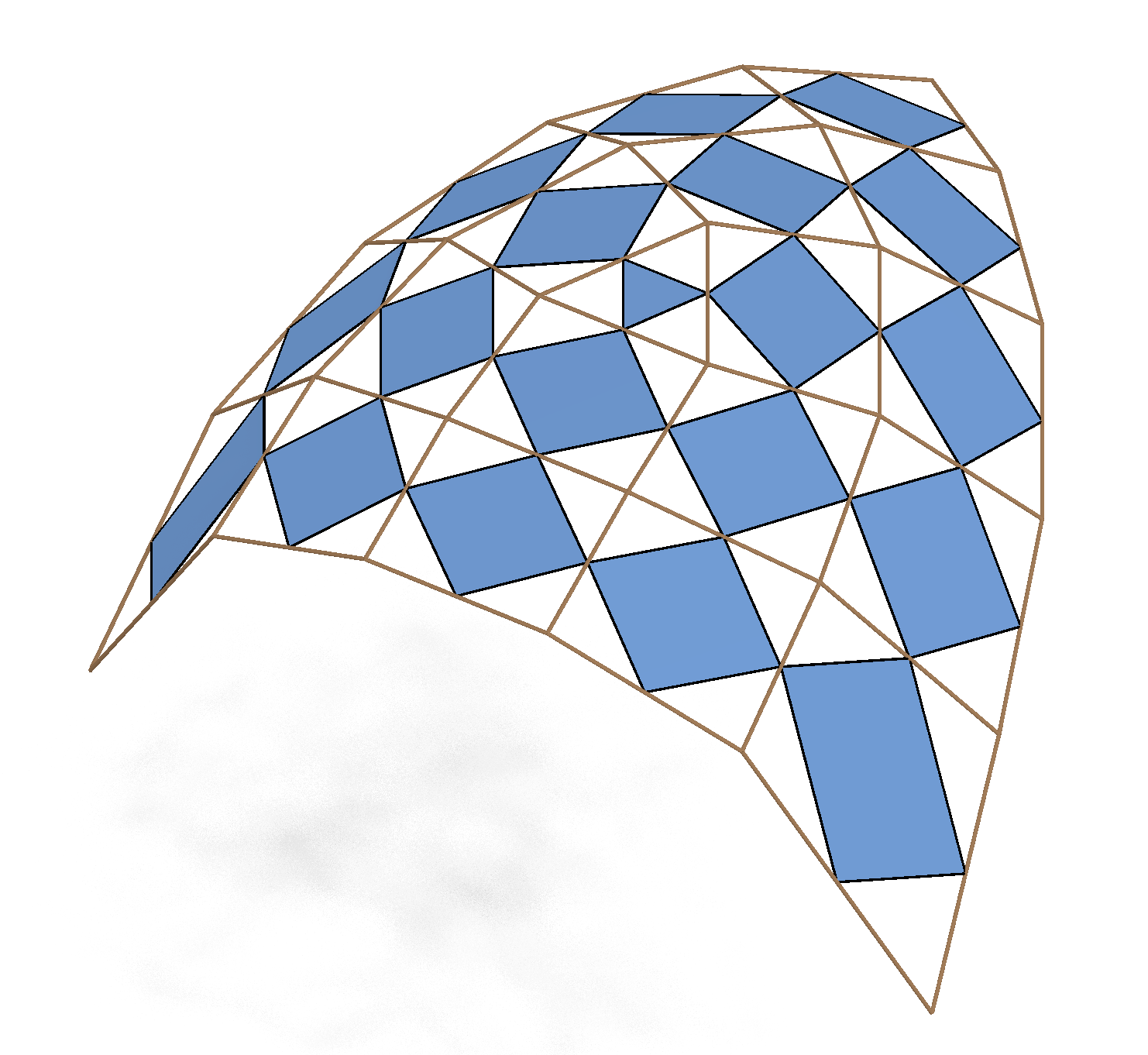}
			\put(52,46.5){$\whiteface$}
			\put(63,50){$\blackface$}
			\put(68,39){$\whiteface_1$}
			\put(58,59){$\whiteface_2$}
			\put(80,40){$\blackface_1$}
			\put(67,63){$\blackface_2$}
		\end{overpic}
	\end{subfigure}
	\caption{\textit{Left:} The first order faces of the checkerboard pattern are the blue parallelograms $\blackface[\net]$ inscribed in the faces of the control net $\net$. The white quadrilaterals $\whiteface[\net]$ in between are the second order faces. \textit{Right:} Control net and associated checkerboard pattern with a combinatorial singularity.}
	\label{Fig:CBPAndCBPWithSing}
\end{figure}
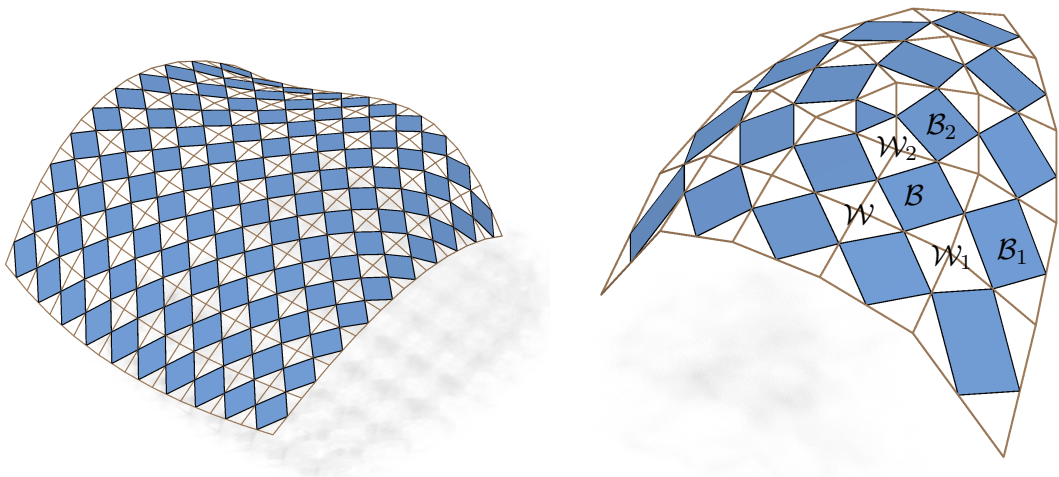

The second type of face is formed by the midpoints of edges emanating from a common vertex of $\net$. Those faces are, in general, non-planar quadrilaterals. We will refer to them as second order faces, because we associate properties related to second order derivatives with them. The second order face associated to the vertex $f(k,l)$ will be denoted by $\whiteface[f](k,l)$, compare Figure \ref{Fig:CBPAndCBPWithSing}.

If no confusion can arise we will drop the index $\net$ in all quantities. Following \citet{CheckerboardIntro} we call $\cbp[\net]$ the checkerboard pattern of $\net$ and $\net$ the control net of $\cbp[\net]$, see Figure \ref{Fig:CBPAndCBPWithSing}. Note that for a given checkerboard pattern there is a three parameter family of control nets. A control net is uniquely determined after the choice of an initial vertex as all other vertices can be obtained through iterated reflection at the vertices of the checkerboard pattern.

\begin{bem}
	The checkerboard pattern approach can be extended to nets with combinatorial singularities. For each $n$-gon, midpoint subdivision creates an inscribed $n$-gon, see e.g. an inscribed triangle in Figure \ref{Fig:CBPAndCBPWithSing}, right.
\end{bem}


For $\epsilon>0$, let the net $\net: \Z^2 \to \R^3$ sample a smooth surface parametrization $\phi: \R^2 \to \R^3$ i.e., $f(k,l) = \phi(\epsilon k,\epsilon l)$. We define the directions $u=\frac{1}{\sqrt{2}}(1,1)^T$ and $v=\frac{1}{\sqrt{2}}(-1,1)^T$. Intuitively speaking, the parameter lines of $f$ and $\cbp[f]$ enclose an angle of $45$ degrees. So, we can think of $\cbp[f]$ as being parameterized along the directions $u$ and $v$ in the coordinate plane.

The \emph{edge vectors} $\discdifu{f}:=\frac{1}{\sqrt{2}}(f_{12} - f)$ and $\discdifv{f}:=\frac{1}{\sqrt{2}}(f_2-f_1)$ of $\blackface[f](k,l)$ approximate the directional derivatives $\partial_u \phi$ and $\partial_v \phi$ at the point $\left(\epsilon (k+\frac{1}{2}),\epsilon (l+\frac{1}{2})\right)$ up to second order. Indeed,
\begin{align*}
	\frac{1}{\epsilon}\discdifu{f}(k, l) &= \partial_u \phi \Big(\epsilon \big(k+\frac{1}{2}\big),\epsilon \big(l+\frac{1}{2}\big)\Big) + \mathcal{O}(\epsilon^2),\\
	\frac{1}{\epsilon}\discdifv{f}(k, l) &= \partial_v \phi\Big(\epsilon \big(k+\frac{1}{2}\big),\epsilon \big(l+\frac{1}{2}\big)\Big) + \mathcal{O}(\epsilon^2),
\end{align*}
as a simple Taylor expansion shows. Moreover it can be shown by Taylor expansion that the difference of opposite edge vectors in a second order face $\whiteface[f](k,l)$ approximates $\partial_{uv} \phi(\epsilon k, \epsilon l)$ by first order. This motivates the notation of $\discdifu{f}$ and $\discdifv{f}$ for the edge vectors of $\blackface[f]$ and gives rise to the following definition.

\begin{defi}\label{Defi:OrthoConjuPrinciNet}
	We call a checkerboard pattern \emph{orthogonal} if its first order faces are rectangles. We call it \emph{conjugate} if its second order faces are planar. A checkerboard pattern is \emph{principal} if it is both conjugate and orthogonal.
\end{defi}

\begin{figure}
	\begin{overpic}[width=0.45\linewidth,trim= 5cm 2cm 6cm 9cm, clip]{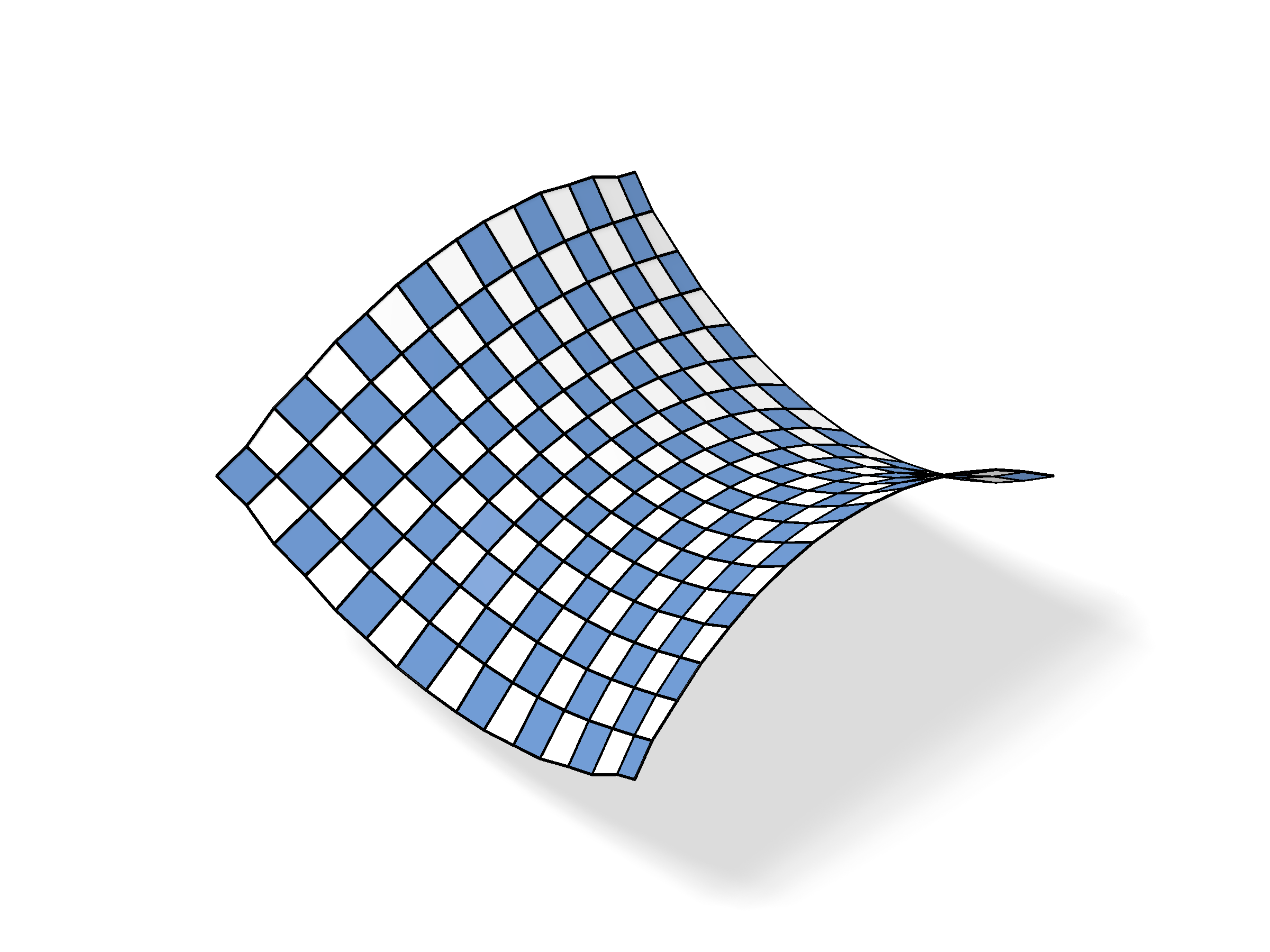}
	\end{overpic}
%
	\begin{overpic}[width=0.45\linewidth,tics=10]
	{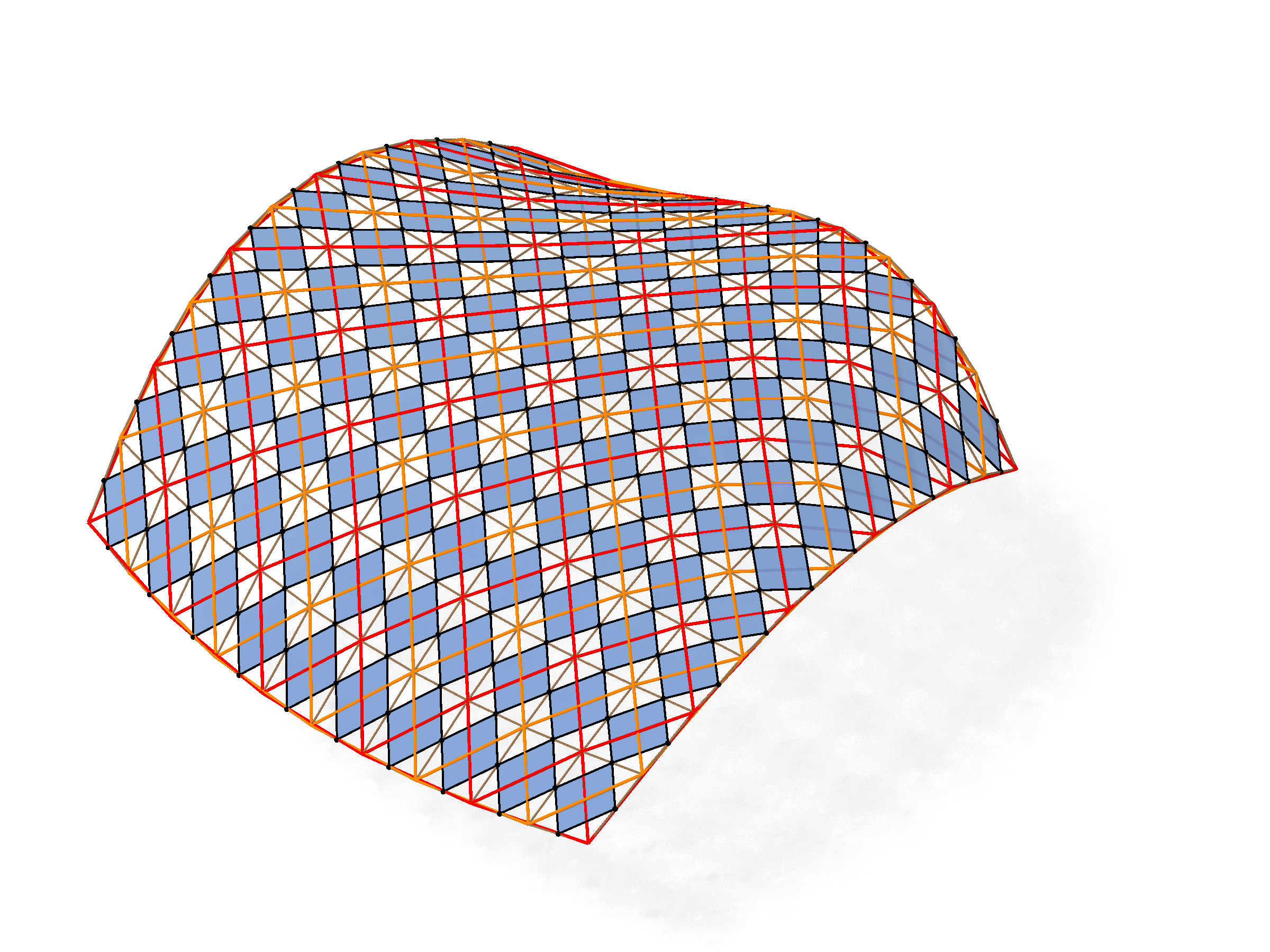}
	\end{overpic}
\vskip -4mm
\caption{\textit{Left}: A principal checkerboard pattern. All the white faces are planar and all blue faces are rectangles. \textit{Right:} The two nets defined by the diagonals of the control net have planar faces if and only if the checkerboard pattern is conjugate.}
\label{Fig:PrincipalNetAndDiagonalNets}
\end{figure}

\begin{bem}
	Conjugacy of a checkerboard pattern $\cbp[f]$ is already determined by its control net $\net$ and so are orthogonality and principality. Indeed, second order faces of $\cbp[f]$ are planar if and only if the two nets defined by the diagonals of $f$ have planar faces, compare Figure \ref{Fig:PrincipalNetAndDiagonalNets}, right. Thus the class of conjugate checkerboard patterns is invariant under projective transformations applied to the vertices of the control net.  
\end{bem}

\section{Curvature Theory}\label{Section:ShapeOperator}
\begin{figure}[ht]
	\begin{subfigure}[l]{0.49\linewidth}
		\includegraphics[trim = 10cm 0cm 5cm 10cm, clip,
		width=\linewidth]{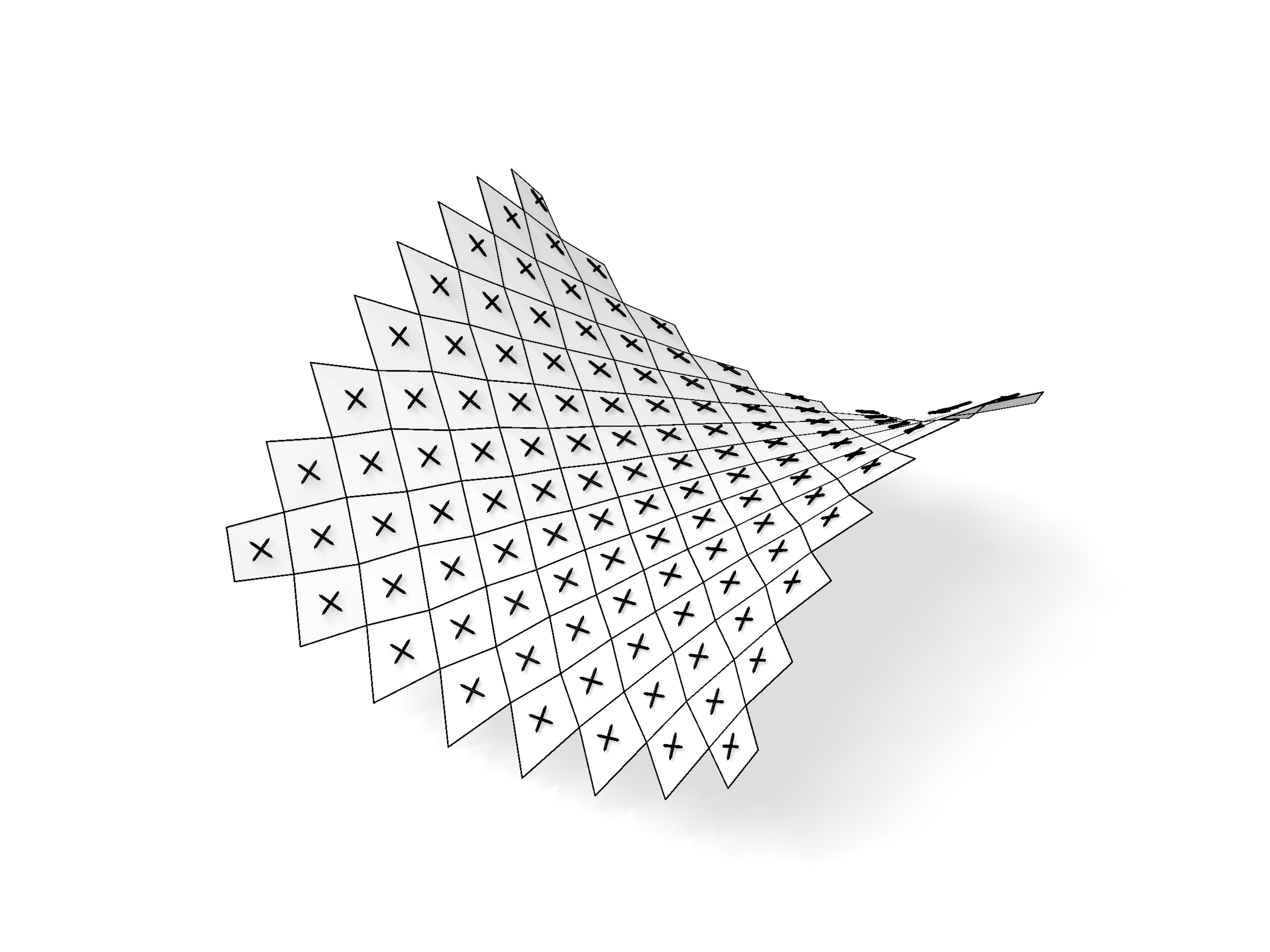}
	\end{subfigure}
	\begin{subfigure}[c]{0.49\linewidth}
		\includegraphics[trim = 10cm 0cm 5.5cm 10cm, clip,
		width=\linewidth]{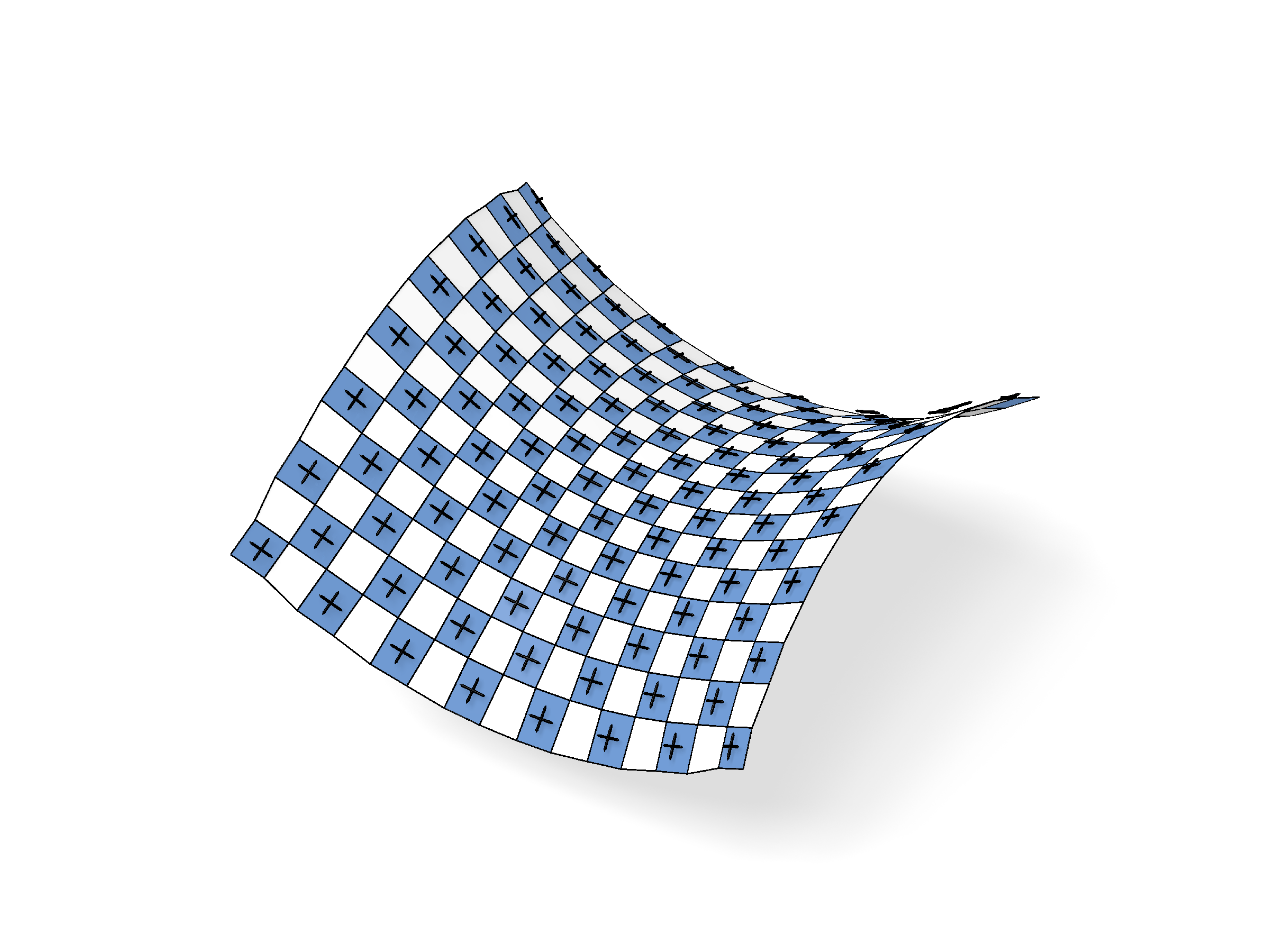}
	\end{subfigure}
	\vskip -5mm
	\caption{\textit{Left:} A control net of a principal checkerboard pattern and the eigenvectors of the shape operator. \textit{Right:} The checkerboard pattern of the same net. We see that the first order faces are aligned with the eigenvectors of the shape operator as stated by Corollary \ref{Kor:PrincipalDirections}.}
	\label{Fig:ConsistencyOfPrincDir}
\end{figure}

\begin{figure}[ht]
	\begin{overpic}[width=0.45\linewidth, trim = 2cm 2cm 5cm 5cm, clip ]{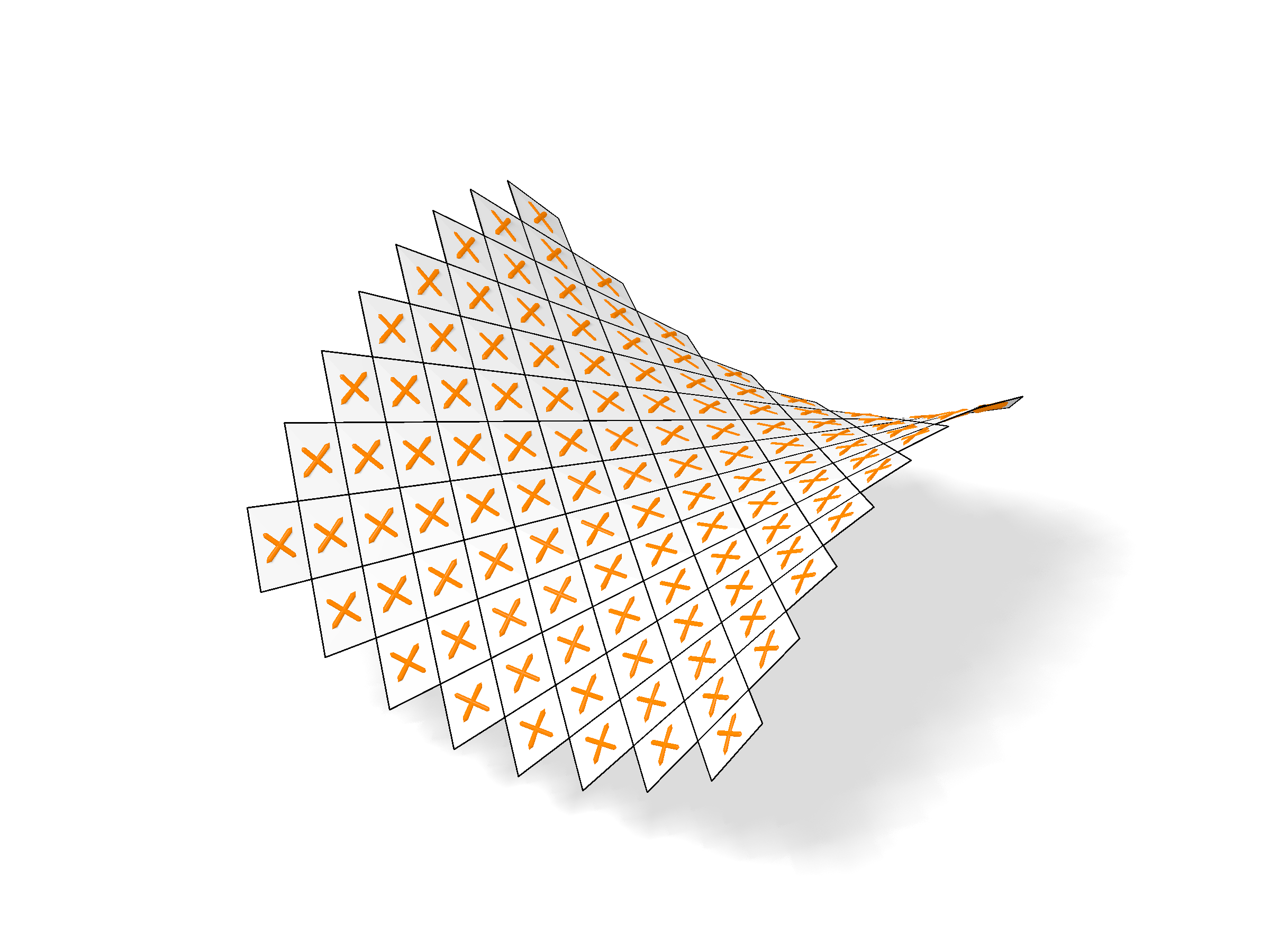}
	\end{overpic}
	\begin{overpic}[width=0.45\linewidth, trim = 2cm 2cm 5cm 5cm, clip]{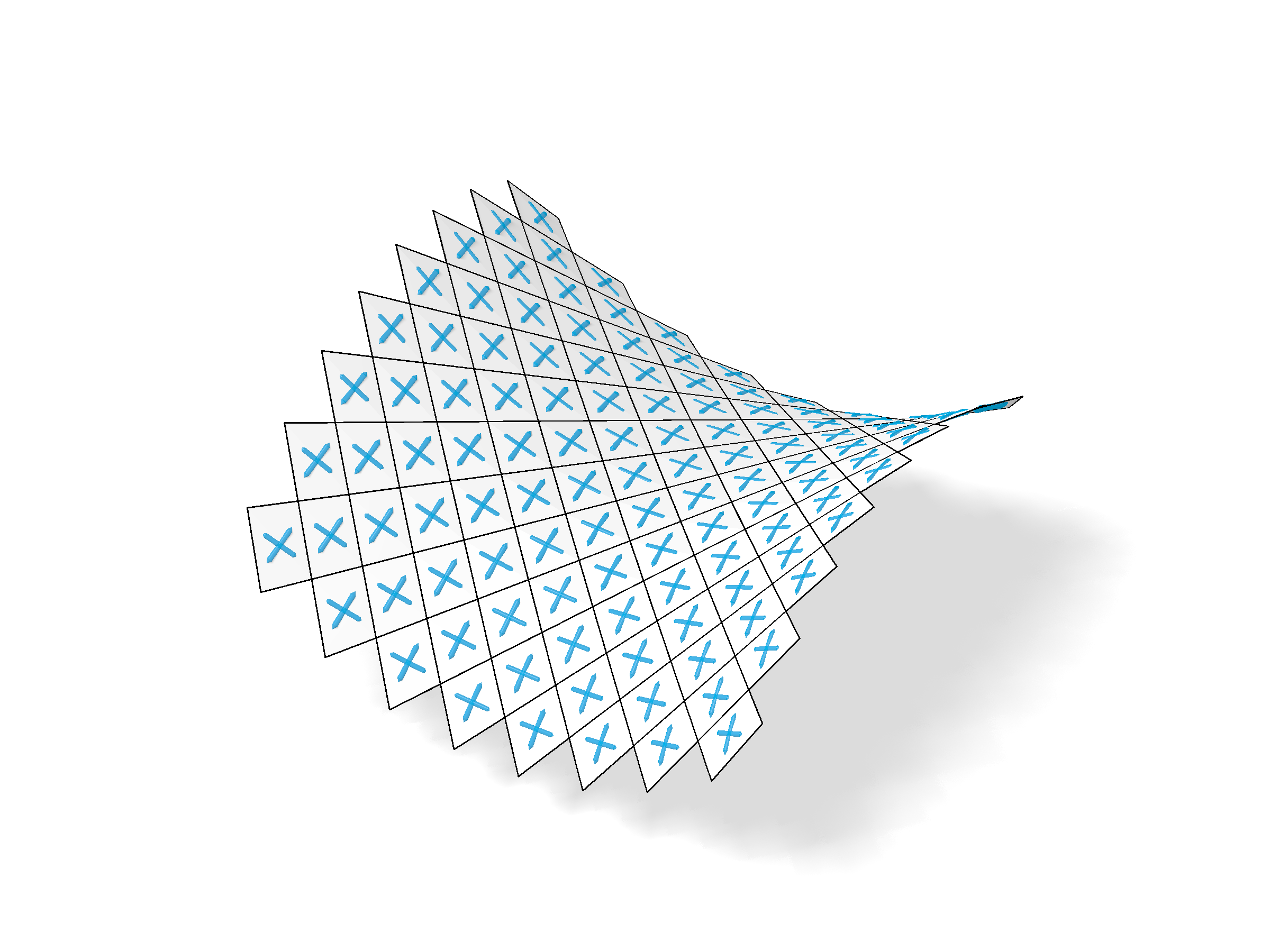}
	\end{overpic}
	\vspace*{-1.5em}
	\caption{If the net $\net$ samples a smooth parametric surface $\phi$, the underlying checkerboard pattern can be used to compute the discrete principal curvature directions \textit{(left)} which are visually not distinguishable from the analytically computed directions \textit{(right)}.}
	\label{Fig:ApproxOfPrincDir}
\end{figure}

In this section we define a discrete version of the shape operator connecting nets to their Gauß images. We find that the properties of the shape operator for conjugate or principal nets are consistent with the smooth theory, see Figure \ref{Fig:ConsistencyOfPrincDir}. Moreover, the discrete shape operator provides a way to numerically approximate smooth principal curvature directions, compare Figure \ref{Fig:ApproxOfPrincDir}. We start by defining the Gauß image of a net.

\begin{defi}
	\label{Defi:GeneralizedSurfAreaAndNormVector}
	Let $\net$ be a net. Then
	\begin{align*}
		n = \frac{(f_{1} - f_{\bar{1}})\times (f_2 - f_{\bar{2}})}{\| (f_{1} - f_{\bar{1}})\times (f_2 - f_{\bar{2}}) \|}
	\end{align*}
	is a net with vertices on the unit sphere $\S^2$. We call $n$ the \emph{Gauß image} or \emph{vertex normals} of $\net$. Additionally, for the face $\quadf_f = (f,f_1,f_{12},f_2)$ we define the \emph{face normal} $N$ by
	\begin{align}\label{Eq:FaceNormal}
		N = \frac{(f_{12}-f)\times(f_2-f_1)}{\|(f_{12}-f)\times(f_2-f_1)\|}.
	\end{align}
	The \emph{generalized surface area} of $\quadf_f$ is the surface area of the orthogonal projection of $\quadf_f$ in direction of $N$,
	\begin{align}
		\area(\quadf) = \det(\discdifu{f},\discdifv{f},N).
	\end{align}
\end{defi}

\begin{bem}
	For planar quadrilaterals without self-intersections the generalized surface area is the same as the surface area. The face normal $N$ is a normal vector to $\blackface[\net]$ and for a planar face $\quadf$ it coincides with a normal vector to $\quadf$. The vertex normal $n$ at $\net$ is also the face normal of the corresponding second order face $\whiteface[\net]$ in the sense of formula \eqref{Eq:FaceNormal}.
\end{bem}


%

Having defined a Gauß image $n$ for a net $\net$, we can relate the discrete derivatives $(\discdifu{\net}, \discdifv{\net})$ and $(\discdifu{n}, \discdifv{n})$ with the help of the corresponding checkerboard patterns $\cbp[\net]$ and $\cbp[n]$. The idea is to define the shape operator as the linear mapping $(\discdifu{\net}, \discdifv{\net}) \mapsto (\discdifu{n}, \discdifv{n})$. However, we face the problem that $(\discdifu{\net}, \discdifv{\net})$ and $(\discdifu{n}, \discdifv{n})$ not necessarily span the same two-dimensional subspace. This is overcome by projecting in direction of $N$, leading to the following definition:

\begin{defi}\label{Defi:ShapeOperator}
	Let $\net$ be a net, let $\gaussimage[\net]$ be its Gauß image and let $P_N$ be the orthogonal projection along the corresponding face Gauß image $N$. We define $\shape$ as the function on $\Z^2$ that maps $(k,l)$ to a linear operator in the space spanned by $(\discdifu{f},\discdifv{f})$ such that
	\begin{align*}
		\shape (\discdifu{\net},\discdifv{\net}) = P_N (\discdifu{n},\discdifv{n})
	\end{align*}
	where all entities are evaluated at a point $(k,l)\in\Z^2$. We call $\shape(k,l)$ the \emph{shape operator} of the face $\quadf_{f}(k,l)$. If no confusion can arise we drop the argument $(k,l)$. The eigenvalues of $\shape(k,l)$ are denoted by the symbols $\kappa_1$ and $\kappa_2$ and are called the \emph{principal curvatures}. The eigenvectors of $\shape(k,l)$ are the \emph{principal curvature directions}.
\end{defi}


\begin{figure}\label{Fig:GaussImage}
	\begin{overpic}[width=0.49\linewidth]{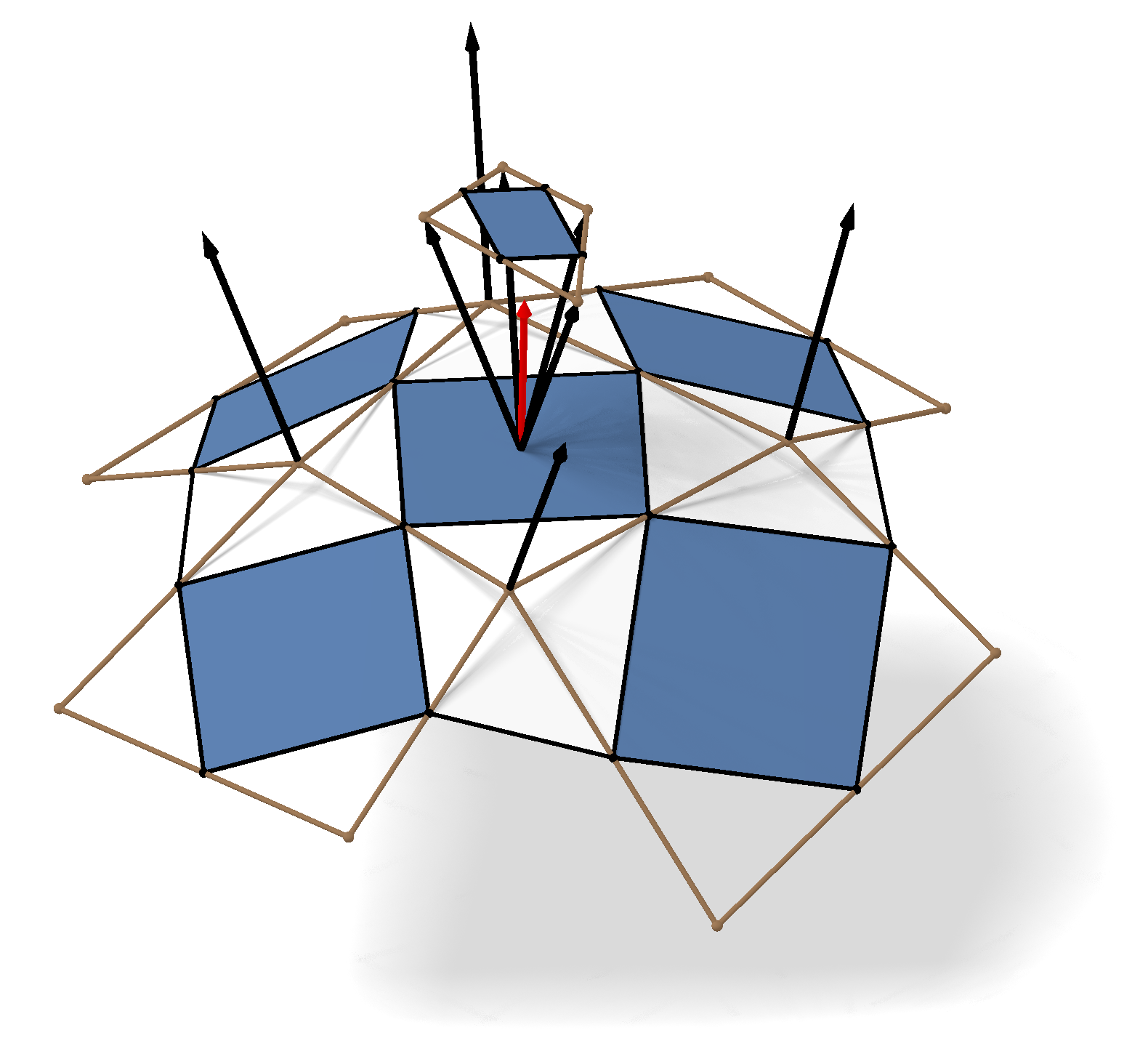}
	\end{overpic}
	\begin{overpic}[width=0.49\linewidth]{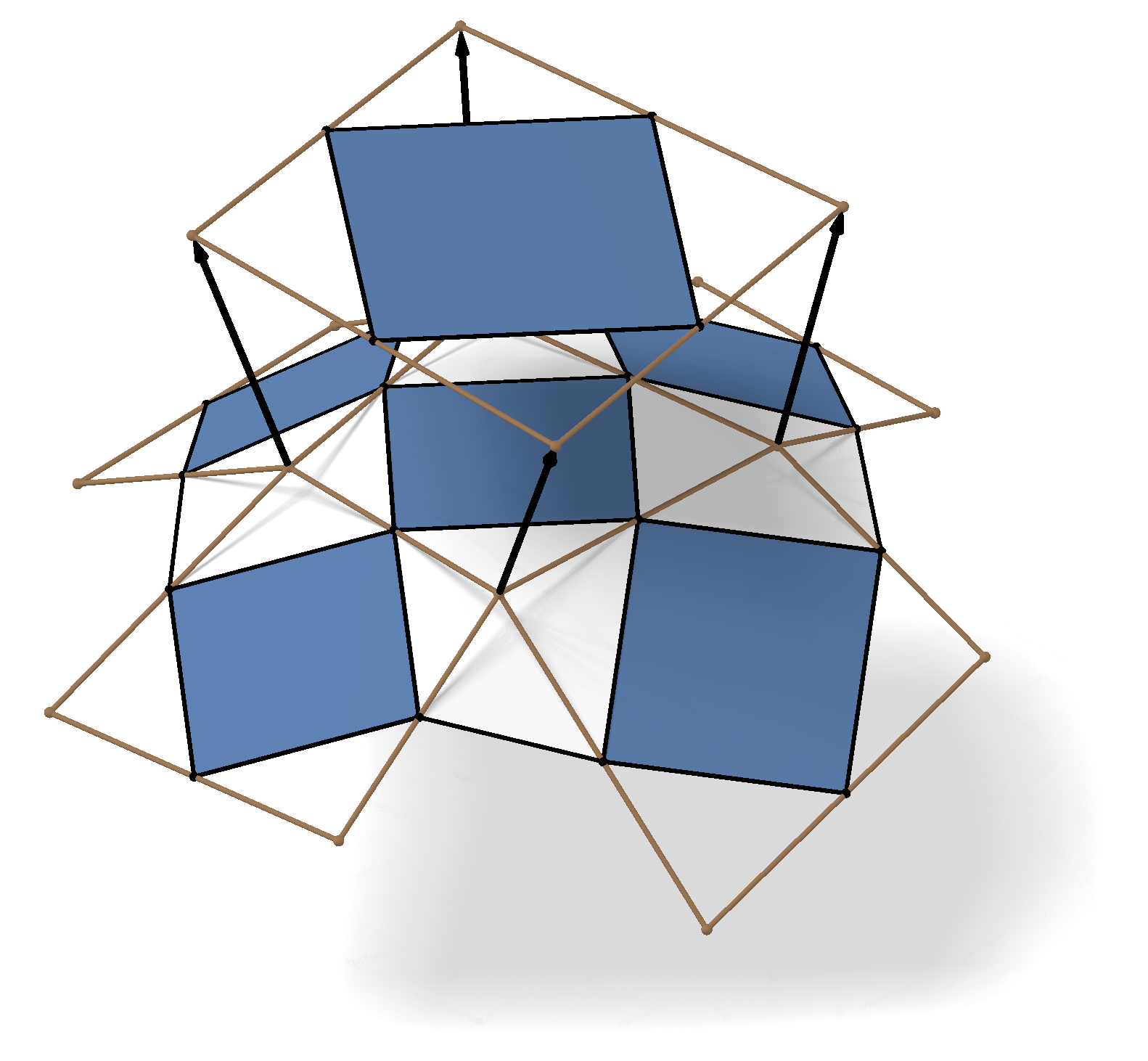}
	\end{overpic}
	\caption{\textit{Left:} A face $\quadf_f$ of $f$ and the corresponding Gauß image $\quadf_n$. The shape operator maps the first order face $\blackface[\net]$ to the first order face $\blackface[n]$ projected into the plane of $\blackface[\net]$. \textit{Right:} Face $\quadf$ and its offset $\quadf^t$.}
\end{figure}

For each face $\quadf$ we can define an offset face $\quadf^t$ by intersecting the plane parallel to $\blackface[\net]$ at distance $t$ with the lines spanned by the vertices of $\quadf$ and their corresponding vertex normals $n$. Similar to \cite{hoffmann2014discrete, MixedAreaHannesChristian, pottmann-2007-pm}, the area of $\quadf^t$ can be expressed by the \textit{Steiner formula}
\begin{align}\label{SteinerFormel}
	\area(\quadf^t) = (1 + t\cdot \tr(\shape) + t^2 \cdot \det(\shape)) \area(\quadf),
\end{align}
which can be shown by short algebraic manipulations.


\begin{lem}\label{Lem:ConjugateNets}
	For a conjugate checkerboard pattern the identities $\innerprod{\shape \discdif{u}{f}}{\discdif{v}{f}} = \innerprod{\discdifu{f}}{\shape \discdif{v}{f}} = 0$ hold. Thus the shape operator is symmetric.
\end{lem}
\begin{bew}
	For a conjugate checkerboard pattern $\cbp[f]$ the Gauß image $n$ is the normal vector of the corresponding second order face $\whiteface[f]$. Thus, it is orthogonal to all the edges that $\whiteface[f]$ shares with neighboring first order faces. As $\blackface[f]$ is a parallelogram, both $n_f$ and $(n_f)_{12}$ are orthogonal to the edge $\discdifv{f}$. We find that
	\begin{align*}
		0 = \innerprod{n - n_{12}}{\discdifv{f}} = 2 \innerprod{\discdifu{n}}{\discdifv{f}} = 2 \innerprod{P_N \discdifu{n}}{\discdifv{f}} = 2 \innerprod{\shape \discdifu{f}}{\discdifv{f}}.
	\end{align*}
	The same argument applies to $\innerprod{\shape \discdifv{f}}{\discdifu{f}}$. As $\discdif{u}{f},\discdif{v}{f}$ constitute a basis of the domain of the shape operator, the shape operator is symmetric.
\end{bew}

\begin{kor}\label{Kor:PrincipalDirections}
	For a principal checkerboard pattern 
	the edge vectors $(\discdifu{f},\discdifv{f})$ of $\blackface[\net]$ are eigenvectors of the shape operator.
\end{kor}
\begin{bew}
	This follows immediately from $\innerprod{\shape \discdif{u}{f}}{\discdif{v}{f}}= 0 = \innerprod{\discdif{u}{f}}{\discdif{v}{f}}.$ 
\end{bew}

As the partial derivatives can be observed in first order faces, so can the first fundamental form $\firstFundi$. By using the first order face $\blackface[\gaussimage]$ of the Gauß image and the corresponding derivatives $\discdifu{n}$ and $\discdifv{n}$ we can analogously define a second fundamental form.

\begin{defi}
	\label{Defi:FundamentalForms}
	Consider a net $\net$ and its Gauß image $n$. We define the \emph{first} and \emph{second fundamental forms} by letting
	\begin{align*}
		\firstFundi := \begin{pmatrix}
			\innerprod{\discdifu f}{\discdifu f} & \innerprod{\discdifu f}{\discdifv f} \\
			\innerprod{\discdifv f}{\discdifu f} & \innerprod{\discdifv f}{\discdifv f}
		\end{pmatrix}, \qquad
		\secondFundi := \begin{pmatrix}
		\innerprod{\discdifu f}{\discdifu n} & \innerprod{\discdifu f}{\discdifv n} \\
		\innerprod{\discdifv f}{\discdifu n} & \innerprod{\discdifv f}{\discdifv n}
		\end{pmatrix}.
	\end{align*}
\end{defi}

\begin{lem}
	A matrix representation $\shapeMat$ of the shape operator with respect to the basis $(\discdifu{f}, \discdifv{f})$ is given by
	\begin{align*}
		\shapeMat = \firstFundi^{-1} \secondFundi.
	\end{align*}
\end{lem}

\begin{bew}
	When using coordinates with respect to $(\discdifu{\net},\discdifv{\net})$ the inner product is $\innerprod{\cdot}{\cdot}$ is represented by the coordinate matrix $\firstFundi$. For any vector $v \in \operatorname{span}(\discdifu{f},\discdifv{f})$ we have $\innerprod{v}{\discdifu{n}} = \innerprod{v}{P_N \discdifu{n}}$ and likewise for $\discdifv{n}$. Thus the bilinear form $\innerprod{\cdot }{\shape \cdot }$ is represented by the coordinate matrix $\secondFundi$. For two vectors $\fw_1$ and $\fw_2$ with coordinates $w_1$ and $w_2$ we find that
	\begin{align*}
		w_1^T \secondFundi w_2 = \innerprod{\fw_1}{S \fw_2} = w_1^T \firstFundi \shapeMat w_2.
	\end{align*}
	It follows that $\secondFundi = \firstFundi \shapeMat$.
\end{bew}

\begin{bem}
	Due to Lemma \ref{Lem:ConjugateNets}, in a conjugate checkerboard pattern the second fundamental form is a diagonal matrix.
\end{bem}

In analogy to \cite{MixedAreaHannesChristian} and \cite{hoffmann2014discrete} the area defined in Definition \ref{Defi:GeneralizedSurfAreaAndNormVector} can be computed by a mixed area form. This motivates the following definition.

\begin{defiLemma}
	Let $A(\cdot,\cdot)$ be the \emph{mixed area form} defined by
	\begin{align}
		A(\quadf_f,\quadf_g) = \frac{1}{2}\Big( \det\big(\discdifu f, \discdifv g, N_f \big) + \det\big(\discdifu g, \discdifv f, N_f\big) \Big)
	\end{align}
	for two quadrilaterals with the same normal $N_f = N_g$. Then
	\begin{align*}
		\area(\quadf_f) = A(\quadf_f,\quadf_f)
	\end{align*}
	holds.
\end{defiLemma}

The mixed area form is closely related to the mean and Gaußian curvatures.

\begin{lem}
	For a net $\net$ and its Gauß image $n$ we define $\tilde{\quadf}_n$ as the orthogonal projection of $\quadf_n$ onto the supporting plane of $\blackface[f]$. The following identities hold
	\begin{align*}
		\det(\firstFundi) = A(\quadf_f, \quadf_f)^2,\quad 
		\kappa_1 \kappa_2 = \det(\shapeMat) = \frac{A(\tilde{\quadf}_n,\tilde{\quadf}_n)}{A(\quadf_f,\quadf_f)}, \\
		\frac{\kappa_1 + \kappa_2}{2} = \frac{1}{2} \tr(\shapeMat) = \frac{A(\quadf_f,\tilde{\quadf}_n)}{A(\quadf_f,\quadf_f)}.
	\end{align*}
\end{lem}
\begin{bew}
	These identities can be shown by algebraic manipulations in particular making use of the Lagrange identity $\innerprod{a\times b}{c \times d} = \det \left( {
	\innerprod{a}{c} \atop \innerprod{b}{c}} {\innerprod{a}{d} \atop \innerprod{b}{d}} \right)$.
\end{bew}

\begin{bem}\label{Remark:PCNormal_Forshadowing}
	Definition \ref{Defi:GeneralizedSurfAreaAndNormVector} requires that every normal vector lies exactly on the unit sphere. For principal nets one can relax this requirement and instead adapt the lengths of normal vectors, such that the first order face $\blackface[\gaussimage]$ of $\gaussimage$ is parallel to the first order face $\blackface[\net]$ of $\net$, as we will see in Section \ref{SubSection:PrincipalGauss}. This does not change principal directions and the Steiner formula \eqref{SteinerFormel} still holds. 
\end{bem}

\section{Möbius Transformations of Checkerboard Patterns}\label{Section:MoebiusTrafo}

This section discusses a way to apply Möbius transformations to orthogonal nets. This was originally introduced by \citet{techter2021discrete}, who showed that the orthogonality of a net is equivalent to the existence of a sphere congruence of orthogonally intersecting spheres. A Möbius transformation can then be applied to these spheres, and from the transformed congruence we can obtain the transformed orthogonal net. We show that the class of principal nets is invariant under such Möbius transformations. Moreover, the orthogonal sphere congruence allows us to embed principal nets in the projective model of Möbius geometry $\Mink$. This turns out to be a powerful tool for studying principal nets and gives rise to a non-Euclidean generalization of discrete principal nets.

\begin{figure}
\begin{overpic}[width=0.45\linewidth]{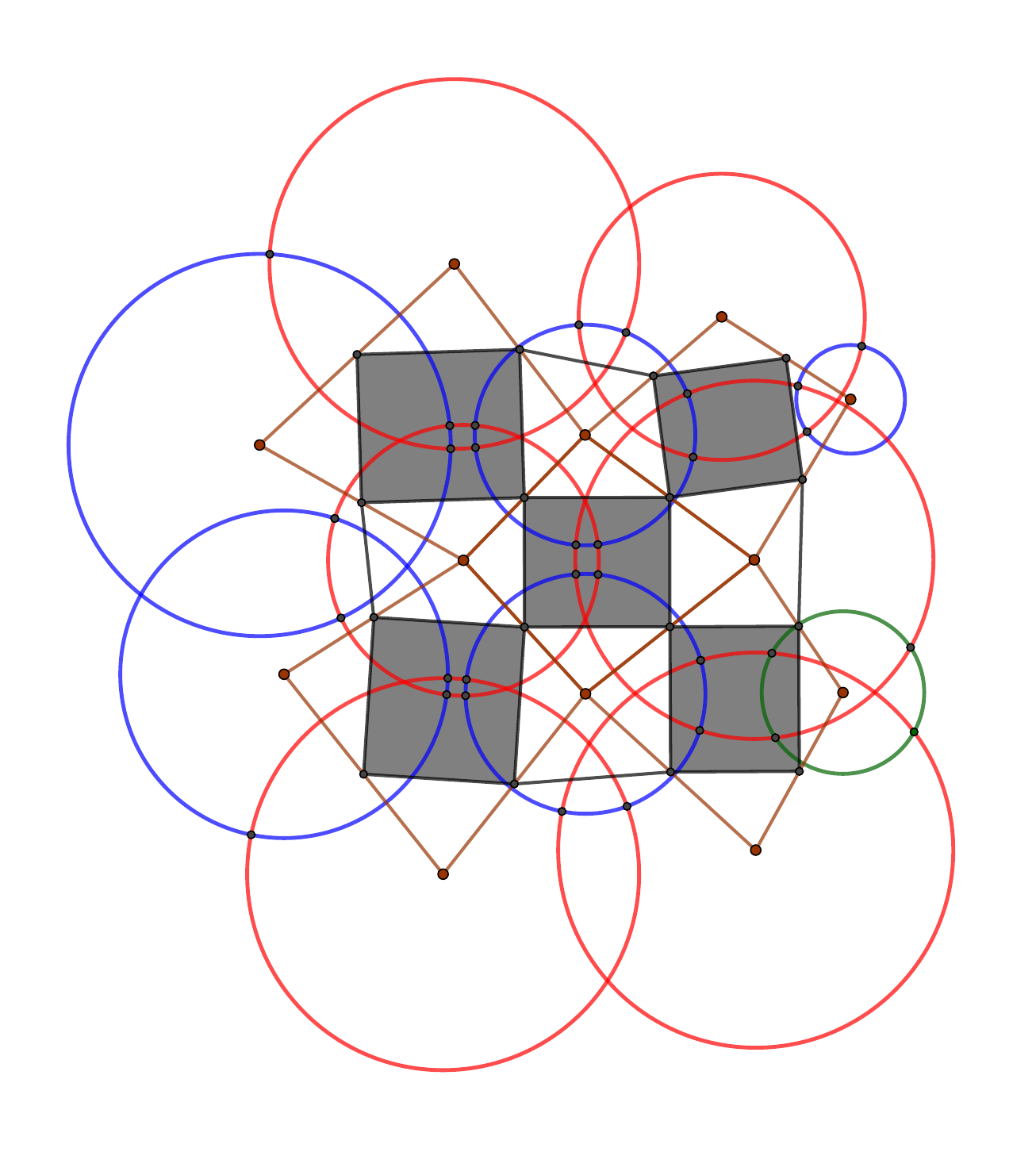}
\end{overpic}
\begin{overpic}[width=0.45\linewidth,tics=10]
	{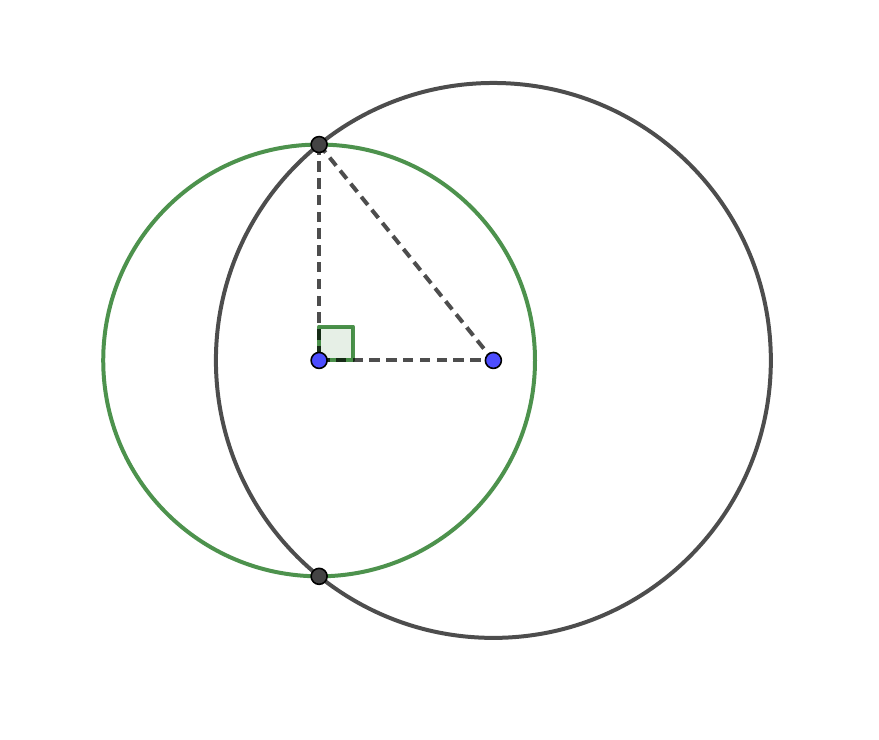}
	\put(30.5,39){\textcolor{geogebraBlue}{$c_2$}}
	\put(56,39){\textcolor{geogebraBlue}{$c_1$}}
	\put(42,52){$r_1$}
	\put(28.5,52){$|r_2|$}
	\put(11,61){\textcolor{geogebraGreen}{$\tilde{s}_2$}}
	\put(78,71){$s_1$}
\end{overpic}
\caption{\textit{Left}: The Möbius representation of a two-dimensional orthogonal net. Every red circle intersects every neighboring blue circle orthogonally and vice versa. The green circle represents a circle with negative squared radius. \textit{Right}: The circle $s_1 = (c_1,r_1^2)$ intersects the circle $s_2 = (c_2,r_2^2)$ with $r_2^2<0$ orthogonally if and only if it intersects the circle $\tilde{s}_2 = (c_2, -r_2^2)$ in diametrically opposite points.}
\label{Fig:ImagonaryOrthogonalIntersection}
\end{figure}
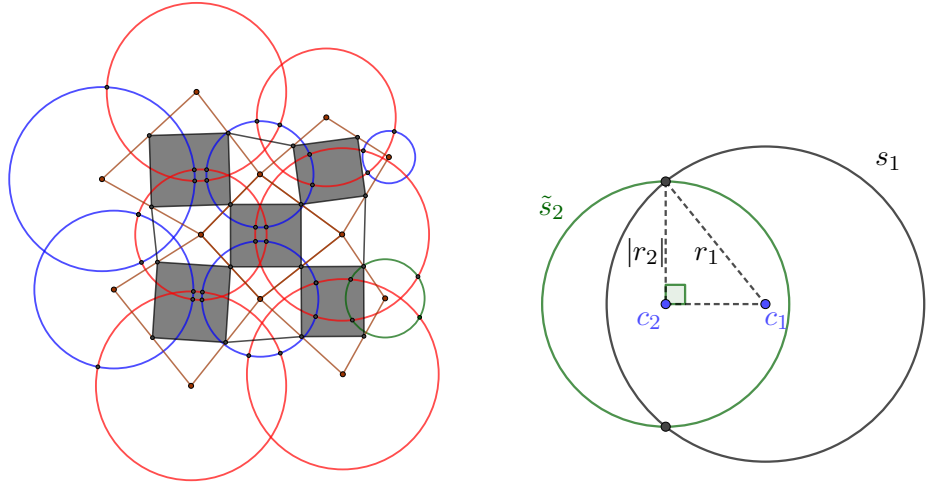

%
We write $s = (c,r^2)$ for a sphere $s$ with center $c\in\R^3$ and squared radius $r^2 \in \R$. Two spheres $s_1 = (c_1,r_1^2)$ and $s_2 = (c_2, r_2^2)$ intersect orthogonally, if and only if
\begin{align}\label{eq:Orthogonality}
	\innerprod{c_1-c_2}{c_1 -c_2} = r_1^2 + r_2^2.
\end{align}
Note that by definition this extends to spheres of negative squared radius. We can interpret this in the projective model of Möbius geometry by including the points inside the light cone as will be explained in more detail later in this section. Geometrically, the orthogonal intersection of a sphere with negative squared radius can be understood as illustrated by Figure \ref{Fig:ImagonaryOrthogonalIntersection}. This setup allows for the following lemma and definition.


\begin{defiLemma}
	Let $f$ be a net and $r^2:\Z^2 \to \R$. We call the function $s = (f,r^2)$ a sphere congruence and interpret it as a family of spheres with centers in $f$ and \textit{possibly complex} radius $r$. If and only if the checkerboard pattern of $\cbp[\net]$ is orthogonal, there exists a one-parameter family of sphere congruences $s=(f,r^2)$ such that neighboring spheres intersect orthogonally. We call such a sphere congruence the \emph{Möbius representation} $\moeb[\net]$ of $\net$ and $\cbp[f]$. If the checkerboard pattern associated to a sphere congruence $s$ is principal, we call $s$ a \emph{principal sphere congruence}. 
\end{defiLemma}

\begin{bew}
	Consider a quadrilateral $\quadf = (f, f_1, f_{12}, f_2)$. We fix the squared radius $r^2$ of $s=(f,r^2)$ at an initial point $(k,l) \in \Z^2$. This uniquely determines the radii $r_1$ and $r_2$ since
	\begin{align*}
		r_i^2 = \innerprod{f - f_i}{f - f_i} - r^2, \quad \text{for } i \in \{1,2\}.
	\end{align*}
	Now an easy computation shows that 
	\begin{align*}
		\innerprod{f_{12} - f_1}{f_{12} - f_1} - r_1^2 = \innerprod{f_{12} - f_2}{f_{12} - f_2} - r_2^2 \quad \iff \quad \innerprod{f - f_{12}}{f_1 - f_2} = 0.
	\end{align*}
	Hence, the radius $r_{12}$ is well defined if and only if the checkerboard pattern $\cbp[\net]$ is orthogonal. This process can be continued unambiguously, so every radius only depends on the choice of the initial radius.
\end{bew}

\begin{bem}
	If the domain of the net $\net$ is not simply connected, the orthogonal sphere congruences $\moeb[\net]$ do not exist in general. It is an interesting question for further research which properties might guarantee the existence of a Möbius representation for more complex topology or for combinatorial singularities.
\end{bem}


\begin{defiLemma}
	Let $\net$ be the control net of an orthogonal checkerboard pattern $\cbp[\net]$ and let $\moeb[\net]$ be its Möbius representation. The image of $\moeb[\net]$ under a Möbius transformation is again an orthogonal sphere congruence with a corresponding net $\net'$ and checkerboard pattern $\cbp[\net']$. We call $\net'$ resp $\cbp[\net']$ a \emph{Möbius transformation} of $\net$ resp $\cbp[\net]$.
\end{defiLemma}

The invariance of orthogonal nets under such Möbius transformations was observed in \citep{techter2021discrete}. The extension to complex spheres and the next theorem are contributions of this paper.

\begin{satz}\label{MoebiusInvariance}
	Principal checkerboard patterns are mapped to principal checkerboard patterns under Möbius transformations.
\end{satz}

\begin{bew}
	This follows directly from Theorem \ref{Theorem:ProjectiveModel} as will be explained later on in this section.
\end{bew}

\subsection{Projective model of Möbius geometry}

To prove Theorem \ref{MoebiusInvariance} we embed the Möbius representation $\moeb[\net]$ of a principal checkerboard pattern $\cbp[\net]$ into the projective model of Möbius geometry, see \cite{BobenkoBuch}. Let $\fe_1,\dots, \fe_5$ be the canonical basis vectors of the five-dimensional Minkowski space $\R^{4,1}$. It is equipped with the inner product
\begin{align*}
	\llangle \fe_i, \fe_j \rrangle = \begin{cases}
	1 & i = j \le 4, \\
	0 & i \neq j,\\
	-1 & i = j = 5.
	\end{cases}
\end{align*}
For $x \in \R^{4,1}\backslash \{0\}$ we write $[x]$ for its projective equivalence class, i.e.
\begin{align*}
	[x] = \left\{ y \in \R^{4,1} : y = \lambda x, \lambda \neq 0 \right\}.
\end{align*}
We write $\Mink$ for the space of these equivalence classes. Any sphere $s = (c,r^2)$ can be identified with a point of $\Pro\R^{4,1}$ by the mapping
\begin{align*}
	\iota(s) = \left[\left(c, \frac{1}{2}(|c|^2 - r^2 - 1), \frac{1}{2}(|c|^2 - r^2 + 1) \right) \right] \in \Pro\R^{4,1}.
\end{align*}
We view $c$ as a vector in $\R^{4,1}$ where the fourth and fifth components are zero and we define the vectors $\fe_0 := \frac{1}{2}(\fe_5 - \fe_4)$ and $\fe_{\infty} := \frac{1}{2}(\fe_4 + \fe_5)$. Then we can write
\begin{align*}
	\iota(s) = [c + \fe_0 + (|c|^2 - r^2)\fe_{\infty}].
\end{align*} 
Points can be seen as spheres with radius zero, so $\iota$ extends to points in $\R^3$. Observe that $\llangle \iota(s),\iota(s) \rrangle = r^2$. Thus the set of spheres with radius zero is identified with the light cone $\LightCone:= \{ x \in \Mink : \llangle x,x\rrangle = 0\}$. The points inside the light cone are those with $\llangle x, x\rrangle < 0$ and correspond to spheres with negative squared radius.

From a Möbius geometric point of view, planes in $\R^3$ are spheres with infinite radius and center at infinity. We write $\epsilon = (n,d)$ for the plane defined by the equation $\langle n, x \rangle = d$. The mapping $\iota$ can now be extended to spheres with infinite radius (i.e., planes) by
\begin{align}
	\iota(\epsilon) = [n +  0\cdot \fe_0 + 2d \fe_{\infty}].
\end{align}

The advantage of the projective model of Möbius geometry lies in the well known linearization of orthogonal intersection and Möbius transformations \cite{BobenkoBuch}.

\begin{satz}\label{Theorem:IntersectionMinkowski}
	Two spheres $s_1$ and $s_2$ in $\R^3$ with squared radii in $\R \cup \{\infty\}$ intersect orthogonally if and only if $\llangle \iota(s_1) , \iota(s_2) \rrangle = 0$. If one sphere has radius $0$, orthogonal intersection is equivalent to just intersection.	Möbius transformations in $\R^3$ canonically extended to spheres an planes are exactly the orthogonal transformations in $\Mink$.
\end{satz}


\begin{defi}\label{Defi:PCNetsMoebModel}
	Let $\nesh$ be a net $\Z^2 \to \Pro\R^{4,1}$. If adjacent vertices are orthogonal, i.e., $\llangle \nesh, \nesh_1\rrangle = \llangle \nesh, \nesh_2 \rrangle = 0$, and the corresponding checkerboard pattern $\cbp[\nesh]$ is conjugate, we call $\nesh$ a \emph{pseudo-principal net in} $\Pro\R^{4,1}$. In order to avoid confusion we will denote nets in $\Pro\R^{4,1}$ by $\nesh$, while we use $\net$ for nets in $\R^3$.
\end{defi}

Let $\net$ be a net with orthogonal checkerboard pattern and let $\moeb[\net]$ be a corresponding sphere congruence. Then $\iota \circ \moeb[\net]$ is a net $\Z^2 \to \Mink$, where the vertices are the images of $\moeb[\net]$ under $\iota$. 

\begin{satz}\label{Theorem:ProjectiveModel}
	If $\moeb[\net]$ is a principal sphere congruence in $\R^3$, then $\iota(\moeb[\net])$ is a pseudo-principal net in $\Pro\R^{4,1}$. If $\nesh$ is a pseudo-principal net in $\Pro\R^{4,1}$, then $\iota^{-1}(\nesh)$ is a principal sphere congruence in $\R^3$.
\end{satz}

\begin{figure}[h]
	\includegraphics[width=0.5\linewidth, trim = 0 -8cm 0 0, clip]{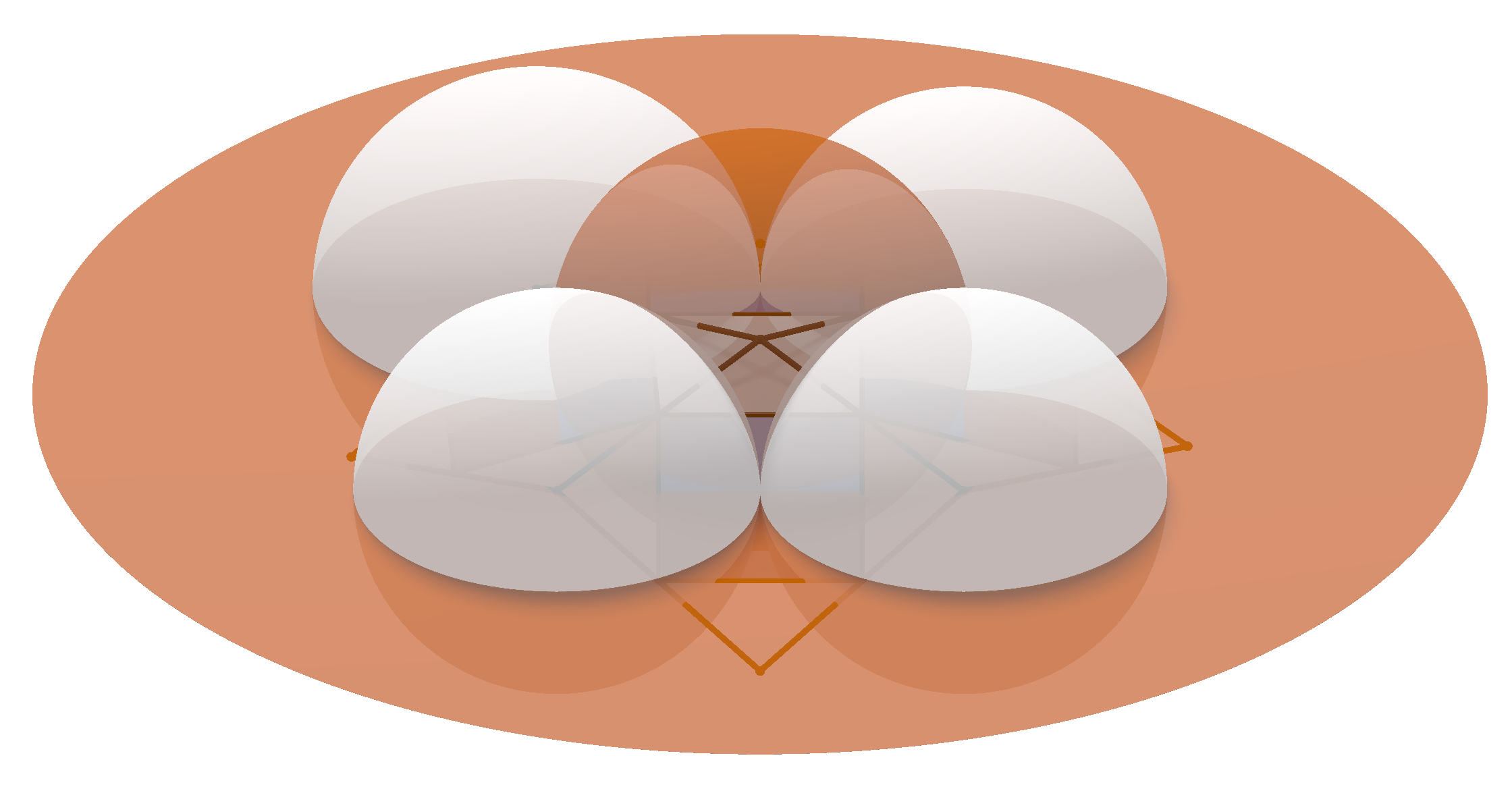}
	\hfill
	\includegraphics[width=0.35\linewidth]{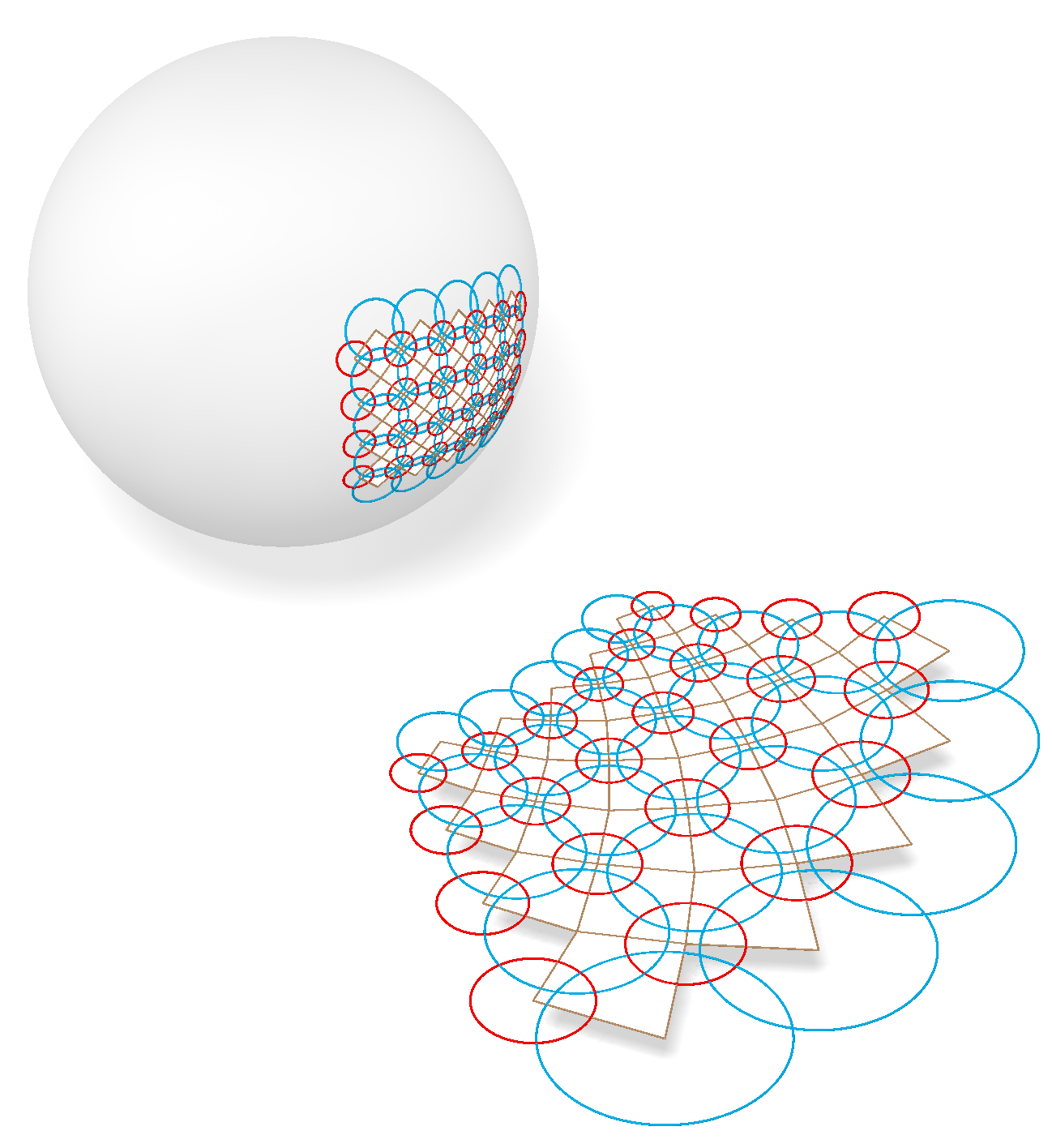}
	\vskip -0.5em 
	\caption{\textit{Left:} This figure illustrates why conjugacy of a checkerboard pattern is preserved under Möbius transformations. The four gray spheres are intersected orthogonally by the pencil spanned by the orange sphere and the orange plane. After applying a Möbius transformation the four gray spheres still intersect a pencil orthogonally which contains a plane. Hence the centers of the transformed spheres are still planar. \textit{Right:} The geometric description of the mapping $\iota$ in $\R^2$. A planar orthogonal circle pattern is stereographically projected onto the unit sphere. A new orthogonal net in space is obtained by the polar points of the circles on the sphere (not shown in the figure).}
	\label{Fig:BeweisMoebTrafo}
\end{figure}

\begin{bew}
	Orthogonality of adjacent vertices of a net in $\Pro\R^{4,1}$ is equivalent to the orthogonal intersection of adjacent spheres in $\R^3$.
	
	Let $\moeb[\net]$ be a principal sphere congruence in $\R^3$. The four spheres $s_{\bar{1}}, s_{\bar{2}}, s_1$ and $s_2$ all intersect both $s$ and the plane $\epsilon$ spanned by the centers $f_{\bar{1}}, f_{\bar{2}}, f_1$ orthogonally, compare Figure \ref{Fig:BeweisMoebTrafo}, left. Consequently, the four points $\iota(s_1), \iota(s_{\bar{1}}), \iota(s_2)$ and $\iota(s_{\bar{2}})$ all lie in the subspace $\iota(s)^\perp \cap \iota(\epsilon)^\perp$. Its dimension is two, since $\iota(\epsilon)$ and $\iota(s)$ are linearly independent. Hence, $\iota(s)$ is a pseudo-principal checkerboard pattern in $\Pro\R^{4,1}$.
	
	Now let $\nesh$ be a pseudo-principal net in $\Pro\R^{4,1}$ and let $U$ be the two-dimensional projective subspace that contains the four vertices $g_{1}$, $g_{\bar{1}}$, $g_2$ and $g_{\bar{2}}$. We denote by $U^\perp$ its orthogonal complement with respect to the Minkowski inner product $\llangle\cdot,\cdot \rrangle$. The space of all points in $\Pro\R^{4,1}$ that represent a plane in $\R^3$ is given by $\{\fe_{\infty}\}^\perp$. 
	Referring to the projective space $\Mink$ we have $\dim U^\perp = 1$ and $\dim \{\fe_{\infty}\}^\perp = 3$. It follows that $\dim( U^\perp \cap \{\fe_{\infty}\}^\perp) \ge 0$ and thus contains at least one point $\epsilon$. Since $\epsilon$ is a plane that intersects all points in $U$ orthogonally, we conclude that the centers of $g_{1}$, $g_{\bar{1}}$, $g_2$ and $g_{\bar{2}}$ all lie in $\epsilon$ and thus $\iota^{-1}(g)$ is a principal sphere congruence.
\end{bew}

Now Theorem \ref{MoebiusInvariance} easily follows from Theorem \ref{Theorem:ProjectiveModel}.

\begin{bew}[Proof of Theorem \ref{MoebiusInvariance}]
	As Möbius transformations in $\Mink$ are given by orthogonal transformations of $\R^{4,1}$, they preserve both orthogonality and $k$-dimensional subspaces. Thus pseudo-principal nets are mapped to pseudo-principal nets in $\Mink$ and by Theorem \ref{Theorem:ProjectiveModel} this translates to principal nets in $\R^3$ as well.
\end{bew}

\begin{bem}
	In classical differential geometry, a principal net $f$ can be characterized by the fact that its lift to the light cone $\hat{f} = f + \fe_0 + |f|^2\fe_{\infty}$ is a conjugate net. The mapping $\iota(\cdot)$ is a natural discretization of $f \mapsto \hat{f}$ as $\iota(s)$ converges to $\hat{f}$ if the radius of the sphere $s$ with center $f$ converges to zero. Like in the classical theory $\iota(s)$ is a conjugate net. However, $\iota(s)$ reveals even more structure, namely the orthogonality of spheres, that cannot be observed in the limit anymore.
\end{bem}

\begin{figure}[h]
	\begin{overpic}[width=0.32\linewidth, trim = 12cm 8cm 3cm 10cm, clip]{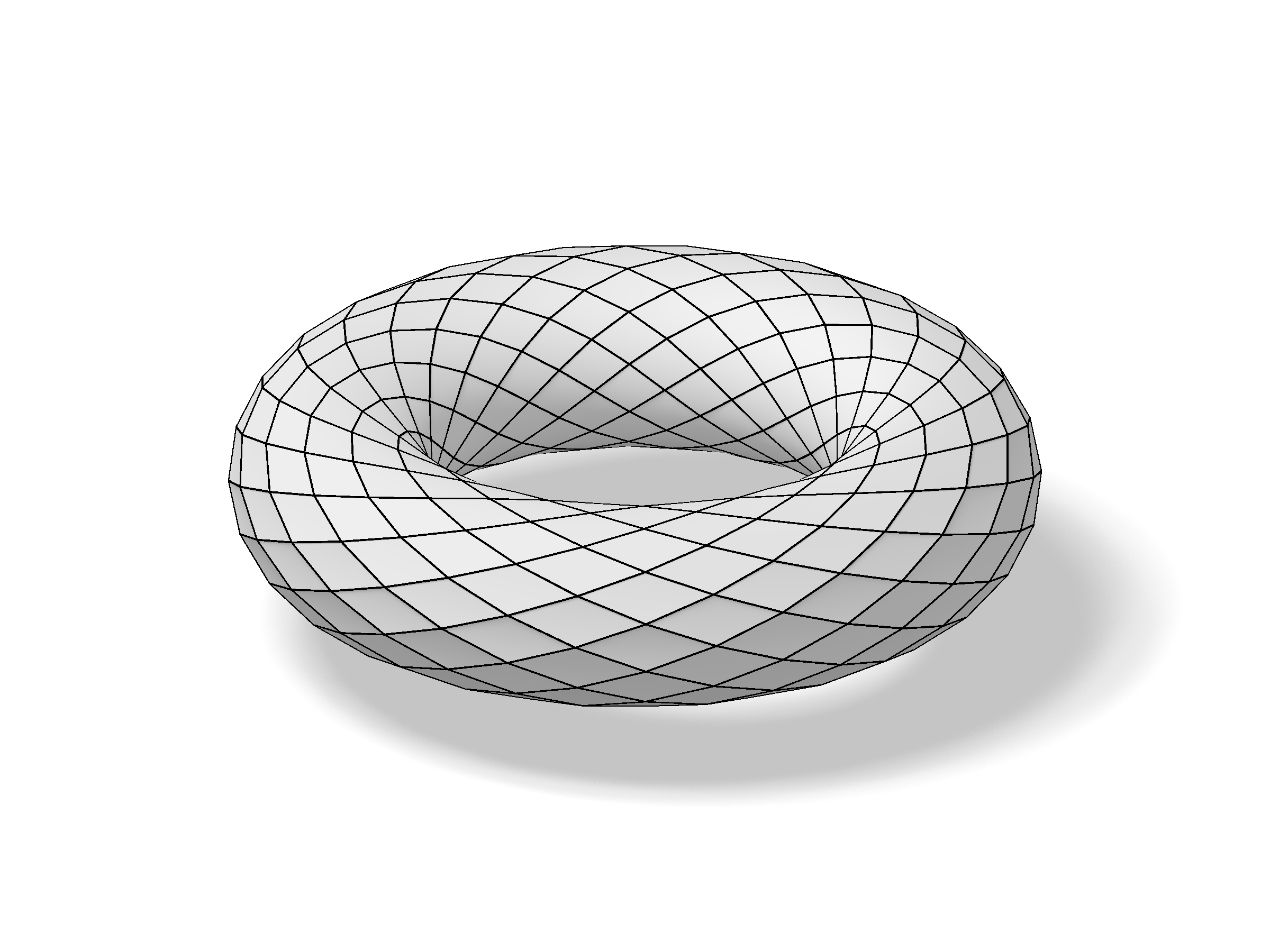}
	\end{overpic}
	\begin{overpic}[width=0.32\linewidth, trim = 12cm 8cm 3cm 10cm, clip]{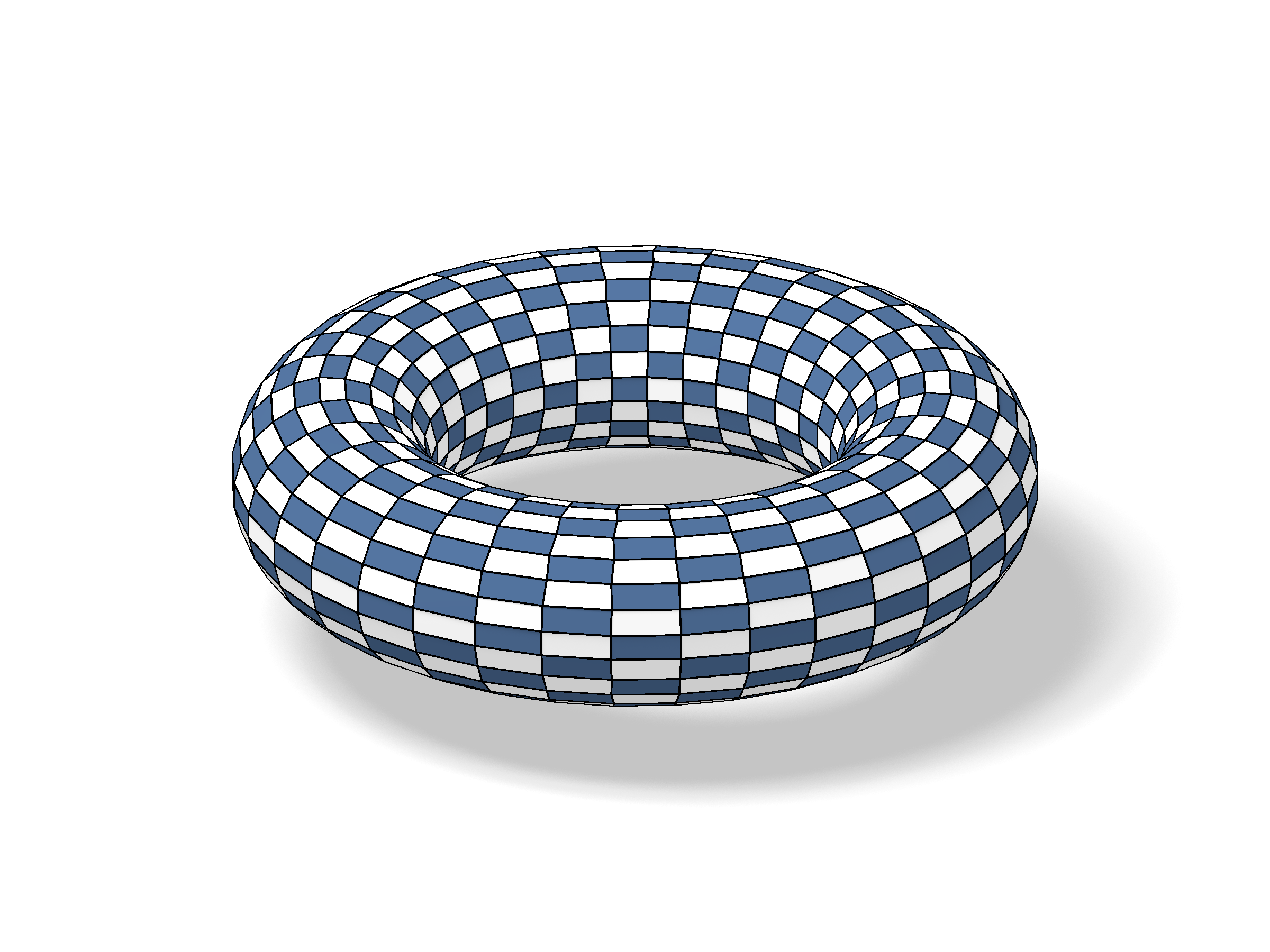}
	\end{overpic}
	\begin{overpic}[width=0.32\linewidth, trim = 12cm 8cm 3cm 10cm, clip]{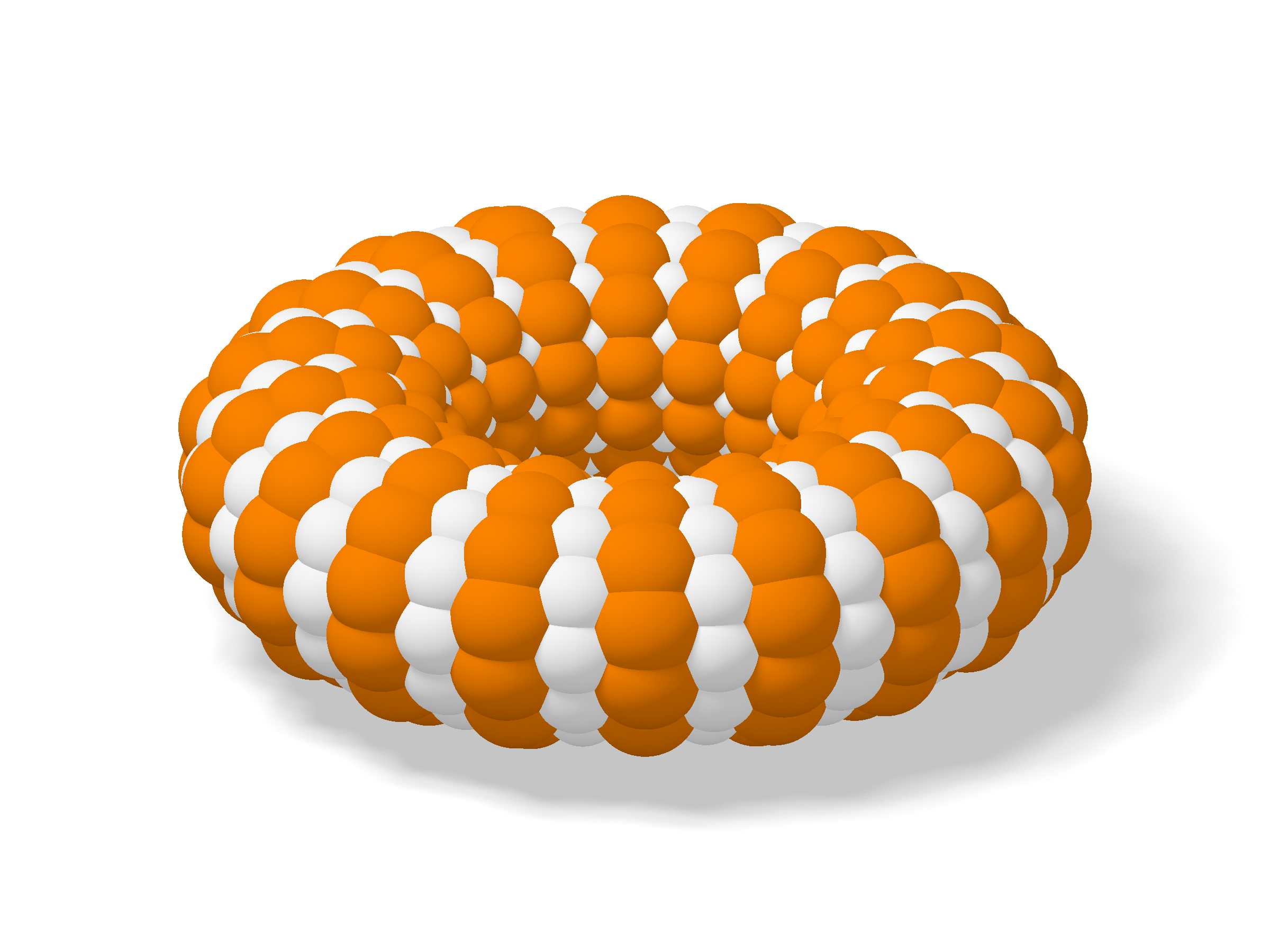}
	\end{overpic}
	\begin{overpic}[width=0.32\linewidth, trim = 8cm 4cm 13cm 5cm, clip]{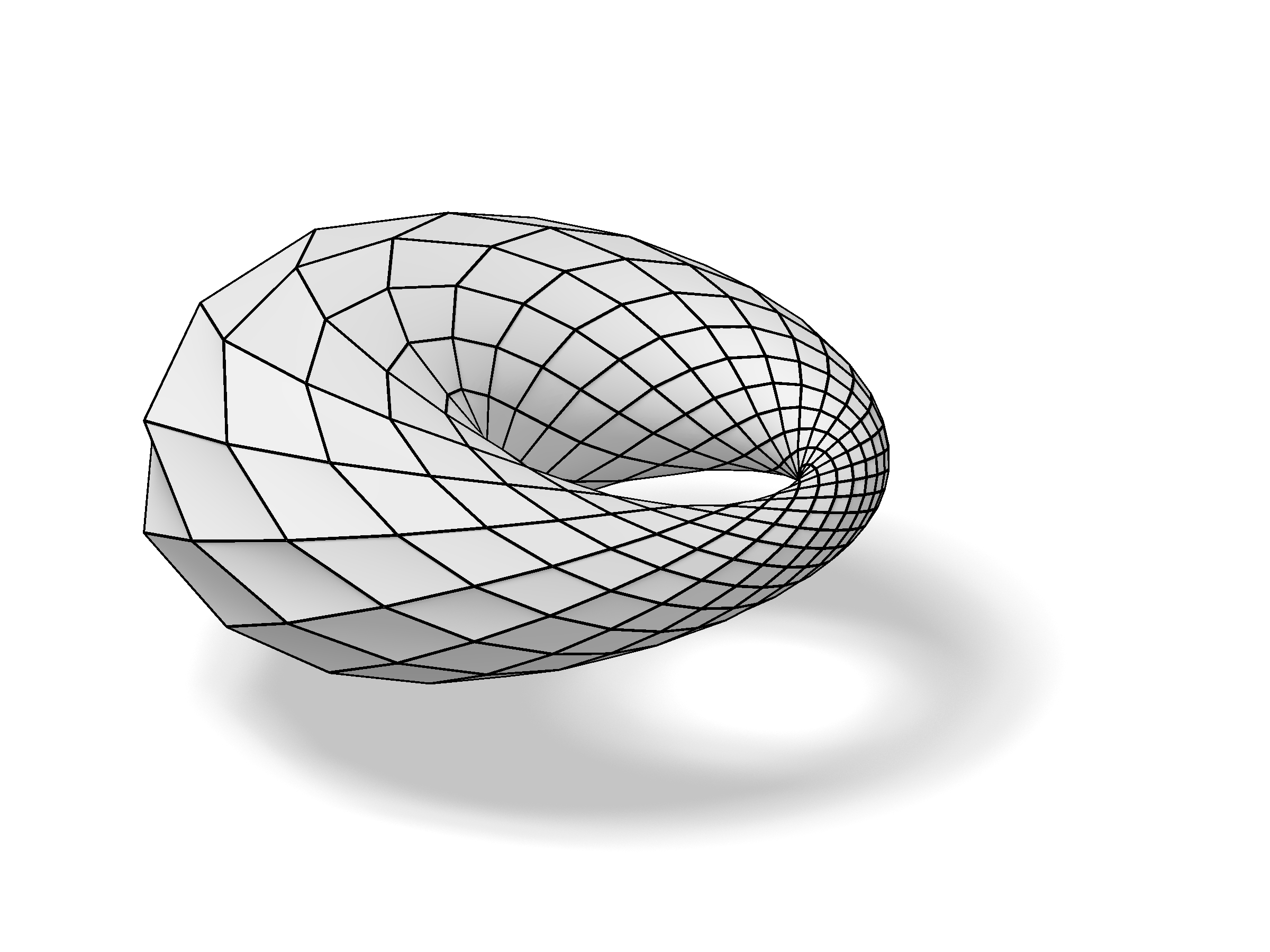}
	\end{overpic}
	\begin{overpic}[width=0.32\linewidth, trim = 8cm 4cm 13cm 5cm, clip]{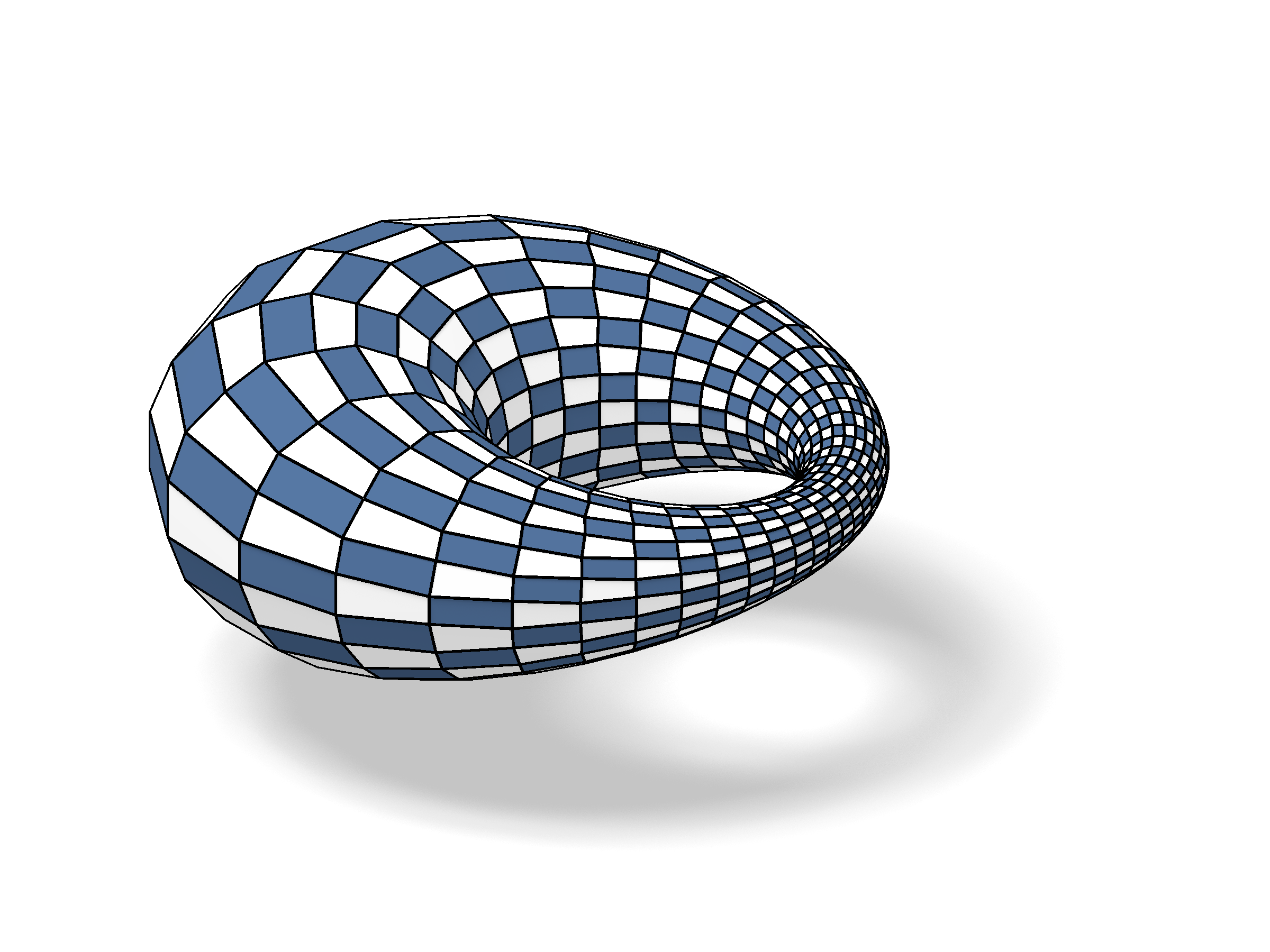}
	\end{overpic}
	\begin{overpic}[width=0.32\linewidth, trim = 0cm 4cm 13cm 5cm, clip]{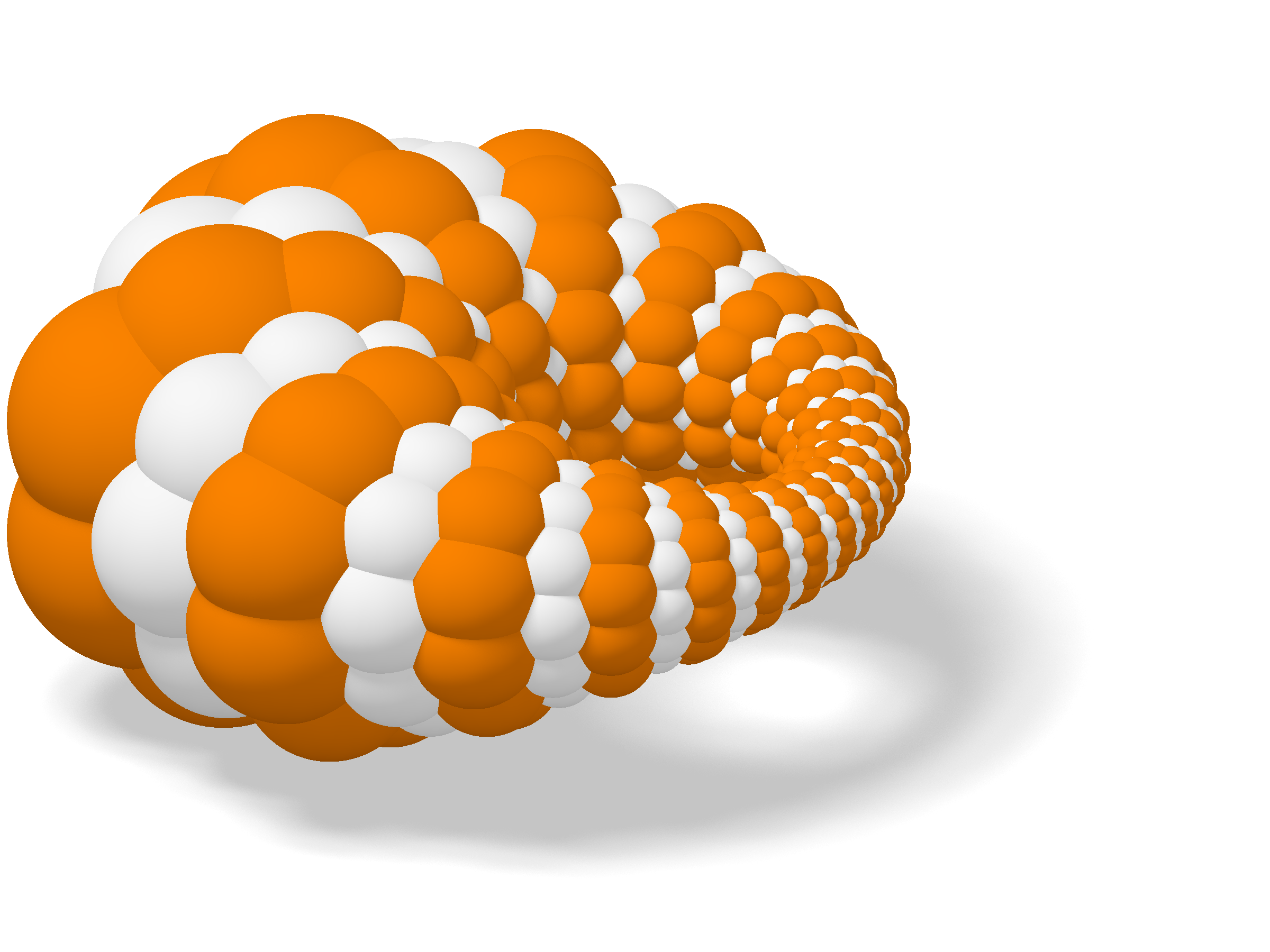}
	\end{overpic}
	\caption{In the first row we see, from left to right, the control net, the checkerboard pattern and a  Möbius representation of a principal net on the torus. The second row shows the image of the first row after a Möbius transformation is applied to the Möbius representation.}
	\label{Fig:MoebTrafoOfTorus}
\end{figure}


\subsection{A projective point of view}


It is enlightening to also study the embedding of the sphere congruence to $\Mink$ from a more geometric perspective.

The mapping $\iota$ can be seen as stereographically projecting a sphere $s$ to the unit sphere $\S^3$ and further mapping the image $s' \subseteq \S^3$ to its polar point $p=\iota(s)$ with respect to $\S^3$, compare Figure \ref{Fig:BeweisMoebTrafo}, right. The polar point $p$ is the apex point of the cone that touches $\S^3$ along $s'$. The polar point of any sphere $s_1'\subseteq \S^3$ that intersects $s'$ orthogonally lies in the polar hyperplane of $p$ and is thus conjugate to $p$. Hence, the diagonals of the quadrilateral $\big(\iota(s), \iota(s_1), \iota(s_{12}), \iota(s_2)\big)$ are not only orthogonal but conjugate with respect to $\S^3$.

The projective approach also gives meaning to the vertices of $\net$ in the projective model. They are the images of $\iota(s)$ under the central projection $\R^4 \to \R^3$ through the north pole of $\S^3$. 

\begin{bem}
	The unique sphere with center $\iota(s)$ that intersects $\S^3$ orthogonally, intersects $\S^3$ along $s'$. Hence, the vertices $\iota(\moeb[\net])$ define a unique sphere congruence $\mathcal{S}$ of three-dimensional spheres, where every sphere intersects its neighbors and also $\S^3$ orthogonally. The stereographic projection $\S^3 \to \R^3$ can be extended to a Möbius transformation $\zeta: \Mink \to \Mink$. The spheres of $\moeb[\net]$ can be directly obtained from the spheres of $\mathcal{S}$ by applying $\zeta$ and then intersecting the image with $\R^3$.
\end{bem}

\begin{bem}
	This geometric approach further allows us to generalize orthogonal sphere congruences to non-Euclidean geometry by replacing the stereographic projection from $\S^3$ to $\R^3$ by a central projection $\psi: \S^3 \to \R^3$. A sphere congruence on $\S^3$ conjugate with respect to $\S^3$ gets mapped to a congruence of non-Euclidean spheres. These non-Euclidean spheres intersect in directions conjugate with respect to $\psi(\S^3)^*$ the contour quadric of $\psi(\S^3)$, compare Figure \ref{Fig:HyperbolicSphereCongruence} and Lemma \ref{Lemma:Conjugacy} in the Appendix.
\end{bem}


\begin{figure}[h]
	\begin{subfigure}[t]{0.48\linewidth}
		\begin{overpic}
			[width=\linewidth]{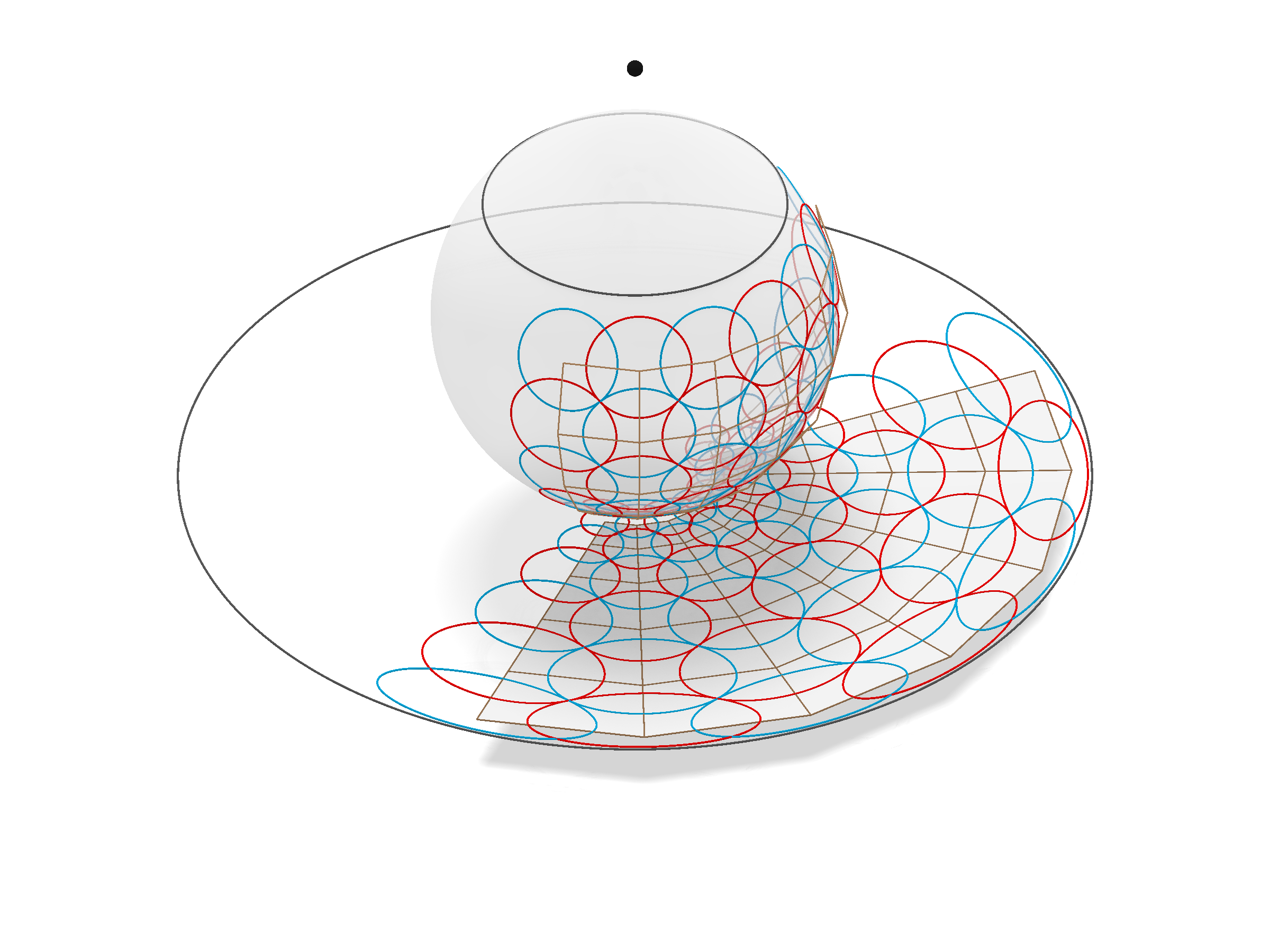}
			\put(22,33){$\mathbb{H}^2$}
		\end{overpic}
	\end{subfigure}
	\begin{subfigure}[t]{0.48\linewidth}
		\begin{overpic}
			[width=\linewidth]{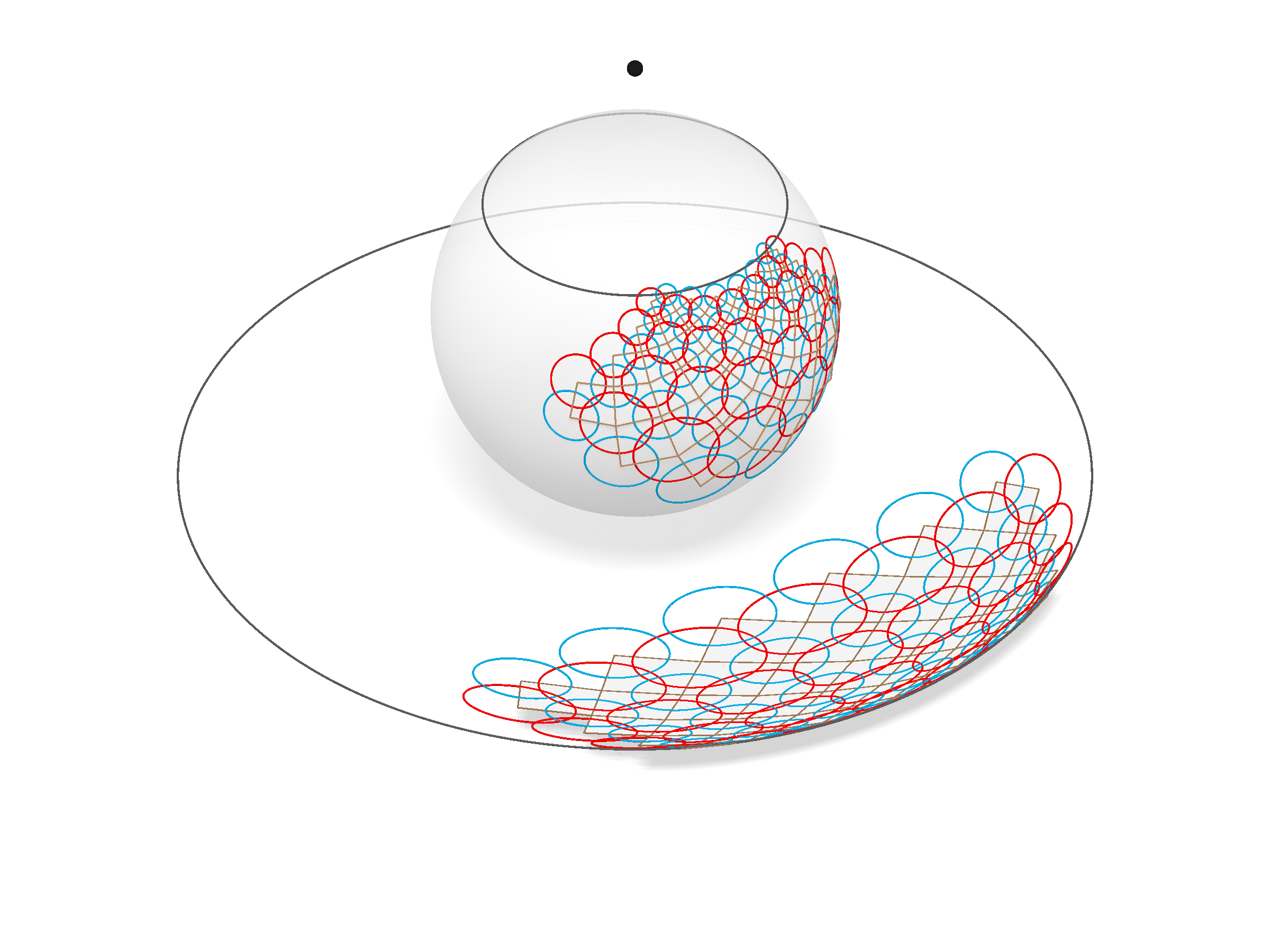}
			\put(22,33){$\mathbb{H}^2$}
		\end{overpic}
	\end{subfigure}
	\vskip -8mm
	\caption{Both images show an orthogonal net of circles on $\S^2$ which by a central projection is mapped to a net of conics. The net of conics is an h-orthogonal net of h-circles in the Cayley-Klein model of hyperbolic geometry.}
	\label{Fig:HyperbolicSphereCongruence}
\end{figure}


\subsection{A Gauß image for principal nets}\label{SubSection:PrincipalGauss}

As mentioned in Remark \ref{Remark:PCNormal_Forshadowing}, we can find an alternative definition of a Gauß image making use of the polarity properties of principal nets. This alternative is particularly interesting in connection with the minimal surfaces described in Section \ref{Section:MinimalSurfaces}.

\begin{defi}
	\label{Defi:NormalVector2}
	If $\net$ is a net with principal checkerboard pattern $\cbp[f]$, then $n$ is a \emph{principal Gauß image} of $\net$ and $\cbp[\net]$, if the edges of $\cbp[\net]$ are parallel to the edges of $\cbp[n]$ and every sphere of $\moeb[n]$ intersects the unit sphere orthogonally.
\end{defi}

The principal Gauß image $\net$ from Definition \ref{Defi:NormalVector2} of $\net$ can be seem as a parallel net of of $\net$ on the unit sphere. The parallelism can be observed in the corresponding checkerboard patterns $\cbp[\net]$ and $\cbp[n]$, while the connection to the unit sphere can be observed in the Möbius representation $\moeb[n]$. Instead of requiring vertices to lie exactly on the unit sphere, we require their corresponding spheres to intersect the unit sphere orthogonally. In the limit of spheres with radius zero the vertices lie exactly on the unit sphere. The principal Gauß image of a principal net $\net$ is determined up to the choice of one initial vertex along a prescribed line.

\begin{lem}
	Let $\net$ be a net with principal checkerboard pattern. Then there exists a one parameter family of principal Gauß images $n$ of $\net$ in the sense of Definition \ref{Defi:NormalVector2}. Diagonals of faces of $n$ are polar to one another with respect to the unit sphere.
\end{lem}

\begin{bew}
	To show the polarity of diagonals we consider a quadrilateral of four spheres $(s,s_1,s_{12}, s_2)$ with centers $(n,n_1,n_{12},n_2)$ that intersect $\S^2$ orthogonally. Additionally every sphere intersects its neighbors orthogonally. The centers of all spheres that intersect both $\S^2$ and $s$ orthogonally lie on a plane that contains the circle $\S^2 \cap s$. This plane is nothing but the polar plane of $n$. The same argument goes for $n_{12}$ and thus the diagonals $(n_1,n_2)$ and $(n,n_{12})$ lie on conjugate lines.
	
	From the polarity the uniqueness follows immediately. Let us fix one vertex $n(k,l)$ of $n$. Due to the parallelism of checkerboard patterns, we know the directions of diagonals emanating from $n(k,l)$. The four corresponding polar lines all lie in the polar plane of $n(k,l)$ and their intersection points determine the neighbors of $n(k,l)$. Thus, the initial vertex $n(k,l)$ corresponding to $f(k,l)$ needs to be chosen on a line orthogonal to $\whiteface[f]$. Note that polar lines are orthogonal and thus the parallelism is preserved. As polarity is a symmetric relation this process can be continued over the entire net.
	
	Now we can choose the initial radius of the sphere $s(k,l)$ at vertex $n(k,l)$ such that it intersects $\S^2$ orthogonally. The neighboring spheres of $s(k,l)$ have their centers in the plane of all centers of spheres that intersect $s(k,l)$ and $\S^2$ orthogonally. Hence all radii can be chosen such that the orthogonal intersection with both, all neighbors and the unit sphere is met.  Hence the so constructed net $n$ is indeed the principal Gauß image of $\net$.
\end{bew}

\begin{bem}
	We could also use the principal Gauß image in Definition \ref{Defi:ShapeOperator} of the shape operator. This only works for principal nets but it allows us to drop the orthogonal projection $P_\faceGaussimage$. Moreover, this approach fits the theory of minimal surfaces very well, as we will discuss in section \ref{Section:MinimalSurfaces}.
\end{bem}



\section{Koenigs nets}\label{Section:KoenigsNets}

In \cite{DoliwasKoenigsNets} Adam Doliwa defined discrete Koenigs nets as those conjugate nets where for every quadrilateral the six focal points lie on a common conic section, the so called conic of Koenigs. We apply the same definition to the checkerboard pattern $\cbp[\net]$ instead of the control net $\net$. This adaptation proves to be very useful as we can naturally dualize checkerboard patterns. Analogous to the smooth theory such a dual checkerboard pattern exists if and only if $\cbp[\net]$ is a Koenigs net. Even though the definition of Koenigs nets is based on checkerboard patterns we find that the class of Koenigs nets is invariant under projective transformations applied to the vertices of the corresponding control nets. Again in \cite{DoliwasKoenigsNets} Adam Doliwa defined discrete analogs of the so called Laplace invariants of a conjugate net. These projective invariants appear, in a slightly adapted way, in the checkerboard approach as well. However, it is only in this setting that Koenigs nets can be characterized as exactly those nets that have equal Laplace invariants analogously to the smooth theory.

\subsection{Characterization of Koenigs nets}

The discretization in both this paper and in \cite{DoliwasKoenigsNets} is based on the smooth characterization of Koenigs nets that can be found in \cite{lane1932projective}.



\begin{figure}[h]
\begin{overpic}[width=0.6\linewidth]
	{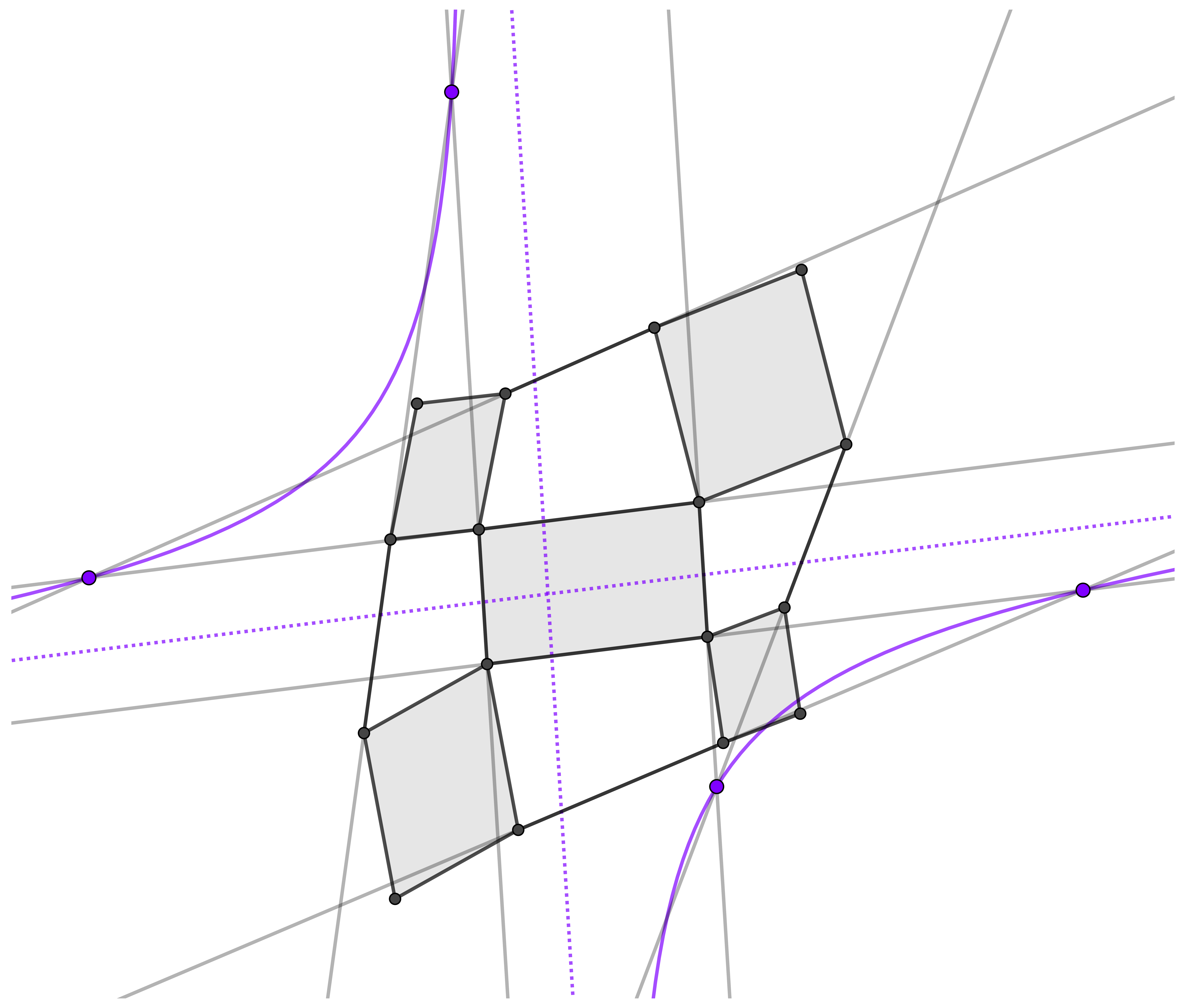}
	\put(41.5,30){$c$}
	\put(55,32){$c_1$}
	\put(60,40.5){$c_{12}$}
	\put(41.5,38){$c_2$}
	\put(68,34){$c_{11}$}
	\put(73,47){$c_{112}$}
	\put(25,24){$c_{\bar{1}}$}
	\put(27,41){$c_{\bar{1}2}$}
	\put(41,53){$c_{22}$}
	\put(55,58){$c_{122}$}
	\put(44,12){$c_{\bar{2}}$}
	\put(61,20){$c_{1\bar{2}}$}
	\put(79,40){\textcolor{geogebraPurple}{$p_1$}}
	\put(83,40){\rotatebox{5}{\textcolor{geogebraPurple}{$\to$}}}
	\put(46,65){\textcolor{geogebraPurple}{$p_2$}}
	\put(46,68){\rotatebox{93}{\textcolor{geogebraPurple}{$\to$}}}
	\put(91,32){\textcolor{geogebraPurple}{$p_3$}}
	\put( 8,33){\textcolor{geogebraPurple}{$p_4$}}
	\put(33,78){\textcolor{geogebraPurple}{$p_5$}}
	\put(62,16){\textcolor{geogebraPurple}{$p_6$}}				
\end{overpic}
\hfill
\begin{minipage}[b]{0.35\linewidth}
	\caption{Definition of a Koenigs net. The supporting lines of neighboring edges in the checkerboard pattern intersect in the six points $p_1,\dots,p_6$. If all of them lie on a common conic section the checkerboard pattern is Koenigs. The points $p_1$ and $p_2$ are always at infinity, here indicated by the dotted line, so the conic section is a hyperbola.
	\vspace*{3em}}
	\label{Fig:KoenigsNetDefi}
\end{minipage}
\end{figure}

\begin{defi}
	\label{Defi:KoenigsNet}
	Let $c$ be a conjugate checkerboard pattern. For the edge $(c,c_i)$ we denote the supporting line by $\suppl{c}{c_i}$. We call the checkerboard pattern $c$ a \emph{Koenigs checkerboard pattern} if for every first order face $(c, c_1, c_{12}, c_2)$ the six points
	\begin{align*}
	\begin{array}{lclcl}
	p_1 = \suppl{c}{c_1} \cap \suppl{c_2}{c_{12}}, & \qquad &
	p_2 = \suppl{c}{c_2} \cap \suppl{c_1}{c_{12}}, & \qquad & p_3 = \suppl{c}{c_1} \cap \suppl{c_{-2}}{c_{1-2}}, \\
	p_4 = \suppl{c_2}{c_{12}} \cap \suppl{c_{22}}{c_{122}}, & \qquad &
	p_5 = \suppl{c}{c_2} \cap \suppl{c_{\bar{1}}}{c_{\bar{1}2}}, & \qquad &
	p_6 = \suppl{c_1}{c_{12}} \cap \suppl{c_{11}}{c_{112}},
	\end{array}
	\end{align*}
	are all different and lie on a common conic section, see Figure \ref{Fig:KoenigsNetDefi}.
\end{defi}

\begin{bem}
	Since in Definition \ref{Defi:KoenigsNet} the points $p_1$ and $p_2$ are always at infinity, we know that the conic of Koenigs is always a hyperbola.
\end{bem}

The Koenigs nets defined in this way are a special instance of a class mentioned in \cite[p.\ 79]{BobenkoBuch}, where a result similar to Theorem \ref{Theorem:ClosedOneForm} is presented. Analogous to \cite{KoenigsNetsBobenkoSuris} the Koenigs property is equivalent to the existence of a closed multiplicative one-form on the edges of the checkerboard pattern. 

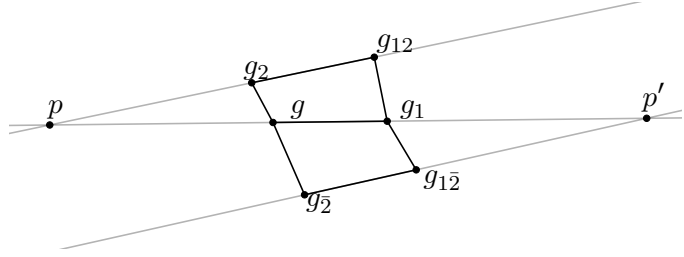
\begin{figure}
\begin{overpic}[width=0.7\linewidth]
	{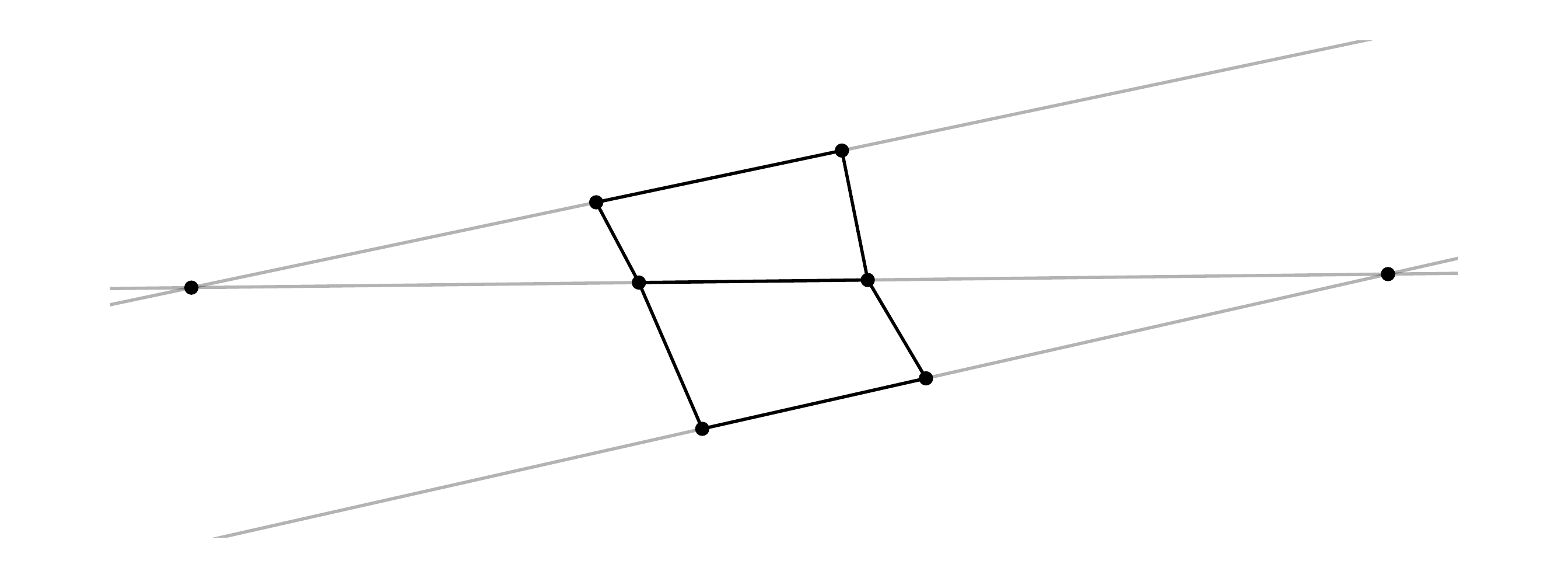}
	\put(43,20){$g$}
	\put(57,20){$g_1$}
	\put(88,21){$p'$}
	\put(12,20){$p$}
	\put(54,28.5){$g_{12}$}
	\put(37,25){$g_2$}
	\put(45,8){$g_{\bar{2}}$}
	\put(60,11){$g_{1\bar{2}}$}
\end{overpic}
\hfill
\begin{minipage}[b]{0.27\linewidth}
	\caption{The multiplicative one-form $q$ is defined on the edge $(g,g_1)$ as $q(g,g1)=\DV(g,g_1,p,p')$.\vspace*{3em}}
	\label{Fig:MultiplicativeOneForm}
\end{minipage}
\end{figure}

\begin{defi}[Multiplicative one-form]
	\label{Defi:MultiplicativeOneForm}
	Let $g$ be a net with planar quadrilaterals. Let further $p = \suppl{g}{g_1} \cap \suppl{g_2}{g_{12}}$ and $p' = \suppl{g}{g_1} \cap \suppl{g_{\bar{2}}}{g_{1\bar{2}}}$, see Figure \ref{Fig:MultiplicativeOneForm}. Then the multiplicative one-form $q$ along this edge $(g,g_1)$ is defined as the cross-ratio of the four points $g,g_1,p$ and $p'$
	\begin{align}
	q(g,g_1) := \frac{(g-p)}{(g_1 - p)} \frac{(g_1 - p')}{(g - p')} =: \DV(g,g_1,p,p').
	\end{align}
\end{defi}

%
%

\begin{satz}
	\label{Theorem:ClosedOneForm}
	Let $c$ be a conjugate checkerboard pattern such that the six points $p_1,\dots,p_6$ from Definition \ref{Defi:KoenigsNet} are all distinct. Let further $q$ be the multiplicative one-form from Definition \ref{Defi:MultiplicativeOneForm} defined on the edges of $\cbp$. Then $q$ is closed if and only if $c$ is Koenigs.
\end{satz}

\begin{bew}
	This theorem can be proven by introducing a projective coordinate system followed by lengthy computations that can be found in detail in the Appendix.
\end{bew}

The multiplicative one-form can also be formulated via the vertices of the control net as the following lemma shows.

\begin{lem}\label{Lemma:OneFormViaControlNet}
	Let $\whiteface[\net] = (c,c_1,c_{12},c_2)$ be a second order face of a checkerboard pattern $c$ with control net $f$. We choose the notation such that $c = \frac{1}{2}(f_{\bar{1}} + f)$, $c_1 = \frac{1}{2} (f + f_{\bar{2}})$, $c_{12} = \frac{1}{2} (f + f_1)$ and $c_2 = \frac{1}{2}(f + f_2)$, see Figure \ref{Fig:MultOneFormViaControlNet}. Let further $p = \suppl{c_2}{c_{12}} \cap \suppl{c}{c_1}$ and $P = \suppl{f_{\bar{1}}}{f_{\bar{2}}} \cap \suppl{f_1}{f_2}$. Then the multiplicative one-form $q$ is computed on the edge $(c_2, c_{12})$ as
	\begin{align*}
	q(c_2,c_{12}) = \frac{c_{12} - p}{c_2 - p} = \frac{f_1 - P}{f_2 - P}.
	\end{align*} 
	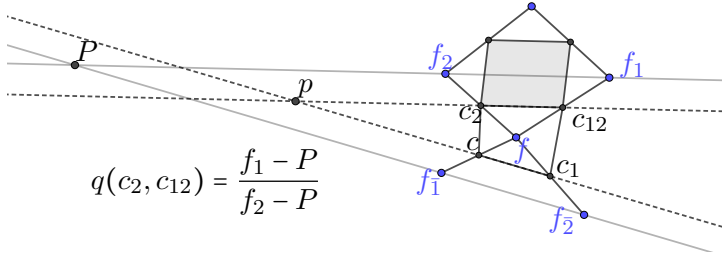
\begin{figure}[h]
	\begin{overpic}[width=0.7\linewidth]
		{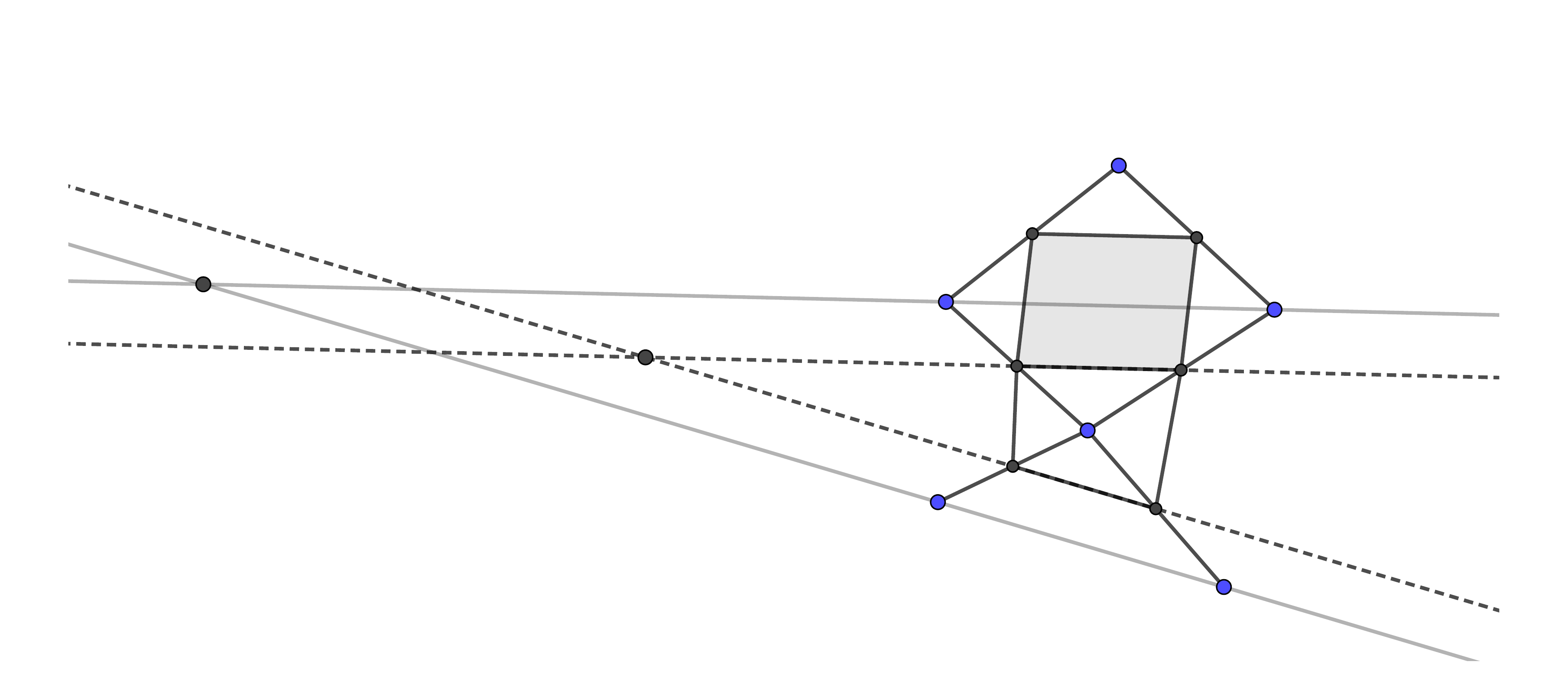}
		\put(69,14){\textcolor{geogebraBlue}{$f$}}
		\put(73.5,5){\textcolor{geogebraBlue}{$f_{\bar{2}}$}}
		\put(82.5,25){\textcolor{geogebraBlue}{$f_1$}}
		\put(58,26){\textcolor{geogebraBlue}{$f_2$}}
		\put(56.5,10){\textcolor{geogebraBlue}{$f_{\bar{1}}$}}
		\put(63,15){$c$}
		\put(74.5,12){$c_1$}
		\put(76.5,18){$c_{12}$}
		\put(62,19){$c_2$}
		\put(41.5,22.5){$p$}
		\put(13,26.5){$P$}
		\put(15,10){$\displaystyle q(c_2,c_{12}) = \frac{f_1 - P}{f_2 - P}$}
	\end{overpic}
	\hfill
	\begin{minipage}[b]{0.25\linewidth}
		\caption{The one-form $q$ defined on the edges of the checkerboard pattern can be expressed by the vertices of the control net.\vspace*{0.5em}}
		\label{Fig:MultOneFormViaControlNet}
	\end{minipage}
	\end{figure}
\end{lem}

\begin{bew}
	Since the quadrilateral $(f_{\bar{1}}, f_{\bar{2}}, f_{1}, f_{2})$ is the image of the quadrilateral $(c,c_1,c_{12},c_2)$ under the affine mapping $\alpha(x) = 2x - f$ the equality holds.
\end{bew}

This gives rise to the following lemma and definition.

\begin{figure}[h]
	\rightline{ \includegraphics[width=0.8\linewidth]{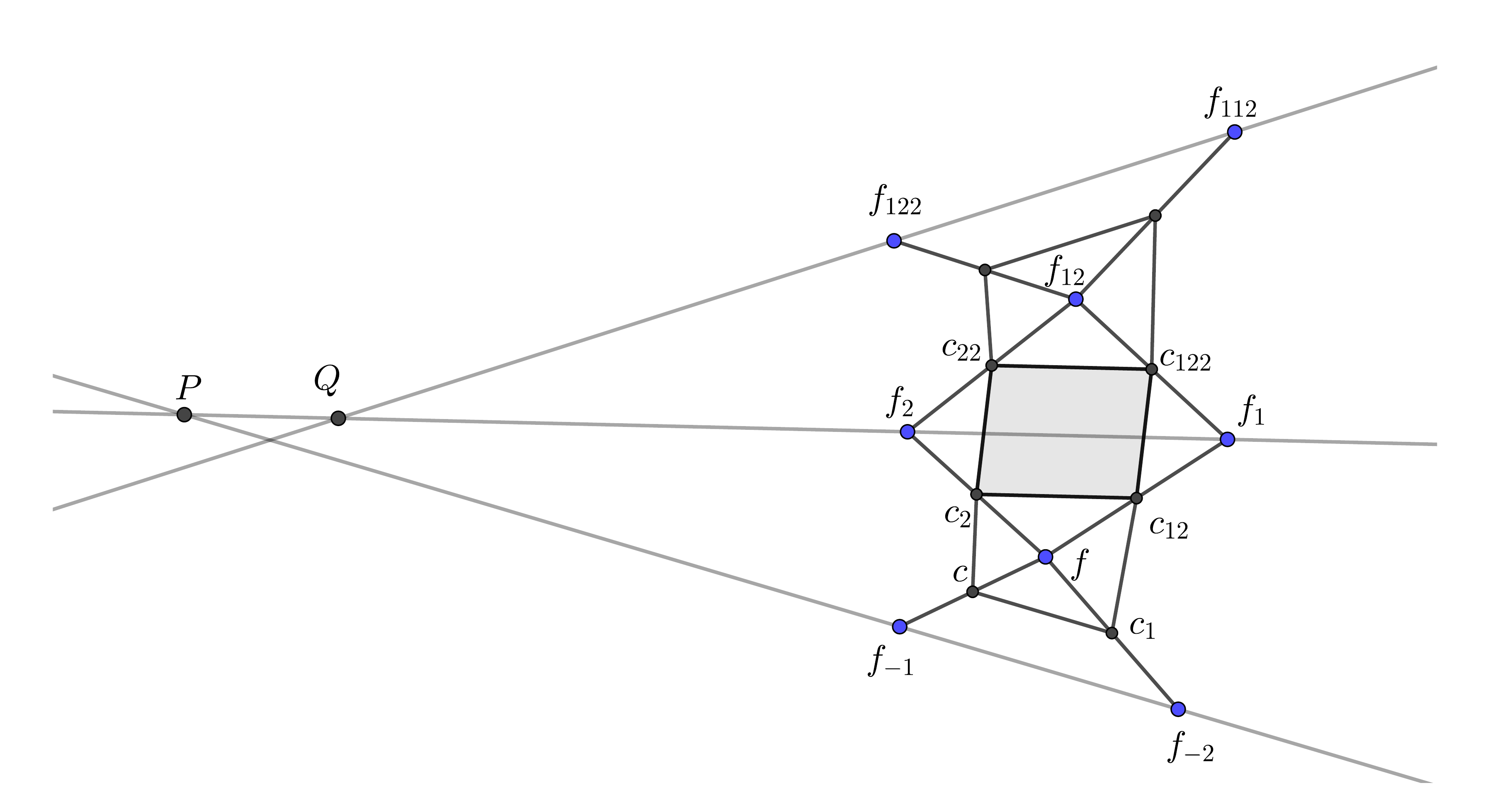}}
	\vskip -2cm
	\leftline{
	\begin{minipage}[c]{0.6\linewidth}
		\caption{The product $q(c_2,c_{12}) q(c_{122},c_{22})$ equals the cross ratio $\DV(f_1,f_2,P,Q)$ and is thus invariant under projective transformations applied to the control net. \vspace*{2em}}
		\label{Fig:MultOneFormViaControlNet2}
	\end{minipage}}
\end{figure}

\begin{defiLemma}\label{Defi:LaplaceInvariant}
	Consider the setting of Figure \ref{Fig:MultOneFormViaControlNet2} with the first order face $(c_2, c_{12}, c_{122}, c_{22})$. The product $q(c_2,c_{12})q(c_{122},c_{22})$ is a projective invariant of the control net. It is called \textit{Laplace invariant} and can be expressed via the control net by
	\begin{align}
	q(c_2,c_{12}) q(c_{122},c_{22}) = \frac{f_1 - P}{f_2 - P} \frac{f_2 - Q}{f_1 - Q} = \DV(f_1,f_2,P,Q).
	\end{align}
	To every face of the control net we can associate two Laplace invariants.
\end{defiLemma}

\begin{satz}\label{Theorem:LapInvariants}
	Let $\cbp[\net]$ be a conjugate checkerboard pattern with control net $\net$ such that the six points $p_1,\dots,p_6$ from Definition \ref{Defi:KoenigsNet} are all distinct. Then $\cbp[\net]$ is Koenigs if and only if the two Laplace invariants defined in each face of the control net are equal.
\end{satz}
\begin{bew}
	The two Laplace invariants of a face of the control net are equal if and only if the multiplicative one-form defined on the edges of the inscribed first order face is closed. Hence the statement follows from Theorem \ref{Theorem:ClosedOneForm}.
\end{bew}

\begin{bem}
	There are special cases where not all points $p_1,\dots,p_6$ are distinct, but the Laplace invariants are still equal. Those cases will turn out to be dualizable as well, so it makes sense to consider these nets to be Koenigs nets as well.
\end{bem}

\begin{bem}
	It is worth noticing that Theorem \ref{Theorem:LapInvariants} is independent of the choice of the control net. So if a checkerboard pattern is Koenigs every associated control net has equal Laplace invariants.
\end{bem}

\begin{kor}\label{Kor:ProjectiveInvariance}
	Koenigs checkerboard patterns are mapped to Koenigs checkerboard patterns under projective transformations applied to the vertices of the control nets.
\end{kor}

\begin{bew}
	The Laplace invariants are defined as cross ratios of vertices and intersection points of lines of the control net. Hence it is invariant under projective transformations and so the property of equal invariants is preserved as well.
\end{bew}

\begin{bem}
	Discrete Laplace invariants are defined for Koenigs nets in \cite[p.\ 5]{DoliwasKoenigsNets} in a similar fashion. The benefit of the checkerboard pattern approach is that now the Koenigs nets can be characterized as \grqq nets with equal invariants\grqq, compare Theorem \ref{Theorem:LapInvariants}, like one would expect coming from the smooth theory.
\end{bem}

\subsection{Dualization}

\begin{defi}
	\label{Defi:DualNet}
	Let $\cbp$ be a checkerboard pattern. We call $\cbp'$ a \emph{dual checkerboard pattern} of $\cbp$, if it is edgewise parallel and corresponding first order faces are similar but have reversed orientation. If such a dual checkerboard pattern $\cbp'$ exists, we call $\cbp$ \emph{dualizable}.
\end{defi}

In analogy to the smooth case we find that the dualizable checkerboard patterns are precisely the Koenigs checkerboard patterns. The following theorem holds.

\begin{figure}[h]
	\begin{overpic}[width=\linewidth]
		{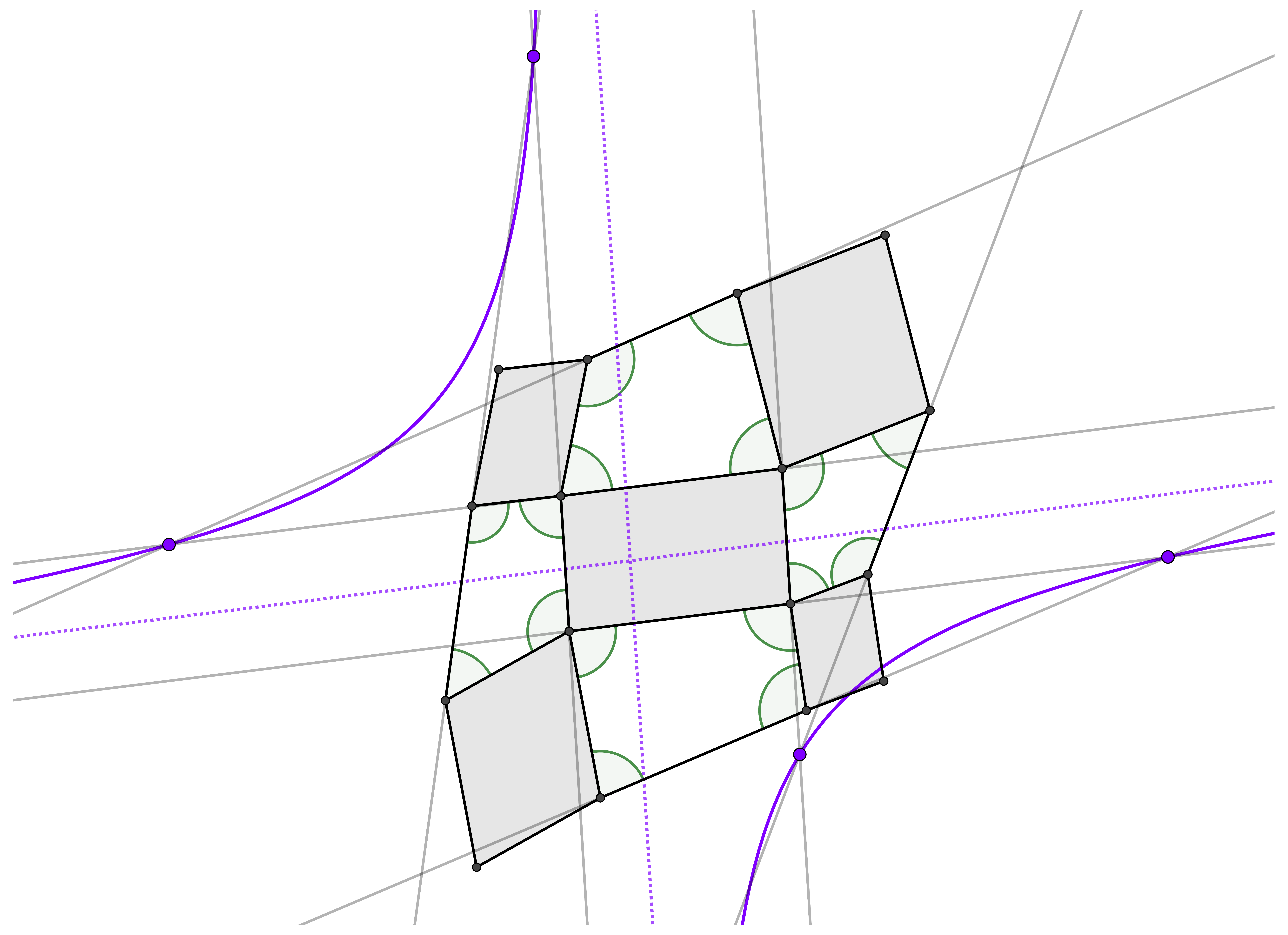}
		\put(79,34.5){\textcolor{geogebraPurple}{$p_1$}}
		\put(81.5,34.7){\rotatebox{5}{\textcolor{geogebraPurple}{$\to$}}}
		\put(48,65){\textcolor{geogebraPurple}{$p_2$}}
		\put(48.5,66.6){\rotatebox{93}{\textcolor{geogebraPurple}{$\to$}}}
		\put(90,31){\textcolor{geogebraPurple}{$p_3$}}
		\put(11,32){\textcolor{geogebraPurple}{$p_4$}}
		\put(37.5,68){\textcolor{geogebraPurple}{$p_5$}}
		\put(63,13){\textcolor{geogebraPurple}{$p_6$}}
		\put(44.8,22){\textcolor{geogebraGreen}{$\alpha_2$}}
		\put(46.7,12){\textcolor{geogebraGreen}{$\beta_2$}}
		\put(59.8,18){\textcolor{geogebraGreen}{$\gamma_2$}}
		\put(59.4,23.2){\textcolor{geogebraGreen}{$\delta_2$}}
		\put(61.2,26.9){\textcolor{geogebraGreen}{$\alpha_3$}}
		\put(65.1,28.5){\textcolor{geogebraGreen}{$\beta_3$}}
		\put(68.7,37.7){\textcolor{geogebraGreen}{$\gamma_3$}}
		\put(61.5,34.6){\textcolor{geogebraGreen}{$\delta_3$}}
		\put(57.5,37){\textcolor{geogebraGreen}{$\alpha_4$}}
		\put(55,46.7){\textcolor{geogebraGreen}{$\beta_4$}}
		\put(46,43){\textcolor{geogebraGreen}{$\gamma_4$}}
		\put(44.5,35){\textcolor{geogebraGreen}{$\delta_4$}}
		\put(41,32){\textcolor{geogebraGreen}{$\alpha_1$}}
		\put(36.7,31.5){\textcolor{geogebraGreen}{$\beta_1$}}
		\put(35.2,20){\textcolor{geogebraGreen}{$\gamma_1$}}
		\put(41.7,23.6){\textcolor{geogebraGreen}{$\delta_1$}}
		\put(44.5,29.5){$a$}
		\put(52.5,25){$b$}
		\put(46,16.5){$a_1$}
		\put(40,19){$b_1$}
		\put(62.5,22){$a_2$}
		\put(64,24.5){$b_2$}
		\put(64.8,39.2){$b_3$}
		\put(56,42.5){$a_3$}
		\put(45,39){$a_4$}
		\put(39,34.5){$b_4$}
	\end{overpic}
	\caption{The configuration of Theorem \ref{Theorem:KoenigsNets}.}
	\label{Figure:KoenigsNetAngles}
\end{figure}
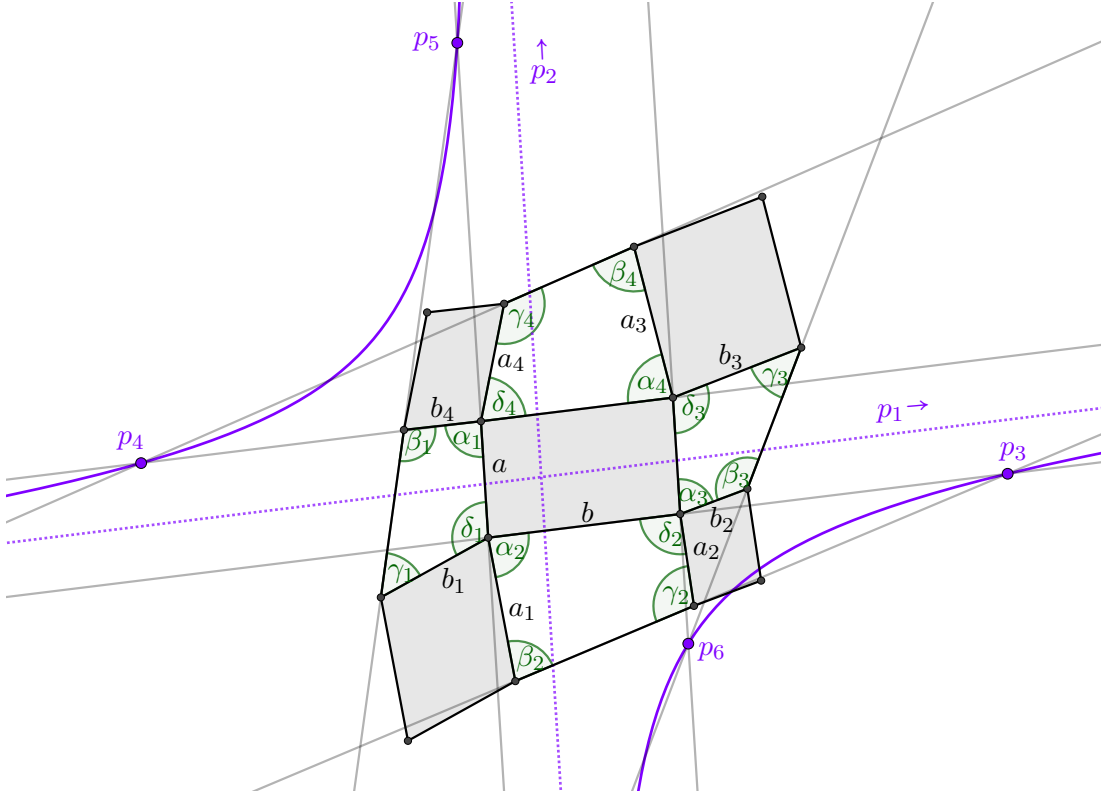

\begin{satz}
	\label{Theorem:KoenigsNets}
	Let $\cbp$ be a conjugate checkerboard pattern. We introduce the following local notation in the face patch of a given first order face, see Figure \ref{Figure:KoenigsNetAngles}:
	\begin{itemize}
		\item Let $a = \|\discdifv{f}\|$ and $b = \|\discdifu{f}\|$ be the edge lengths of the central first order face.
		\item We enumerate the exterior first order faces counterclockwise and denote their edge lengths with $a_i$ and $b_i$ accordingly.
		\item For every second order face in the patch we denote its interior angles by $\alpha_i,\beta_i,\gamma_i$ and $\delta_i$ in counterclockwise order.
	\end{itemize}
	Let $r_i = \frac{a_i}{b_i}$ be the ratio of edge lengths for each first order face. If no two of the six points $p_1,\dots, p_6$ from Definition \ref{Defi:KoenigsNet} are equal, the following conditions are equivalent:
	\begin{enumerate}[label=(\alph*)]
		\item $\prod_{i=1}^{4} \frac{\sin(\gamma_i)}{\sin(\beta_i)} r_i^{(-1)^i} = 1$ \label{Produkt}
		\item There exists a nontrivial conformal Combescure transformation of $\cbp$. This means that a checkerboard pattern with parallel edges exists where corresponding first order faces differ only by a similarity transformation. If one non-trivial conformal Combescure transformation exists, an entire two-parameter family of such transformations exists.\label{ParaNet}
		\item $\cbp$ is dualizable.\label{CFDual}
		\item $\cbp$ is a Koenigs checkerboard pattern.\label{Doliwa}
	\end{enumerate}
\end{satz}

\begin{bew}
	First we show that $\ref{Produkt} \iff \ref{ParaNet}$. If $j$ is even, the edge lengths $a_{j-1}$ and $a_{j}$ are related by the formula
	\begin{align*}
	a_{j-1} &= k_j a_j + c_j b,
	\end{align*}
	where
	\begin{align*}
	k_j := \frac{\sin(\gamma_j)}{\sin(\beta_j)} \quad \text{and} \quad c_j := \frac{\sin(\gamma_j + \delta_j)}{\sin(\beta_j)}.
	\end{align*}
	Analogously we find	$b_2 = k_3 b_3 + c_3 a$ and $b_4 = k_1 b_1 + c_1 a$. Using $a_i = b_i r_i$ we find the closing condition
	\begin{align*}
	a_1 &= k_1 k_2 k_3 k_4 \frac{r_2 r_4}{r_1 r_3} a_1 + a\left( k_2 r_2  \left( k_3 \left( \frac{k_4 r_4}{r_3}c_1 + c_4 r\right) + c_3 \right) + c_2 r\right).
	\end{align*}
	Hence we can compute $a_1$ only from the given angles and ratios if and only if
	\begin{align}
	\label{GeneralizedAngleCondition_neg}
	\prod_{i=1}^{4} k_i r_i^{(-1)^i} \neq 1.
	\end{align}
	Consequently a nontrivial parallel net with the same ratios $r$ and $r_i$ exists if and only if $\ref{Produkt}$ holds.
	
	Next we show that \ref{Produkt} $\iff$ \ref{CFDual}. When we dualize the net $\net$ all angles are replaced by their respective complements, i.e. $\alpha_i^*  = \pi - \alpha_i$, $\beta_i^* = \pi - \beta_i$, $\gamma_i^* = \pi - \gamma_i$ and $\delta_i^* = \pi - \delta_i$. Hence the coefficients $k_i$ are invariant under dualization while the coefficients $c_i$ change sign. If we denote the transformed edge lengths by $a_i^*, a^*$ and $b_i^*, b^*$ respectively, then the transformed relations read
	\begin{align*}
	a_{j-1}^* = k_j a_j^* - c_jb^* \quad \text{and} \quad b_{j-1}^* &= k_j b_j^* - c_j a^*.
	\end{align*}
	So the closing condition becomes
	\begin{align*}
	a_1^* = k_1 k_2 k_3 k_4 \frac{r_2 r_4}{r_1 r_3} a_1^* - a^*\left( k_2 r_2  \left( k_3 \left( \frac{k_4 r_4}{r_3}c_1 + c_4 r\right) + c_3 \right) + c_2 r\right).
	\end{align*}
	Again we find that $a_1^*$ can be determined from this equation if Equation \eqref{GeneralizedAngleCondition_neg} holds. However, comparing the potential formulas for $a_1$ and $a_1^*$ we find that $a_1^* = - a_1$. As no negative edge lengths can exist we conclude that $a_1^*$ exists only if \ref{Produkt} holds. On the other hand if \ref{Produkt} holds, we can construct a dual net for any value $a_1^*$ implying $\ref{CFDual}$.
	
	Next we show that \ref{Produkt} $\iff$ \ref{Doliwa}. To do so we use the inscribed angle theorem (see Theorem \ref{InscribedAngleTheoremHyp}). Let $k(\suppl{p_i}{p_j})$ denote the slope of the line $\suppl{p_i}{p_j}$ with respect to a coordinate system aligned with the asymptotes of the hyperbola, compare Theorem \ref{InscribedAngleTheoremHyp}. We find
	
	\begin{align*}
	k(\suppl{p_4}{p_5}) &= \pm \frac{\sin(\beta_1)\sin(\beta_4 + \alpha_4)}{r_4 \sin(\gamma_4)\sin(\gamma_1 + \delta_1)},
	&k(\overline{p_4p_6}) = \mp
	\frac{\sin(\gamma_3) \sin(\alpha_4 + \beta_4)}{r_3 \sin(\beta_4)\sin(\delta_3 + \gamma_3)},\\
	k(\suppl{p_5}{p_3}) &= \mp \frac{\sin(\gamma_1) \sin(\alpha_2 + \beta_2)}{r_1 \sin(\gamma_1 + \delta_1) \sin(\beta_2)},
	&k(\suppl{p_6}{p_3}) = \pm \frac{\sin(\beta_3) \sin(\alpha_2 + \beta_2)}{r_2 \sin(\gamma_2)\sin(\delta_3 + \gamma_3)}.
	\end{align*}
	
	Note that $\sin(\gamma_1 + \delta_1) = 0$ is equivalent to $\suppl{c_{\bar{1}}}{c_{\bar{1}2}} \parallel \suppl{c}{c_2}$ and thus is equivalent to $p_5 = p_2$. So if the points $p_1,\dots,p_6$ are all distinct, the denominators in the above equations are all nonzero.	Computing the quotients yields
	
	\begin{align*}
	\frac{k(\suppl{p_1}{p_6})}{k(\suppl{p_1}{p_4})} &= - \frac{r_3}{r_4} \frac{\sin(\beta_1)\sin(\beta_4)\sin(\delta_3 + \gamma_3)}{\sin(\gamma_3)\sin(\gamma_4)\sin(\gamma_1 + \delta_1)},\\
	\frac{k(\suppl{p_3}{p_6})}{k(\suppl{p_3}{p_4})} &= - \frac{r_2}{r_1} \frac{\sin(\gamma_1) \sin(\gamma_2) \sin(\delta_3 + \gamma_3)}{\sin(\beta_2)\sin(\beta_3)\sin(\gamma_1 + \delta_1)}.
	\end{align*}
	By Theorem \ref{InscribedAngleTheoremHyp} the points $p_1,\dots,p6$ lie on a common hyperbola if and only if
	\begin{align*}
	\frac{r_2}{r_1} \frac{\sin(\gamma_1) \sin(\gamma_2)}{\sin(\beta_2)\sin(\beta_3)} &=
	\frac{r_3}{r_4} \frac{\sin(\beta_1)\sin(\beta_4)}{\sin(\gamma_3)\sin(\gamma_4)}
	\end{align*}
	which is equivalent to $\ref{Produkt}$. This concludes the proof.
\end{bew}

\begin{bem}
	If we find $p_i = p_j$ for some $i \neq j$ everything in the proof of Theorem \ref{Theorem:KoenigsNets} still holds except for the application of Theorem \ref{InscribedAngleTheoremHyp}. So in such a case we still find that $\ref{Produkt} \iff \ref{ParaNet} \iff \ref{CFDual}$. 
\end{bem}

\begin{kor}
	Let $\cbp[\net]$ be a conjugate checkerboard pattern with control net $\net$. Then $\cbp[\net]$ is dualizable if and only if each two Laplace invariants defined in the faces of $\net$ are equal.
\end{kor}

\begin{bew}
	If the six points $p_1,\dots,p_6$ from Definition \ref{Defi:KoenigsNet} are distinct, the statement follows from Theorem \ref{Theorem:KoenigsNets}. Hence condition \ref{Produkt} in Theorem \ref{Theorem:KoenigsNets} is equivalent to the multiplicative one-form $q$ being closed if $p_1,\dots,p_6$ are all distinct. However these terms depend continuously on the vertices of the checkerboard pattern. Hence, any face patch, on witch $q$ is closed, can be approximated with a sequence of dualizable face patches where $p_1,\dots,p_6$ are distinct. Since condition $\ref{Produkt}$ is preserved in the limit, so is the existence of a dual.
\end{bew}

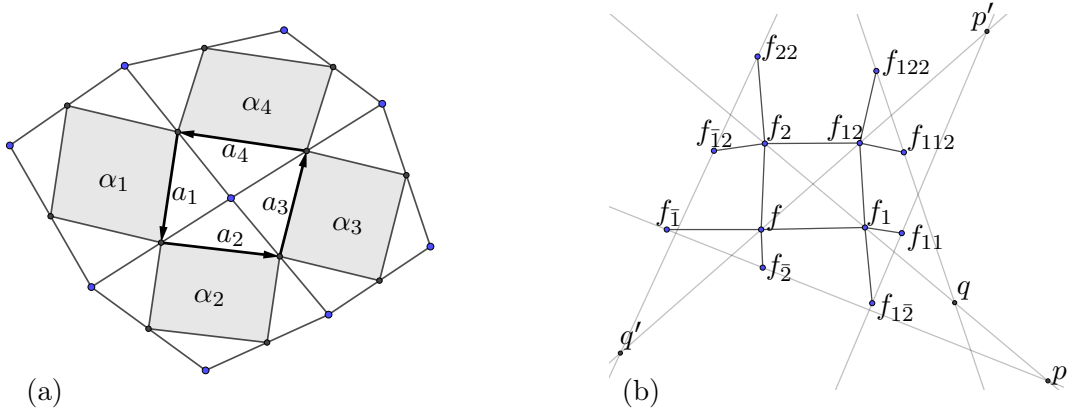
\begin{figure}[h]
	\begin{overpic}[width=0.45\linewidth]
		{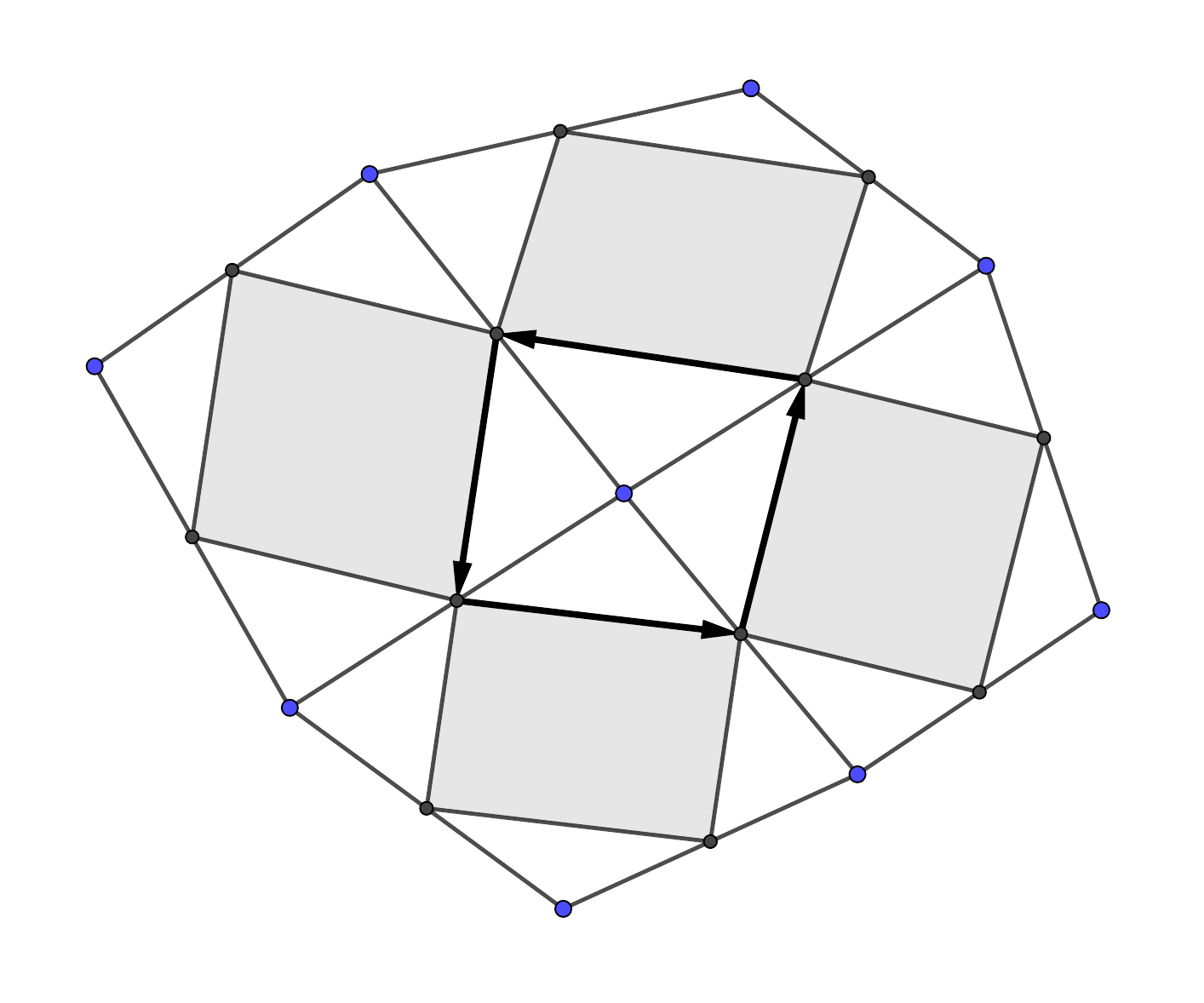}
		\put(40,42){$a_1$}
		\put(49,33.5){$a_2$}
		\put(58,40){$a_3$}
		\put(50,50.5){$a_4$}
		\put(25.5,44.5){$\alpha_1$}
		\put(44.5,21.5){$\alpha_2$}
		\put(72,36){$\alpha_3$}
		\put(54,60){$\alpha_4$}
	\end{overpic}
	\hfill
	\begin{overpic}[width=0.45\linewidth]
		{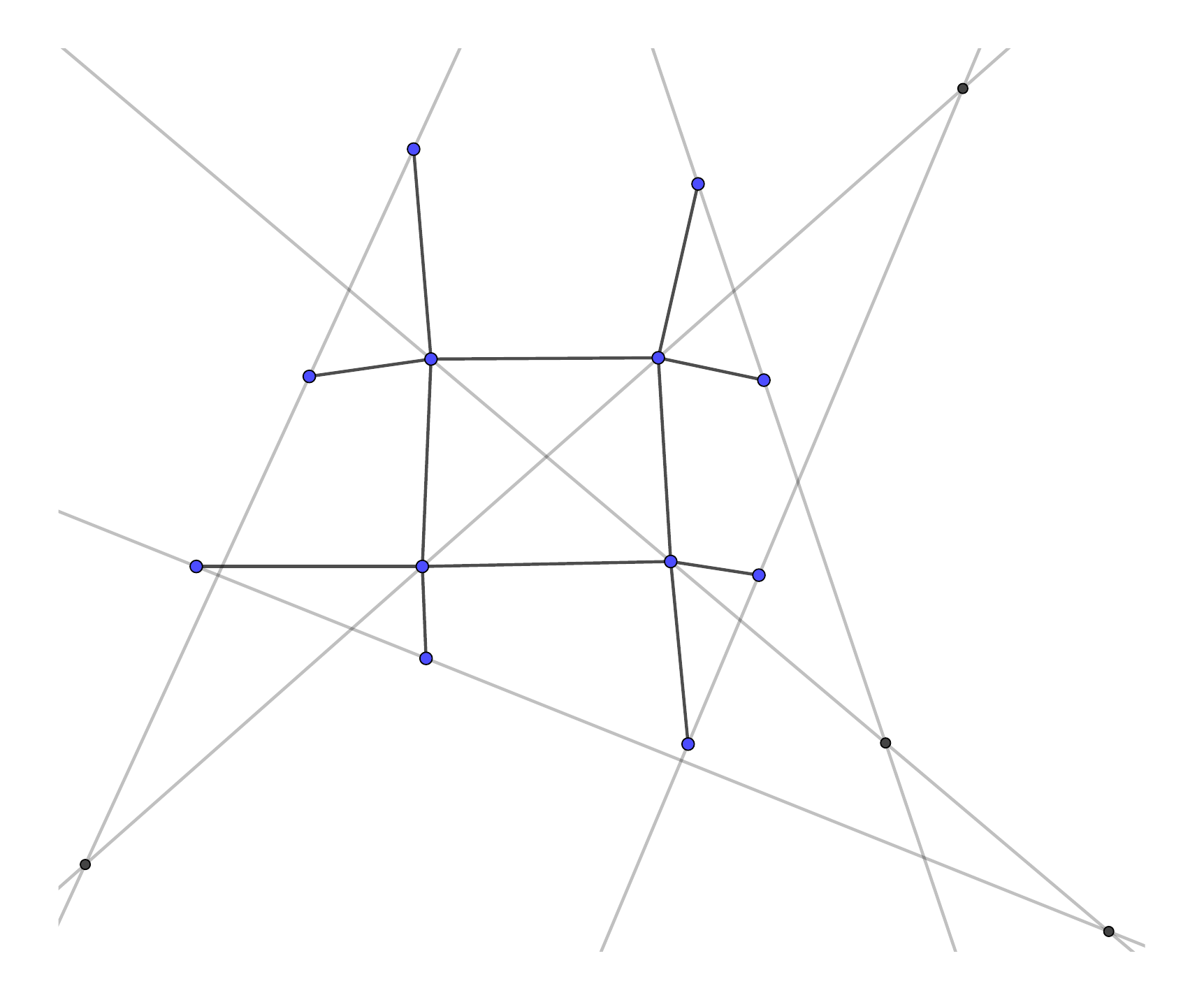}
		\put(36,37){$f$}
		\put(56,38){$f_1$}
		\put(48,55){$f_{12}$}
		\put(36,55){$f_2$}
		\put(14,38){$f_{\bar{1}}$}
		\put(36,27){$f_{\bar{2}}$}
		\put(58,19){$f_{1\bar{2}}$}
		\put(64,33){$f_{11}$}
		\put(64,53){$f_{112}$}
		\put(59,68){$f_{122}$}
		\put(35,71){$f_{22}$}
		\put(22,54){$f_{\bar{1}2}$}
		\put(7,13){$q'$}
		\put(93,6){$p$}
		\put(74,23){$q$}
		\put(77,77){$p'$}
	\end{overpic}
	\vskip -1.5em \leftline{\hskip 0.05\linewidth (a)\hskip 0.5\linewidth (b)}
	\caption{(a): The edges $a_i$ have to close in the initial net as well as in the dualized net. From this condition the scaling factors that guide the dualization can be computed. (b): The setting of Lemma \ref{Lemma:2dConstructionOfKoenigs}.}
	\label{Fig:StableAlgorithm}
\end{figure}

\begin{bem}
	For a given Koenigs checkerboard pattern, there is a two-parameter family of dual checkerboard patterns that differ in the scaling of corresponding first order faces. By choosing the initial scaling factors of two adjacent first order faces $\alpha_1$ and $\alpha_2$, all other scaling factors can be computed recursively by the formulas
	\begin{align*}
	\alpha_2 = \frac{\left\| (\alpha_1 a_1 - \alpha_0 a_0) \times a_3\right\|}{\|a_2 \times a_3 \|} \quad \text{and} \quad \alpha_3 = \frac{\left\| (\alpha_1 a_1 - \alpha_0 a_0) \times a_2\right\|}{\|a_2 \times a_3 \|},
	\end{align*}
	where $a_i$ are the oriented edges of the corresponding first order faces, see Figure \ref{Fig:StableAlgorithm}, left. This permits a stable dualization algorithm.
\end{bem}

The following lemma provides an easy way to generate Koenigs nets in $\Pro\R^2$.

\begin{lem}\label{Lemma:2dConstructionOfKoenigs}
	Let $M$ and $N$ be two commuting projective transformations $\Pro\R^2 \to \Pro\R^2$ and let $P\in\Pro\R^2$. Then the net $\net$ defined by
	\begin{align*}
	f(k,l) = M^k N^l P
	\end{align*}
	is the control net of a Koenigs checkerboard pattern in $\Pro\R^2$.
\end{lem}

\begin{bew}
	We show that the condition of Theorem \ref{Theorem:LapInvariants} is met in the quadrilateral $(f, f_{1}, f_{12}, f_{2})$, compare Figure \ref{Fig:StableAlgorithm}, right. Let
	\begin{align*}
	\begin{array}{lcr}
	p\phantom{'} = \suppl{f_1}{f_2} \cap \suppl{f_{\bar{1}}}{f_{\bar{2}}}, &\qquad  \qquad&
	q = \suppl{f_1}{f_2} \cap \suppl{f_{112}}{f_{112}}, \\
	p'= \suppl{f}{f_{12}} \cap \suppl{f_{1\bar{2}}}{f_{11}}, & \qquad \qquad&
	q'= \suppl{f}{f_{12}} \cap \suppl{f_{22}}{f_{\bar{1}2}}.\hspace*{3pt}
	\end{array}
	\end{align*}
	Let $F: \Z^2 \to \R^3$ be the net of homogeneous coordinates of the vertices of $\net$. We find
	\begin{align*}
	P &= (F_{\bar{1}} \times F_{\bar{2}}) \times (F_1 \times F_2) = \innerprod{F_{\bar{1}} \times F_{\bar{2}}}{F_2} F_1 - \innerprod{F_{\bar{1}} \times F_{\bar{2}}}{F_{1}} F_2 =\\
	&= \det(F_{\bar{1}}, F_{\bar{2}}, F_2) F_1 - \det(F_{\bar{1}}, F_{\bar{2}}, F_1) F_2\\
	Q &= (F_{112} \times F_{122}) \times (F_1 \times F_2) = \innerprod{F_{112} \times F_{122}}{F_2} F_1 - \innerprod{F_{112} \times F_{122}}{F_{1}} F_2 =\\
	&= \det(F_{112}, F_{122}, F_2) F_1 - \det(F_{112}, F_{122}, F_1) F_2\\
	P' &= (F_{11} \times F_{1\bar{2}}) \times (F \times F_{12}) = \innerprod{F_{11} \times F_{1\bar{2}}}{F_{12}} F - \innerprod{F_{11} \times F_{1\bar{2}}}{F} F_{12} =\\
	&= \det(F_{11}, F_{1\bar{2}}, F_{12}) F - \det(F_{11}, F_{1\bar{2}}, F) F_{12}\\
	Q' &= (F_{\bar{1}2} \times F_{22}) \times (F \times F_{12}) = \innerprod{F_{\bar{1}2} \times F_{22}}{F_{12}} F - \innerprod{F_{\bar{1}2} \times F_{22}}{F} F_{12} =\\
	&= \det(F_{\bar{1}2}, F_{22}, F_{12}) F - \det(F_{\bar{1}2}, F_{22}, F) F_{12}.
	\end{align*}
	From this we can formulate the cross ratios as
	\begin{align*}
	\DV(f_1,f_2,p,q) &= \DV(F_1, F_2, P, Q) = \frac{\det(F_{\bar{1}},F_{\bar{2}},F_1) \det(F_{112}, F_{122}, F_2)}{\det(F_{\bar{1}}, F_{\bar{2}}, F_2) \det(F_{112}, F_{122}, F_1)}\\
	\DV(f, f_{12},p',q') &= \DV(F, F_{12}, P', Q') = \frac{\det(F_{11}, F_{1\bar{2}}, F) \det(F_{\bar{1}2}, F_{22}, F_{12})}{\det(F_{11}, F_{1\bar{2}}, F_{12}) \det(F_{\bar{1}2}, F_{22}, F)}.
	\end{align*}
	Now let $\Mat{M}$ and $\Mat{N}$ be the matrix representations of $M$ and $N$ in homogeneous coordinates. Then we can express the cross ratios as
	\begin{align*}
	\DV(f_1, f_2,p,q) &= \frac{\det(\Mat{M}^{-1} F, \Mat{N}^{-1}F, \Mat{M}F) \det(\Mat{M}\Mat{M}\Mat{N}F, \Mat{M}\Mat{N}\Mat{N}F, \Mat{N}F)}{\det(\Mat{M}^{-1}F, \Mat{N}^{-1} F, \Mat{N}F) \det(\Mat{M}\Mat{M}\Mat{N}F, \Mat{M}\Mat{N}\Mat{N}F, \Mat{M}F)} = \\
	&= \frac{\det(F, \Mat{M}\Mat{N}^{-1}F, \Mat{M}\Mat{M}F) \det(\Mat{M}\Mat{N}F, \Mat{N}\Mat{N}F, \Mat{M}^{-1}\Mat{N}F)}{\det(\Mat{N}\Mat{M}^{-1}F, F, \Mat{N}\Mat{N}F) \det(\Mat{M}\Mat{M}F, \Mat{M}\Mat{N}F, \Mat{N}^{-1}\Mat{M}F)}
	\end{align*} 
	and
	\begin{align*}
	\DV(f,f_{12},p',q') &= \frac{\det(\Mat{M}\Mat{M}F, \Mat{M}\Mat{N}^{-1}F, F) \det(\Mat{M}^{-1}\Mat{N}F, \Mat{N}\Mat{N}F, \Mat{M}\Mat{N}F)}{\det(\Mat{M}\Mat{M}F, \Mat{M}\Mat{N}^{-1}F, \Mat{M}\Mat{N}F) \det(\Mat{M}^{-1}\Mat{N}F, \Mat{N}\Mat{N}F, F)}.
	\end{align*}
	So we see that the two Laplace invariants $\DV(f_1, f_2,p,q)$ and $\DV(f,f_{12},p',q')$ are equal.
\end{bew}

\section{Isothermic nets} \label{Section:IsothermicNets}

Discrete isothermic nets can now be defined as principal nets that are also Koenigs nets. Analogous to the smooth case or to other discrete approaches \cite{KoenigsNetsBobenkoSuris} we find that the class of discrete isothermic nets is invariant under dualizations and Möbius transformations. This permits a construction of discrete minimal surfaces and their Goursat transformations as will be described later on.

\begin{defi}
	\label{Defi:AlternativeIsoSurf}
	We call a checkerboard pattern $\cbp$ \emph{isothermic}, if it is principal and Koenigs.
\end{defi}

\begin{figure}[h!]
	\includegraphics[width=0.65\linewidth]{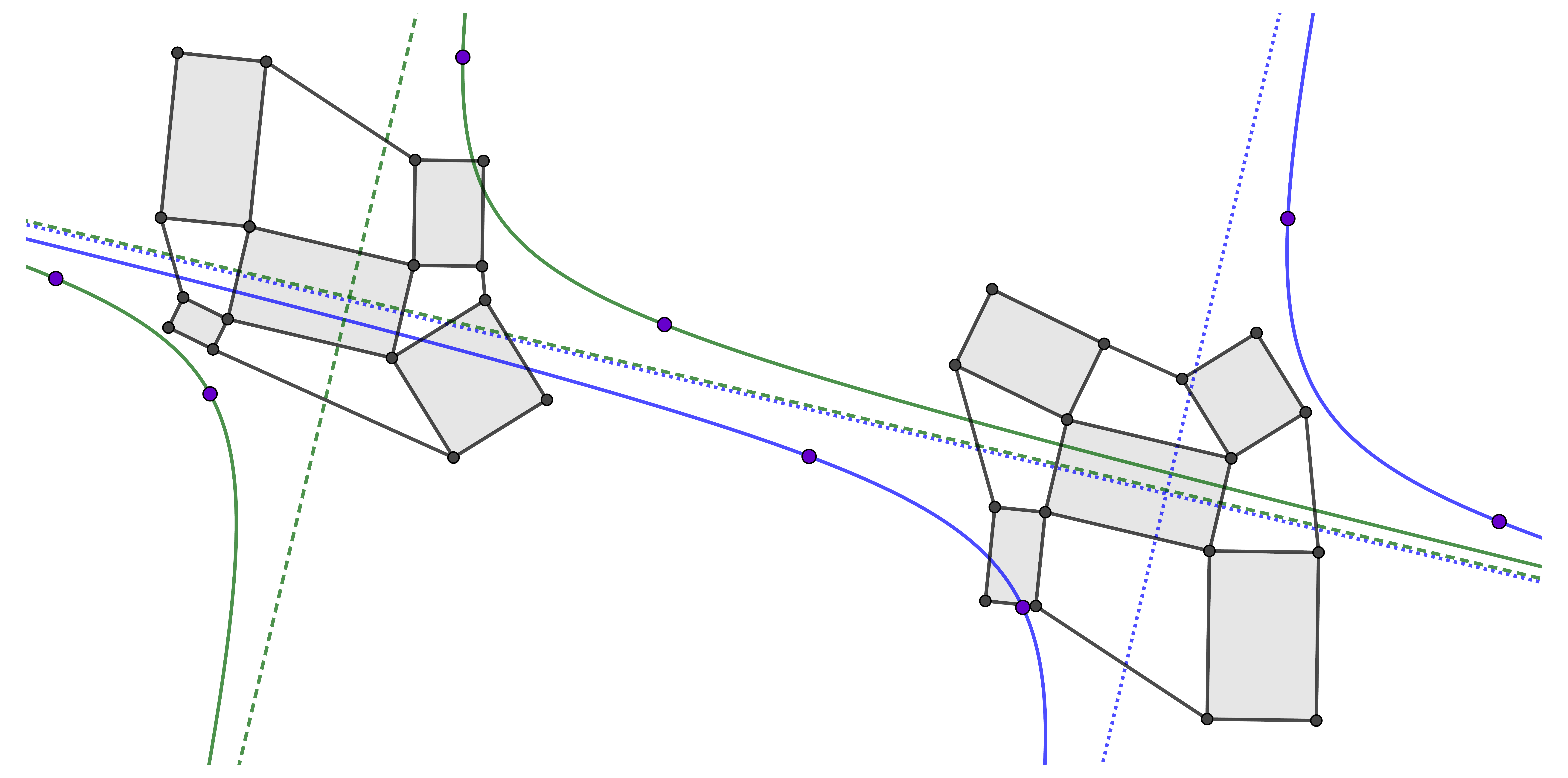}
	\hfill
	\begin{minipage}[b]{0.3\linewidth}
		\caption{An isothermic checkerboard pattern and its dual with the corresponding conics of Koenigs. The points on the hyperbolas are the points of intersecting supporting lines of neighboring edges.\vspace*{1em}}
		\label{Fig:CFTrafoOfGeneralIso}
	\end{minipage}
\end{figure}

As orthogonal first order faces are mapped to orthogonal faces under dualization, the next corollary follows immediately from Theorem \ref{Theorem:KoenigsNets}. See Figure \ref{Fig:CFTrafoOfGeneralIso} for an illustration.

\begin{kor}
	\label{Kor:CFDual}
	Isothermic checkerboard patterns are dualizable. Their dual is again an isothermic checkerboard pattern.
\end{kor}

\begin{figure}
	\includegraphics[width=0.65\linewidth]{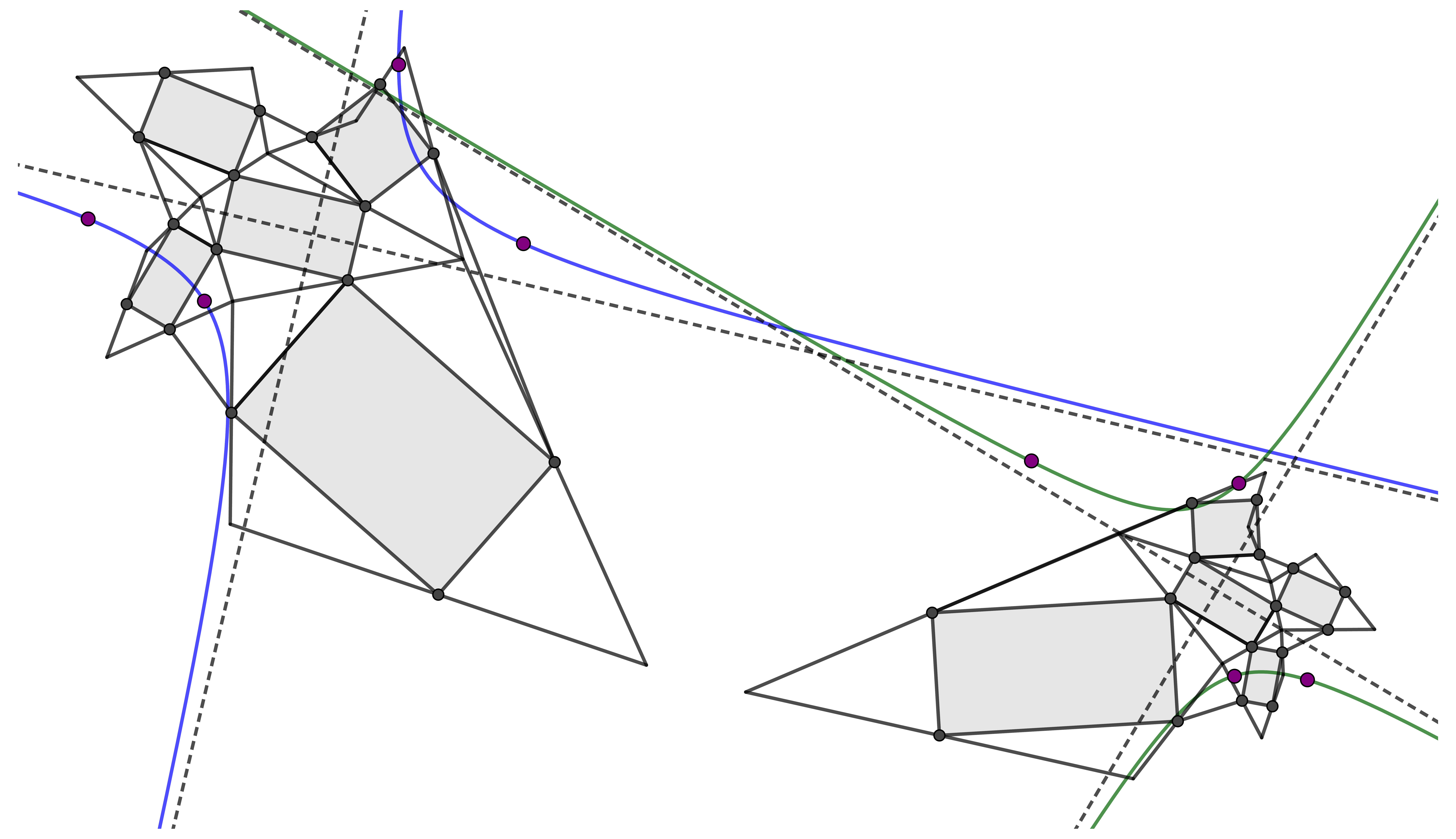}
	\hfill
	\begin{minipage}[b]{0.3\linewidth}
		\caption{An isothermic checkerboard pattern and its Möbius transform together with the corresponding conics of Koenigs. The figure features non-convex quads as the examples were constructed in such a way that the points of intersecting lines are all close to the checkerboard pattern.\vspace*{0em}}
		\label{Fig:IsothermicSurfPlusMoebTrafo}
	\end{minipage}
\end{figure}

\begin{satz}[Möbius invariance]
	\label{Theorem:MoebiusInvarianceOfIsothermic}
	Isothermic checkerboard patterns are mapped to isothermic checkerboard patterns under a discrete Möbius transformation, see Figure \ref{Fig:IsothermicSurfPlusMoebTrafo}.
\end{satz}

The proof will be a direct consequence of Lemma \ref{Lemma:LapInvariantInMink} and is thus postponed for now. In order to the prove Theorem \ref{Theorem:MoebiusInvarianceOfIsothermic}, we study isothermic nets again in the space $\Mink$ under the embedding $\iota$. We have already defined pseudo-principal nets in $\Mink$ and can now extend them to pseudo-isothermic nets.

\begin{defi}
	We call a net $\nesh$ in $\Mink$ pseudo-isothermic if it is pseudo-principal and the two Laplace invariants for each face are equal.
\end{defi}

It turns out that the lift $\iota(\net)$ of an isothermic net $\net$ in $\R^3$ is a pseudo-isothermic net in $\Mink$ as the following lemma shows.

\begin{lem}\label{Lemma:LapInvariantInMink}
	Let $\net$ be the control net of an isothermic checkerboard pattern and $\iota(\net)$ be its lift to $\Mink$. The Laplace invariants of corresponding faces of $\net$ and $\iota(\net)$ are equal.
\end{lem}

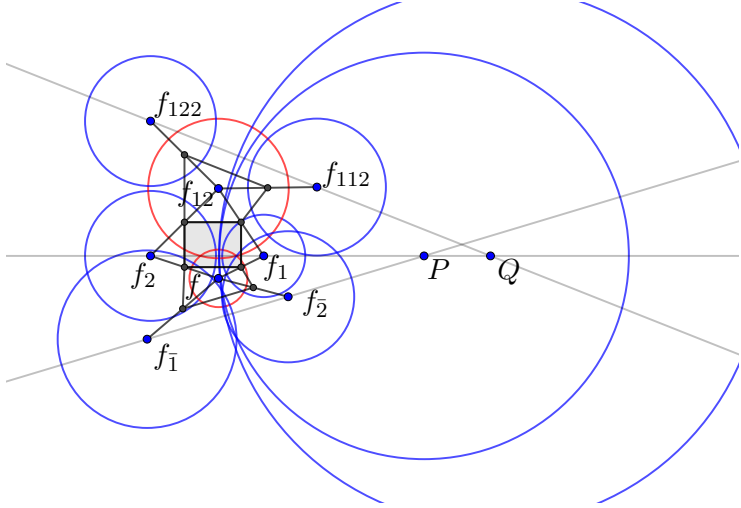
\begin{figure}
\begin{overpic}[width=0.7\linewidth]{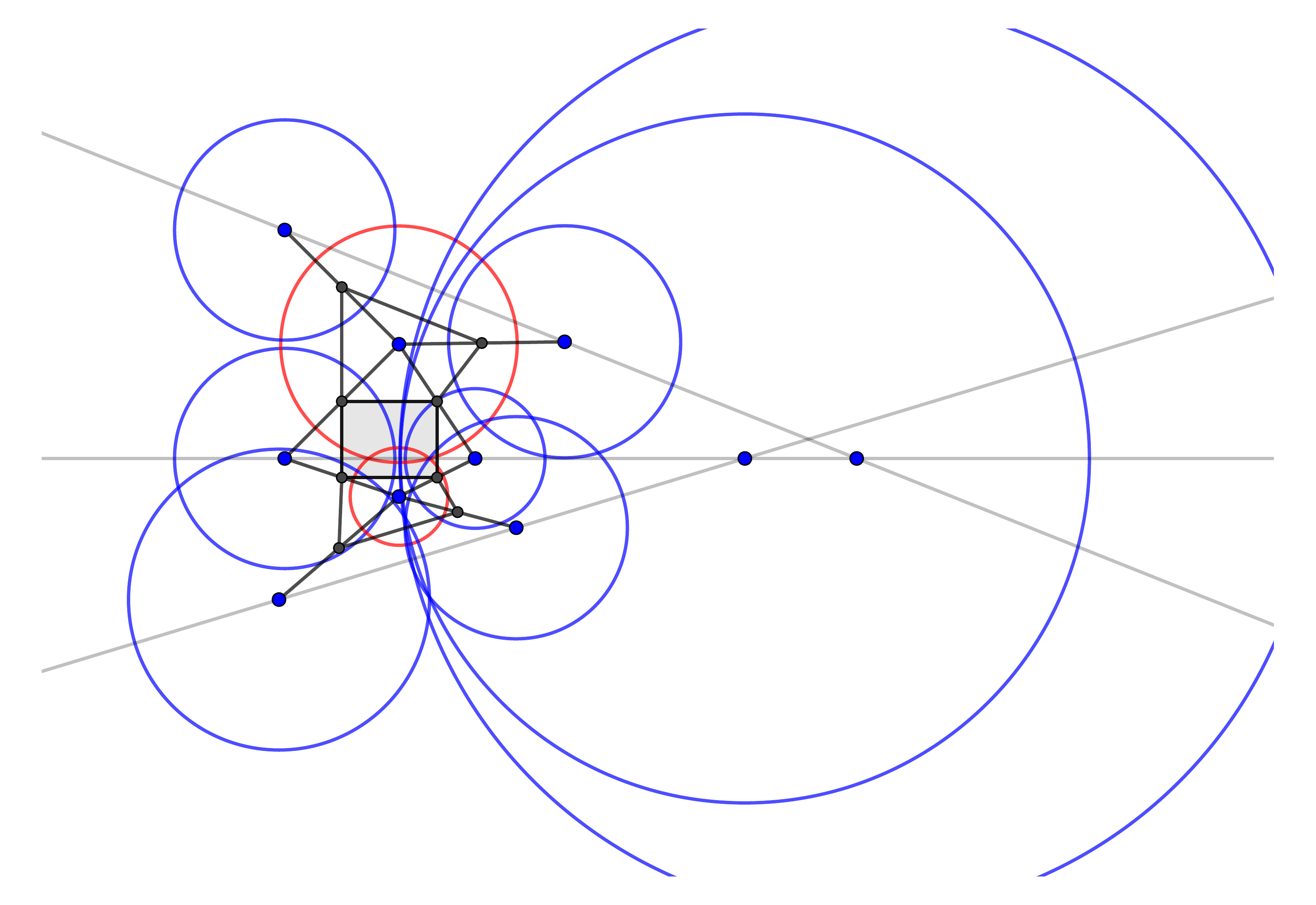}
	\put(36,31.5){$f_1$}
	\put(19,30.5){$f_2$}
	\put(57,31){$P$}
	\put(66,31){$Q$}
	\put(43.5,43.5){$f_{112}$}
	\put(22,52.5){$f_{122}$}
	\put(41,27){$f_{\bar{2}}$}
	\put(22,20){$f_{\bar{1}}$}
	\put(26,29){$f$}
	\put(25,41){$f_{12}$}
\end{overpic}
\hfill
\begin{minipage}[b]{0.25\linewidth}
	\caption{The idea behind the proof of Lemma \ref{Lemma:LapInvariantInMink}: Not only do the lines $\suppl{f_1}{f_2}$ and $\suppl{f_{\bar{1}}}{f_{\bar{2}}}$ intersect in $P$, but also the corresponding pencils of spheres intersect in a sphere with center at $P$. This means that there is a sphere with center at $P$ that intersects both the sphere with center at $f$ and the sphere with center at $f_{12}$ orthogonally.\vspace*{0.5em}}
	\label{Fig:LapInvariantWithSpheres}
\end{minipage}
\end{figure}

\begin{bew}
	First note that $\iota(\net)$ has a conjugate checkerboard pattern and thus the Laplace invariants are well defined. Hence not only the supporting lines $\suppl{f_{\bar{1}}}{f_{\bar{2}}}$ and $\suppl{f_1}{f_2}$ intersect, but also the corresponding pencils of spheres, compare Figure \ref{Fig:LapInvariantWithSpheres}. However, we know that the first three components under the lift $\iota$ are the same as the original centers of spheres and when we compute the cross ratio of points lying on a line it is sufficient to use just one coordinate. So it follows that the Laplace invariants remain unchanged under $\iota$.
\end{bew}

From Lemma \ref{Lemma:LapInvariantInMink} the proof of Theorem \ref{Theorem:MoebiusInvarianceOfIsothermic} follows immediately.

\begin{bew}[Proof of Theorem \ref{Theorem:MoebiusInvarianceOfIsothermic}.]
	Every Möbius transformation can be seen as a projective transformation in $\Mink$ that preserves the inner product. Obviously these transformations preserve the cross ratio and since $\iota$ also preserves the Laplace invariants we can conclude that not only conjugacy and orthogonality, but also the Koenigs property is preserved under Möbius transformations.
\end{bew}

\subsection{Minimal surfaces}\label{Section:MinimalSurfaces}

Minimal surfaces can be constructed by dualizing an isothermic net on the unit sphere, since the theory of minimal surfaces tells us that for any minimal surface its dual and its Gauß image are equal. With the Möbius transformation and dualization at hand we can reproduce this construction in the discrete setting.

\begin{defi}
	Let $\net$ be the control net of an isothermic checkerboard pattern $\cbp[f]$. We call $\cbp[\net]$ minimal if it has a dual checkerboard pattern $\cbp'$ that is also the checkerboard pattern of a principal Gauß image of $f$ in the sense of Definition \ref{Defi:NormalVector2}.
\end{defi}

\begin{defi}
	Let $\net$ and $\tilde{\net}$ be control nets of minimal checkerboard patterns. They are related by a \emph{Goursat transformation} if their principal Gauß images are related by a Möbius transformation.
\end{defi}

\begin{defi}
	We say that a checkerboard pattern $\cbp[\net]$ is \emph{on the unit sphere}, if there is a Möbius representation $\moeb[\net]$ where every sphere intersects the unit sphere orthogonally.
\end{defi}

\begin{kor}
	Let $\cbp[n]$ be an isothermic checkerboard pattern on the unit sphere. The dual checkerboard pattern $\cbp[n]'$ is a minimal checkerboard pattern and $n$ is its principal Gauß image. If $n$ is used to compute the discrete shape operator of $\cbp[n]'$, the mean curvature of $\cbp[n]'$ is zero.
\end{kor}

\begin{bew}
	The first statement follows directly from the definition of minimal checkerboard patterns. The principal curvature $\kappa_1$ and $\kappa_2$ are just the oriented scaling factors between edges of $\cbp[n]$ and $\cbp[n]'$. If the Gauß image is the dual net at the same time the relation $\kappa_1 = - \kappa_2$ holds.
\end{bew}


\begin{figure}
	\includegraphics[width=0.32\linewidth]{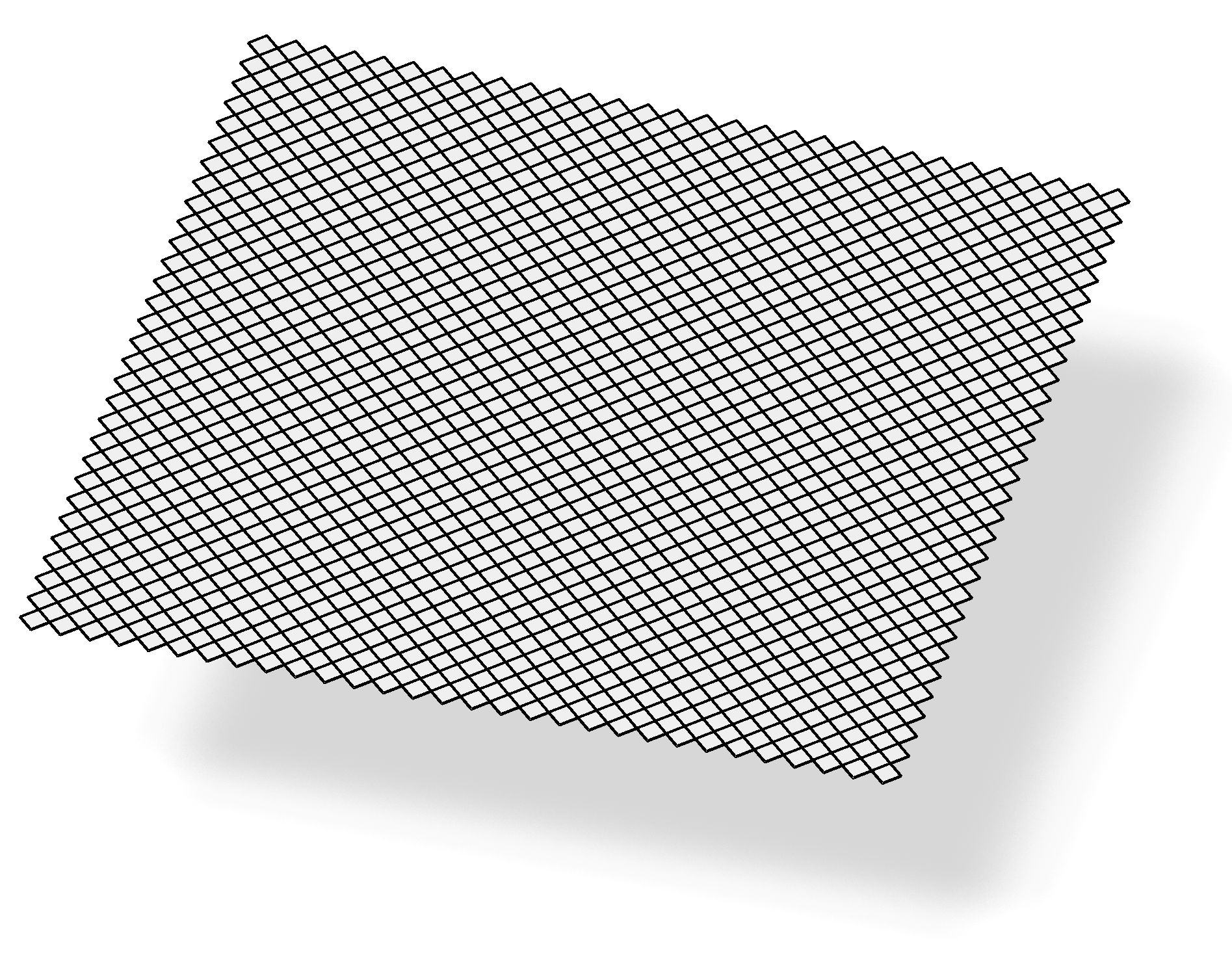}
	\includegraphics[width=0.32\linewidth]{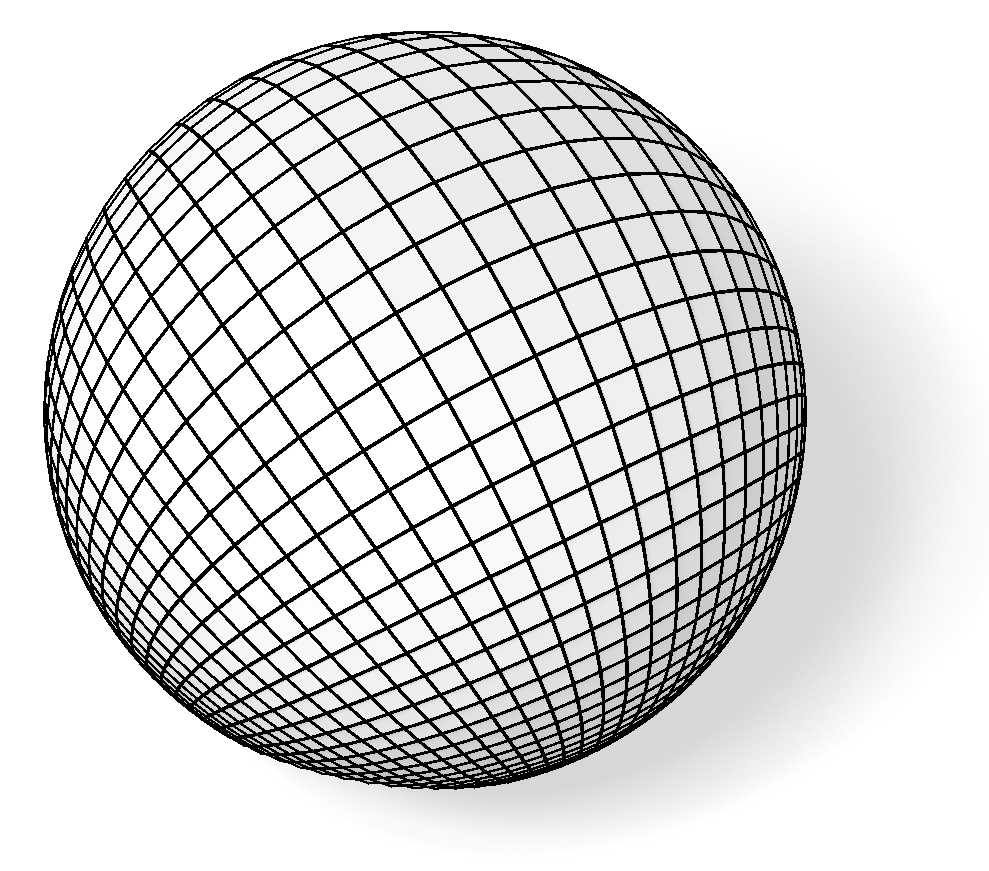}
	\includegraphics[width=0.32\linewidth]{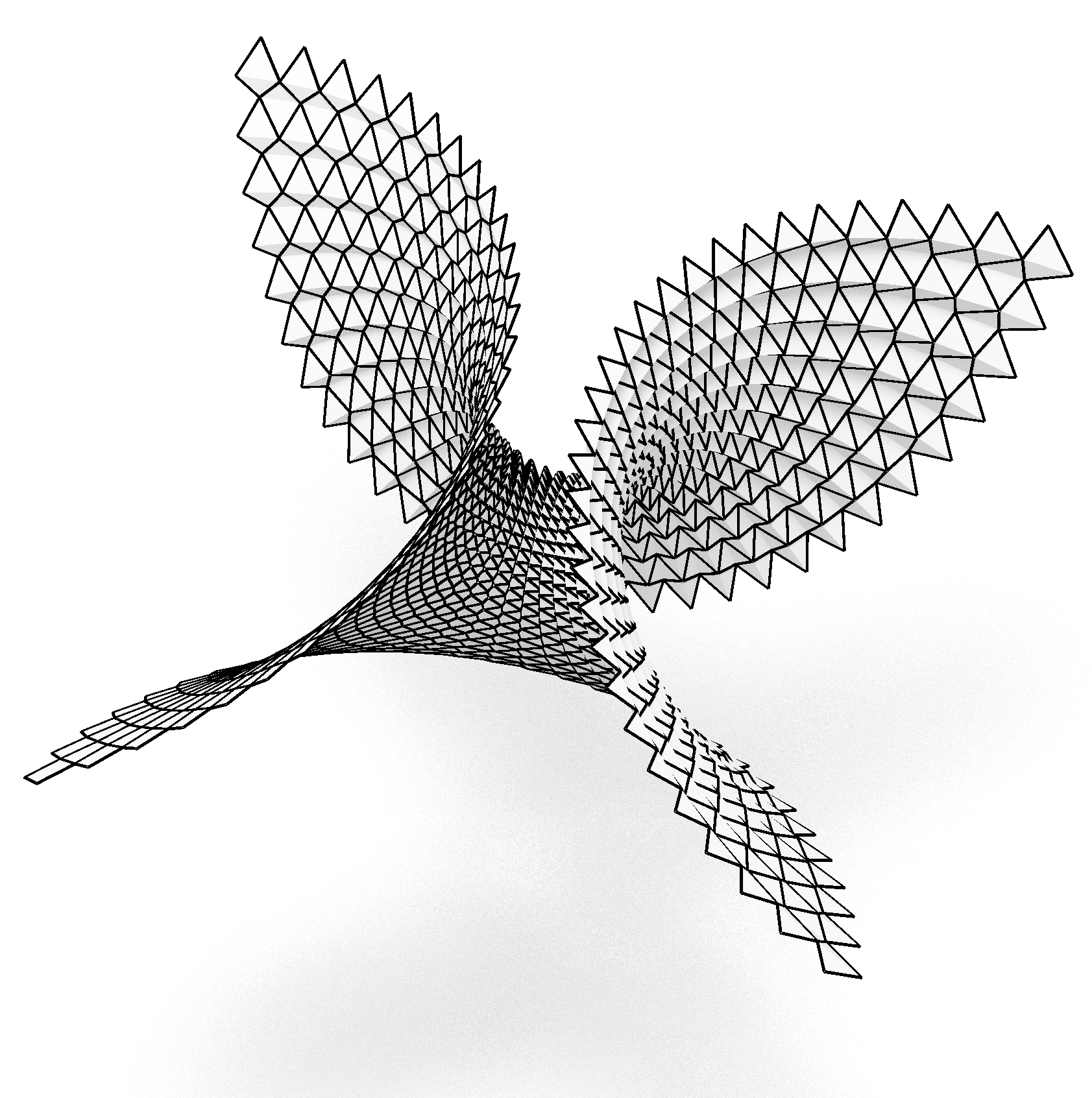}
	\includegraphics[width=0.32\linewidth]{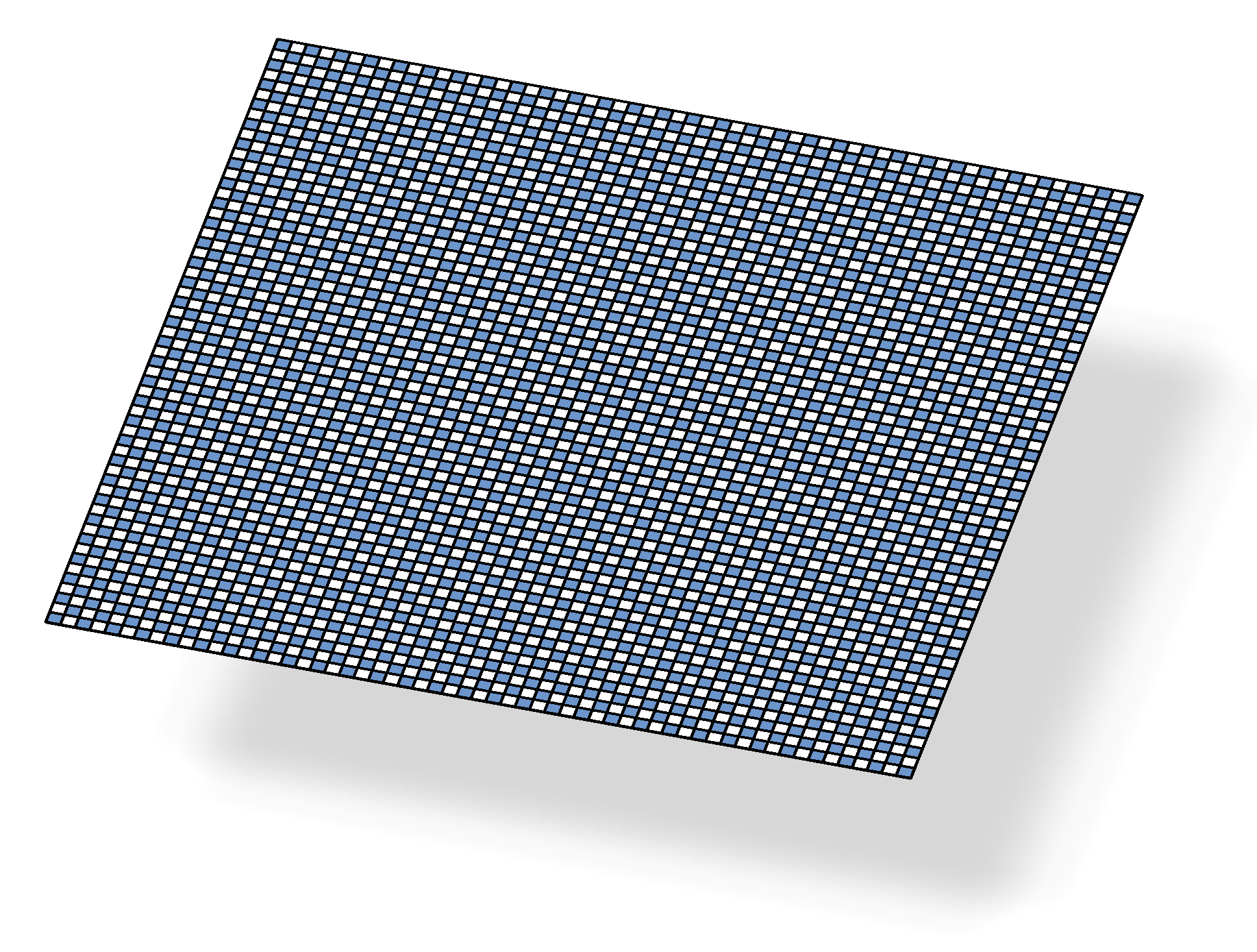}
	\includegraphics[width=0.32\linewidth]{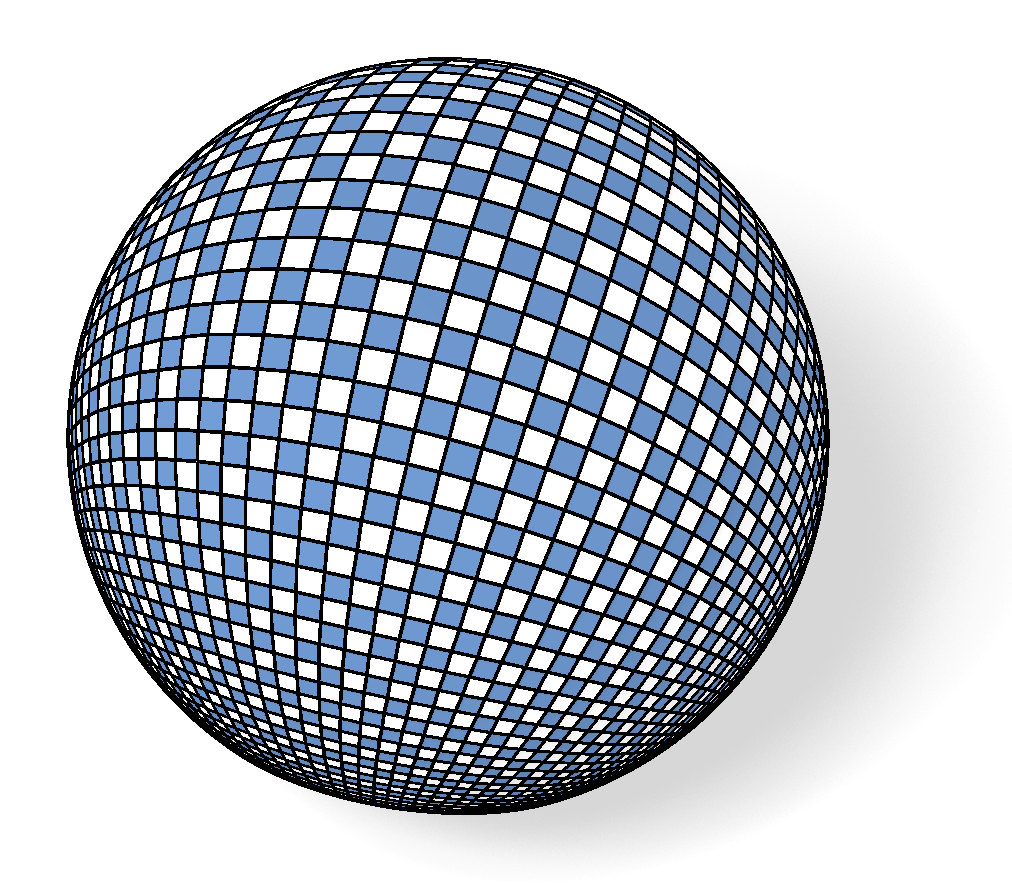}
	\includegraphics[width=0.32\linewidth]{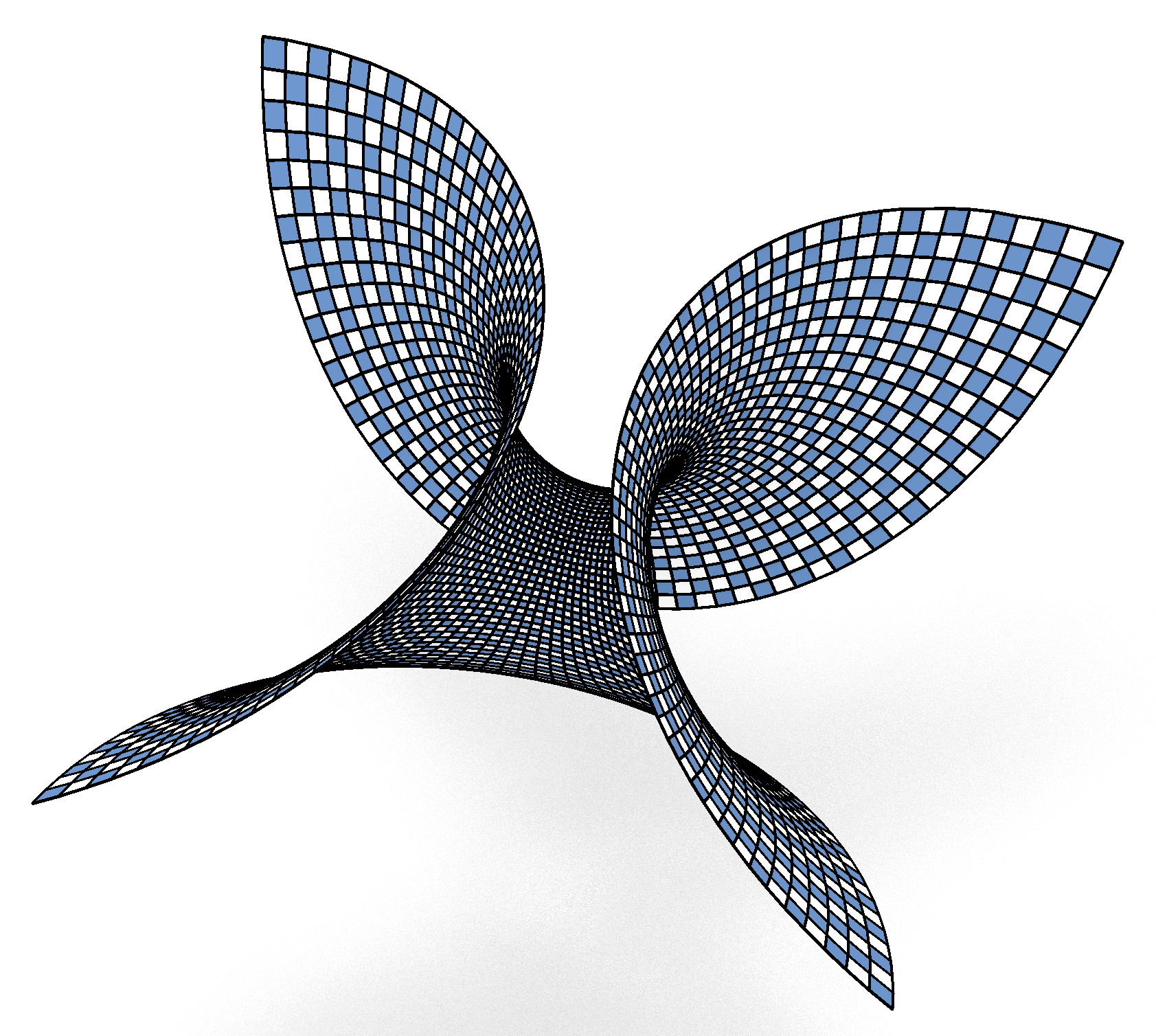}
	\caption{Enepper surface: In the top row we see from left to right the Weierstrass data of the Enepper surface, the Gauß image of the Enepper surface and the Enepper surface itself. In the second row we see the checkerboard patterns of the corresponding nets.}
	\label{Fig:EnepperSurface}
\end{figure}

\begin{figure}
	\includegraphics[width=0.3\linewidth]{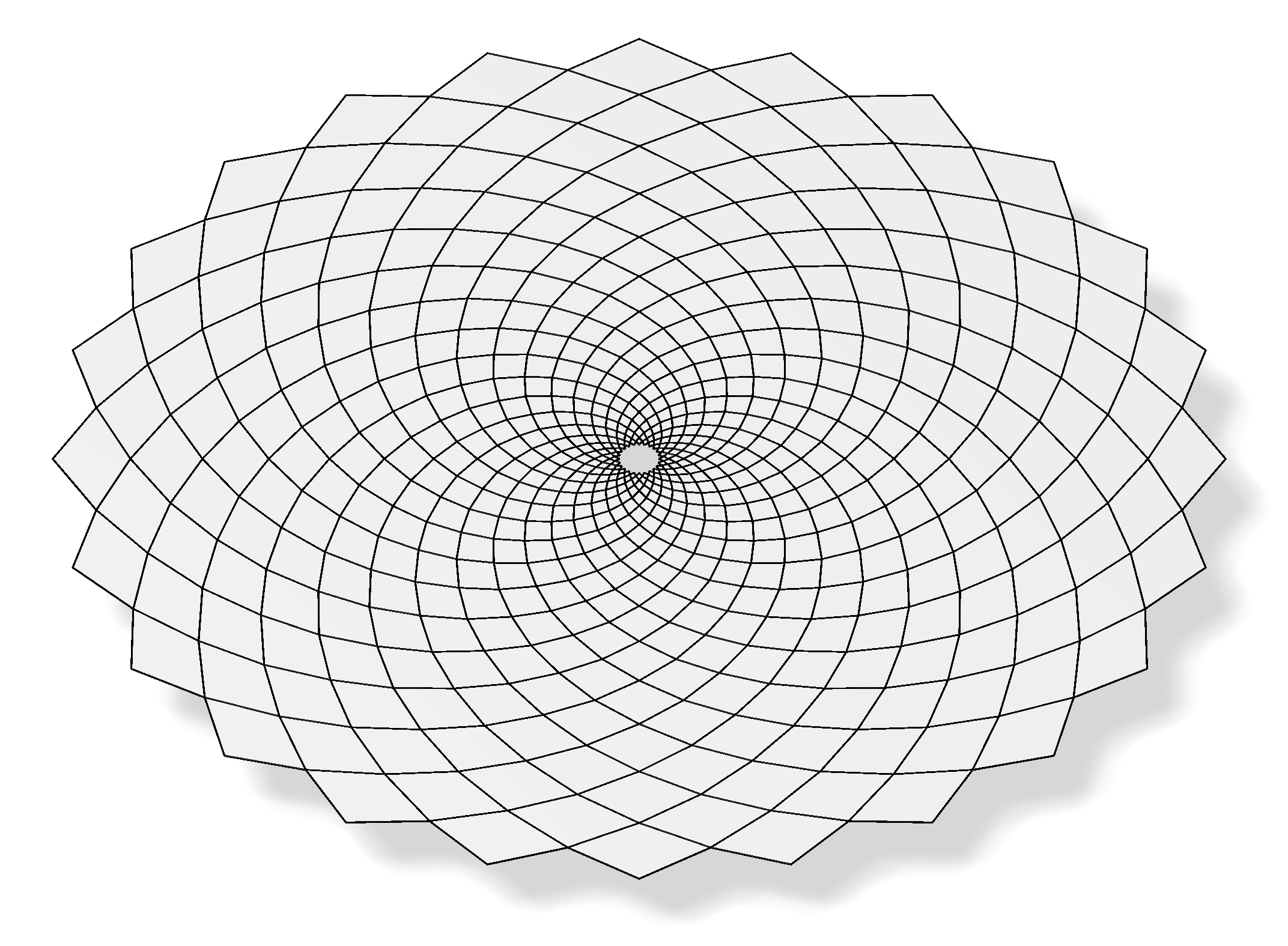}
	\includegraphics[width=0.3\linewidth]{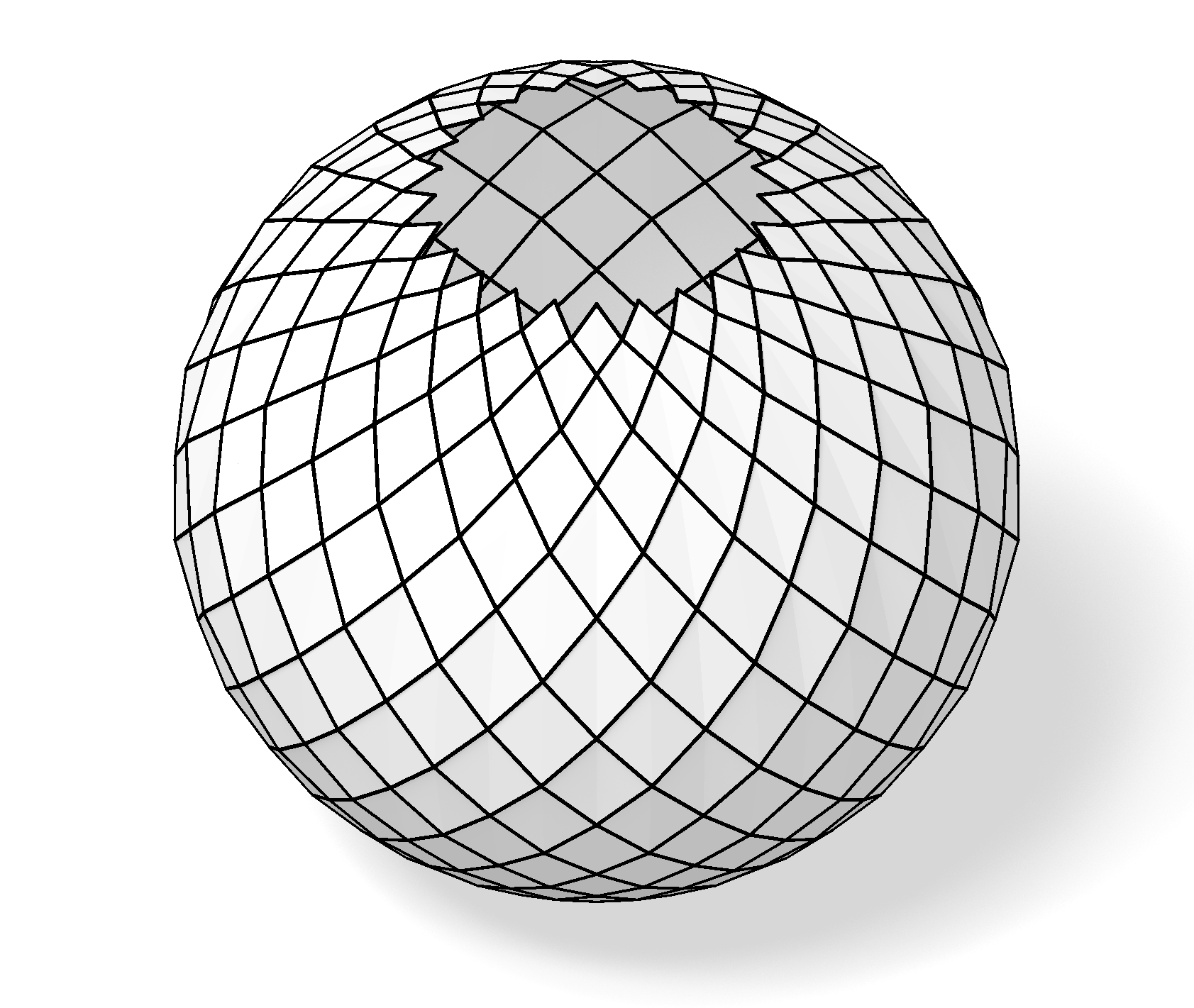}
	\includegraphics[width=0.3\linewidth]{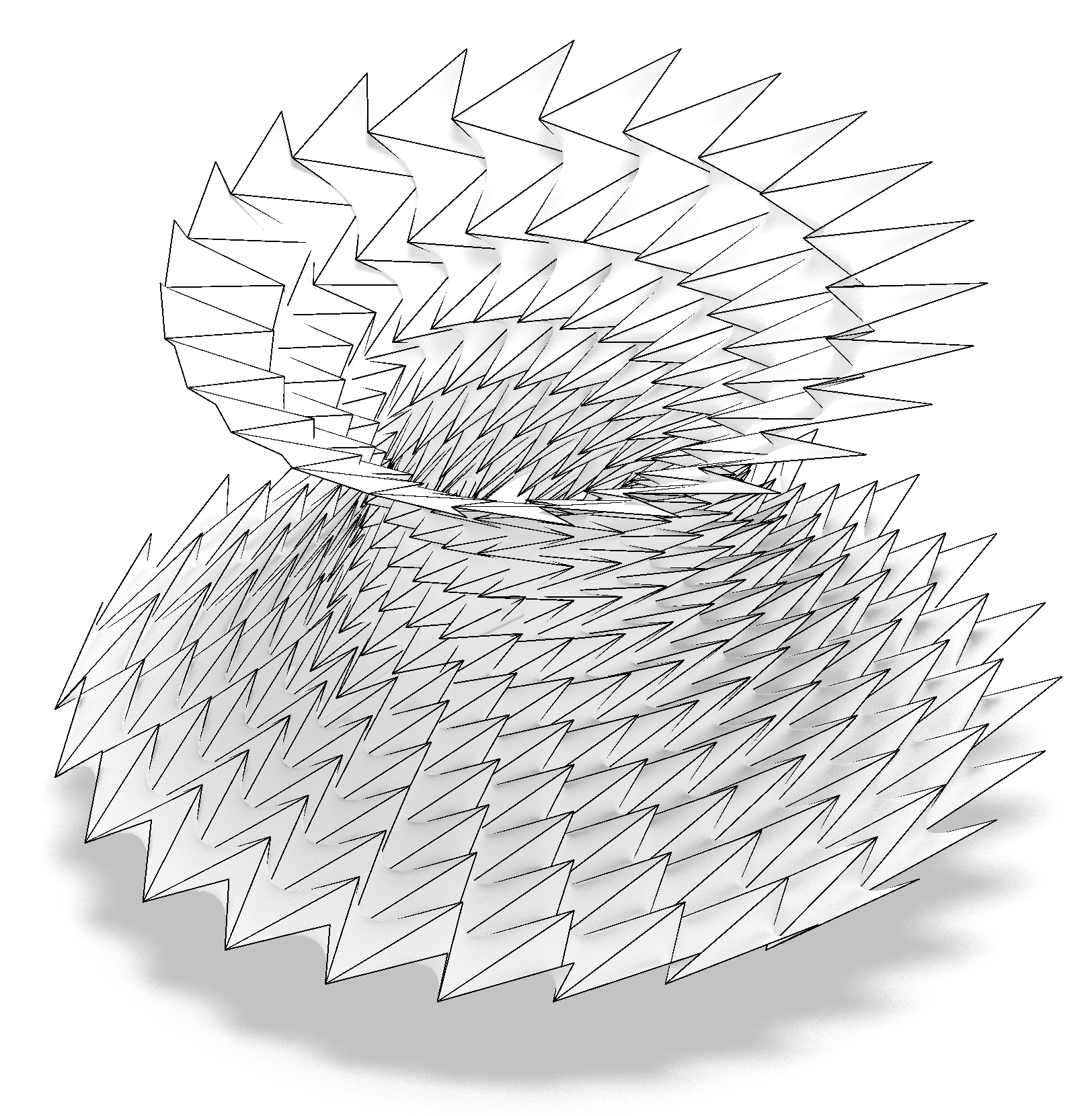}
	\includegraphics[width=0.3\linewidth]{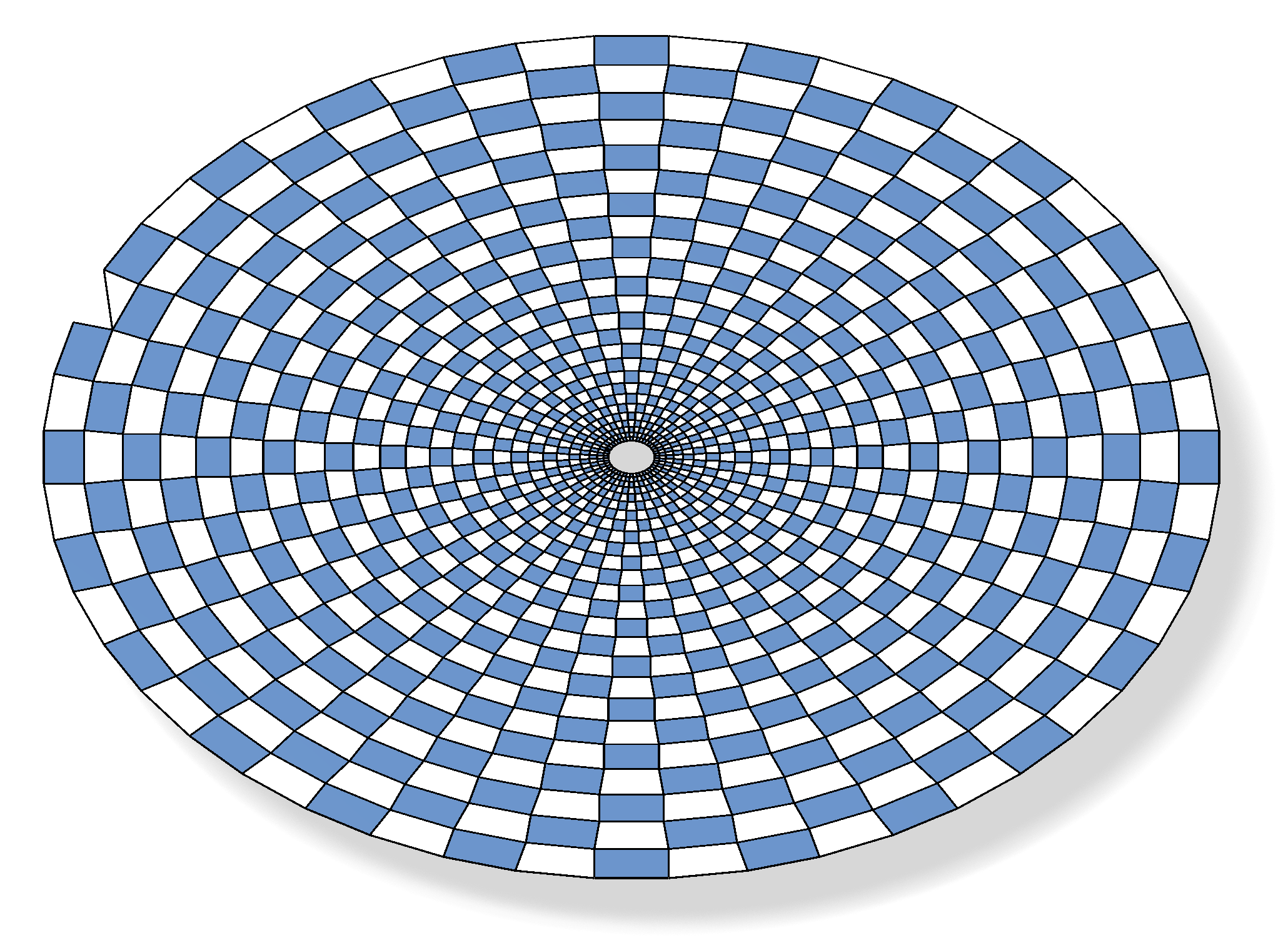}
	\includegraphics[width=0.3\linewidth]{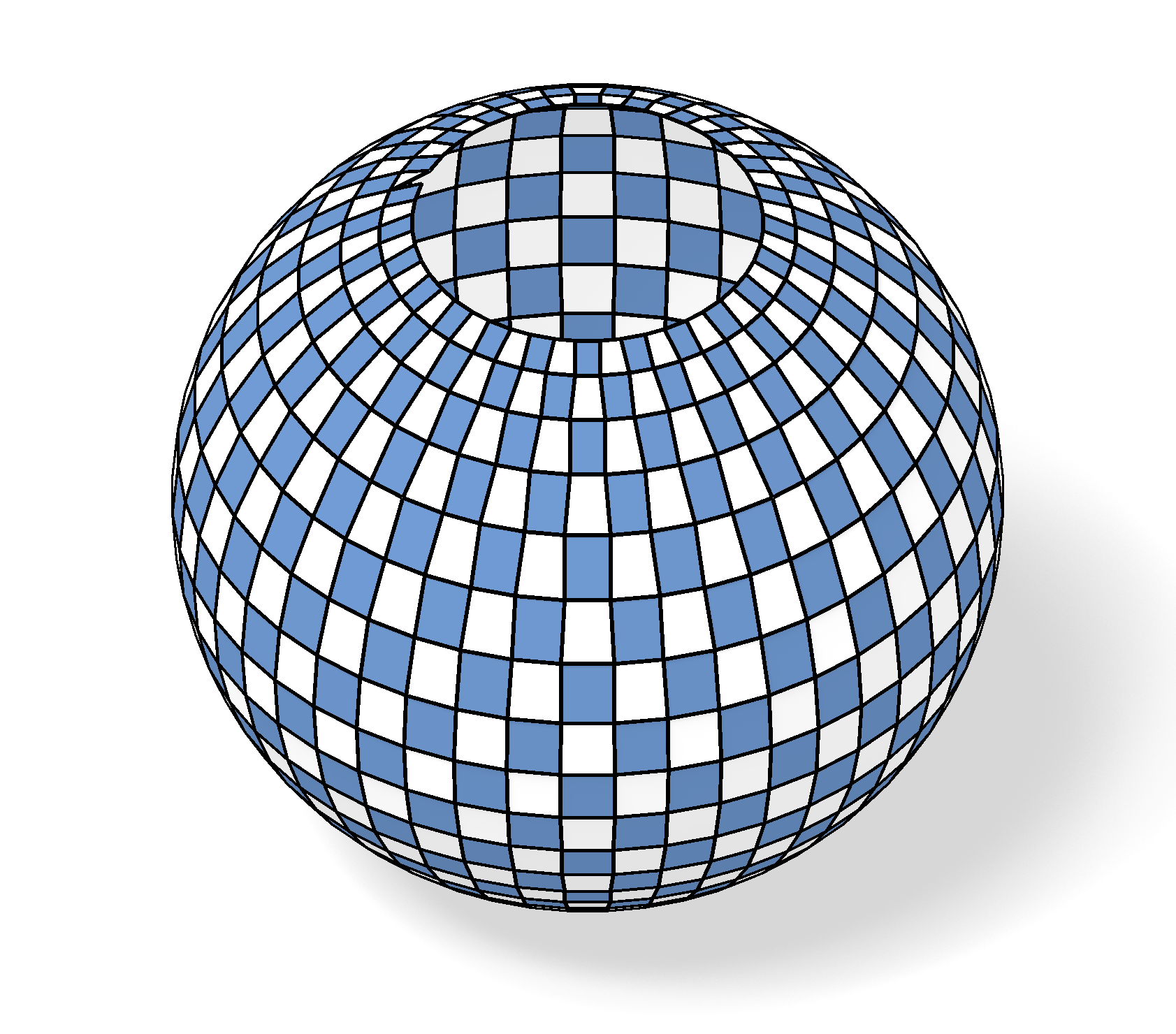}
	\includegraphics[width=0.3\linewidth]{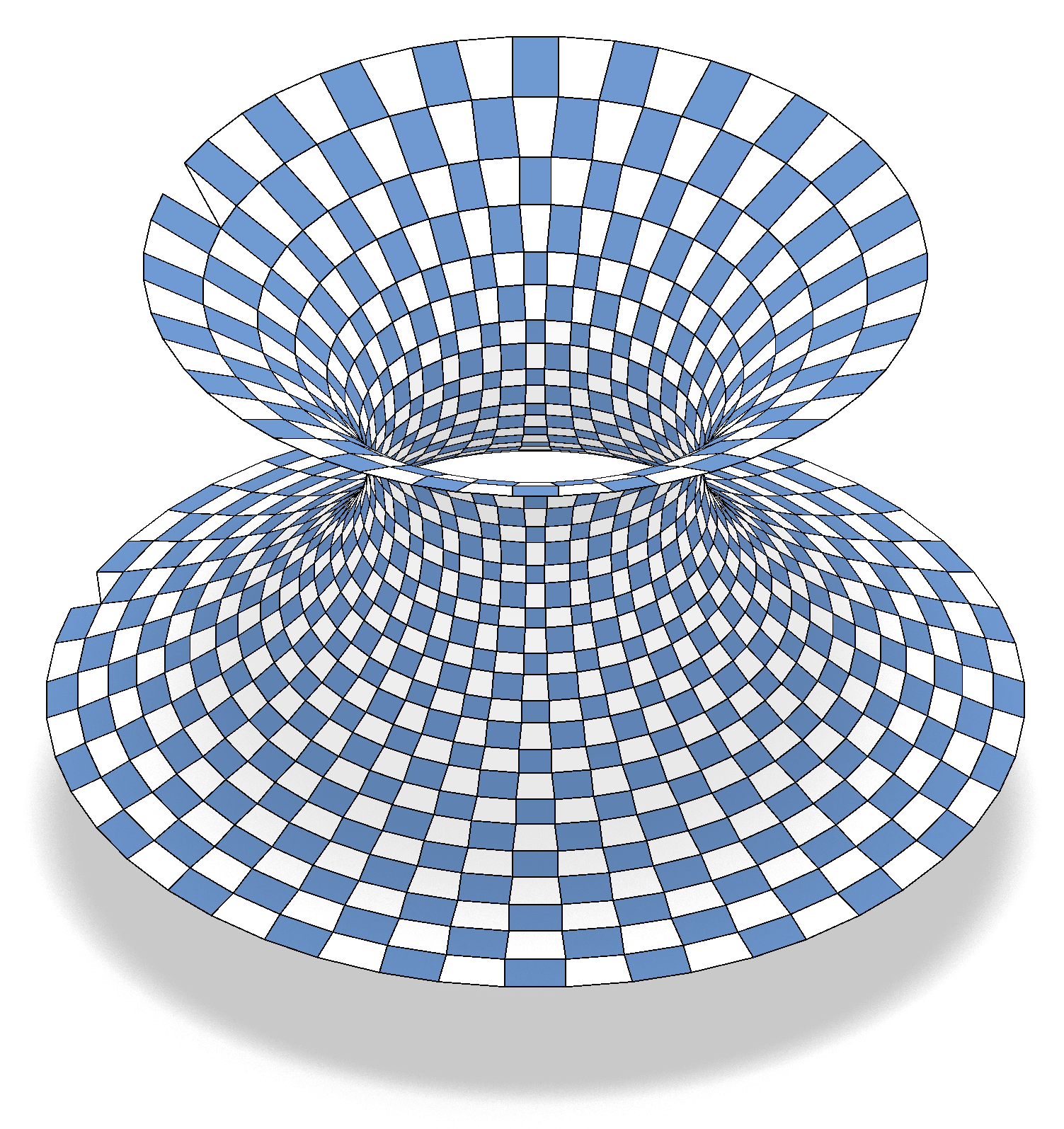}
	\caption{Catenoid: In the top row we see from left to right the Weierstrass data of the Catenoid, the Gauß image of the Catenoid and the Catenoid itself. In the second row we see the checkerboard patterns of the corresponding nets.}
	\label{Fig:Catenoid}
\end{figure}

\begin{figure}[h]
	\includegraphics[width=0.7\linewidth]{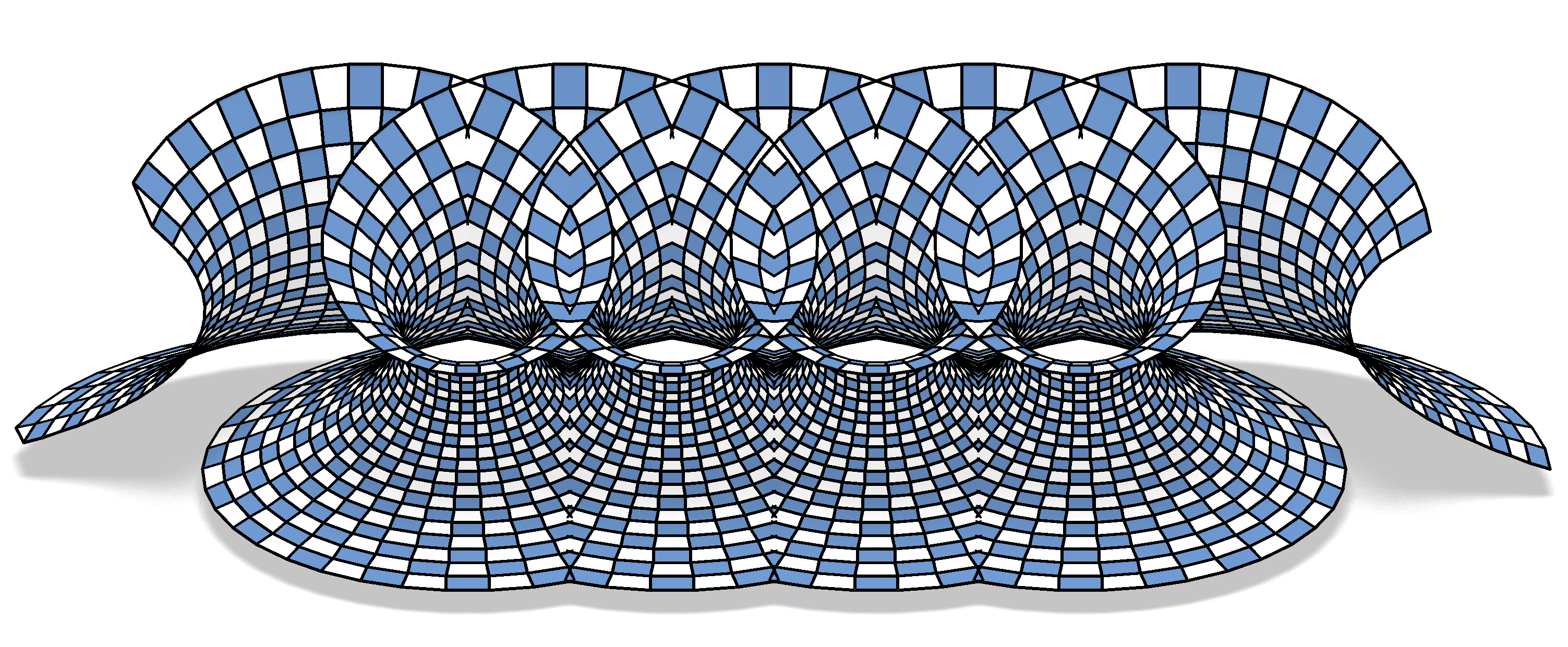}
		\caption{A Goursat transform of a periodically extended Catenoid.\vspace*{2em}}
\end{figure}

\section{Conclusion}

In this paper we presented a novel discretization approach based on the checkerboard pattern inscribed to a quadrilateral net. On the one hand this allows a discrete curvature theory (Definition \ref{Defi:ShapeOperator}) that is compatible with discrete offsets (Formula \ref{SteinerFormel}) similar to \cite{MixedAreaHannesChristian, hoffmann2014discrete}. On the other hand this approach allows a new discretization of conjugate nets, orthogonal nets and principal nets (Definition \ref{Defi:OrthoConjuPrinciNet}). We showed several properties of these nets, most noticeably that principal nets are consistent with the curvature theory (Corollary \ref{Kor:PrincipalDirections}) and are invariant under Möbius transformations (Theorem \ref{MoebiusInvariance}) applied to the corresponding sphere congruence introduced in \cite{techter2021discrete}.

Further the checkerboard pattern could be used to define discrete Koenigs nets using the conic of Koenigs (Definition \ref{Defi:KoenigsNet}) analogous to \cite{DoliwasKoenigsNets}. We find that discrete Koenigs nets are exactly those nets that are dualizable (Theorem \ref{Theorem:KoenigsNets}) which links the approach taken in \cite{DoliwasKoenigsNets} to the approach of \cite{KoenigsNetsBobenkoSuris}. Other characterizations of discrete Koenigs nets that have been found in this paper are the existence of a closed multiplicative one-form defined on the edges of a checkerboard pattern (Theorem \ref{Theorem:ClosedOneForm}) similar to \cite{KoenigsNetsBobenkoSuris}. A new characterization of Koenigs nets that has been found is the equality of Laplace invariants (Theorem \ref{Theorem:LapInvariants}) which fits the original definition of these nets in the classical differential geometry. From the characterization via equal Laplace invariants we could deduce that the class of discrete Koenigs nets is invariant under projective transformations (Corollary \ref{Kor:ProjectiveInvariance}).

Despite the discretization idea of Koenigs nets and principal nets being quite different they work well together for isothermic nets which are defined as principal Koenigs nets. This means that the Koenigs property is preserved upon Möbius transformations (Theorem \ref{Theorem:MoebiusInvarianceOfIsothermic}) and the principality is preserved upon dualization (Corollary \ref{Kor:CFDual}). Consequently we can apply Möbius transformations and dualizations to discrete isothermic nets. This allows a construction of discrete minimal surfaces from an isothermic net in the plane. First we map it to the unit sphere with a Möbius transformation, where it can be interpreted as the Gauß image of a minimal surface. Then it is dualized to gain the corresponding minimal surface from its Gauß image, compare both Figure \ref{Fig:EnepperSurface} and Figure \ref{Fig:Catenoid}.


\section{Acknowledgments}

The author is grateful to C. Müller, J. Wallner and H. Pottmann for many inspiring discussion and extensive proofreading. Further the author gratefully acknowledges the support of the Austrian Science Fund (FWF) through projects P29981, I4868 and W1230 \grq DK Discrete Mathematics\grq.

\printbibliography

@article{CheckerboardIntro,
	year = {2019},
	doi = {10.1145/3355089.3356514},
	journal = {ACM Trans. Graphics},
	Volume=38,
	number =  6,
	pages = "171:1-13",
	author = {Peng, Chi-Han and Jiang, Caigui and Wonka, Peter and Pottmann, Helmut},
	title = {Checkerboard patterns with black rectangles},
	XXXpublisher = {Association for Computing Machinery (ACM)},
	XXXurl = {http://hdl.handle.net/10754/660118}
}

@book{BobenkoBuch,
	series = {Graduate studies in mathematics},
	volume= {98},
	xxxisbn = {9780821847008},
	year = {2008},
	title = {Discrete differential geometry: Integrable structure},
	XXlanguage = {eng},
	author = {Bobenko, Alexander and Suris, Yuri},
	publisher = {American Math. Soc.},
}

@article{DoliwasKoenigsNets,
	author = {Doliwa, Adam},
	year = {2003},
	XXXXmonth = {04},
	title = {Geometric discretization of the {K}oenigs nets},
	volume = {44},
	pages={2234-2249},
	journal = {Journal of Mathematical Physics},
	doi = {10.1063/1.1563041}
}

@article{MixedAreaHannesChristian,
	Author = {Christian {M\"uller} and Johannes {Wallner}},
	Title = {{Oriented mixed area and discrete minimal surfaces}},
	FJournal = {{Discrete \& Computational Geometry}},
	Journal = {{Discrete Comput. Geom.}},
	XXISSN = {0179-5376; 1432-0444/e},
	Volume = {43},
	Number = {2},
	Pages = {303-320},
	Year = {2010},
	XXPublisher = {Springer},
	XXaddress = {US, New York, NY},
	XXLanguage = {English},
	MSC2010 = {52B70 53A10},
	Kommentar = {msc heisst mathematical subject classification},
	Zbl = {1192.52023},
	Kommentar2 = {zbl heisst zentralblatt fuer mathematik},
	XXurl = {https://doi.org/10.1007/s00454-009-9198-7},
	doi = {10.1007/s00454-009-9198-7}
}

@article{hoffmann2014discrete,
	title={A discrete parametrized surface theory in $\mathbb{R}^3$}, 
	author={Tim Hoffmann and Andrew O. Sageman-Furnas and Max Wardetzky},
	year={2017},
	volume={2017},
	number={14},
	pages={4217-4258},
	fulljournal = {International Mathematics Research Notices},
	journal = {Int. Mat. Res. Notices},
	doi = {10.1093/imrn/rnw015}
}

@phdthesis {techter2021discrete,
	author = {Techter, Jan},
	title = {Discrete confocal quadrics and checkerboard incircular nets},
	school = {Technische Universit{\"a}t Berlin},
	year = {2021},
	XXtype = {Doctoral Thesis},
	XXaddress = {Berlin},
	doi = {10.14279/depositonce-11461},
	XXurl = {http://dx.doi.org/10.14279/depositonce-11461},
}

@book{lane1932projective,
	title={Projective differential geometry of curves and surfaces},
	author={Lane, Ernest Preston},
	year={1932},
	publisher={The University of Chicago Press}
}

@article{liu-Pottmann-conical,
	author = {Yang Liu and Helmut Pottmann and Johannes Wallner
	and Yong-Liang Yang and Wenping Wang},
	title = {Geometric modeling with conical meshes and developable
	surfaces},
	journal = {ACM Trans. Graphics},
	volume = 25,
	number = 3,
	year = 2006,
	XXissn = {0730-0301},
	pages = {681--689},
	XXnote = "Proc. SIGGRAPH",
	doi = "10.1145/1141911.1141941",
}

@article{pottmann-2007-pm,
	author = "Helmut Pottmann and Yang Liu and Johannes Wallner
	and Alexander Bobenko and Wenping Wang",
	title = "Geometry of Multi-layer Freeform Structures for Architecture",
	journal = {ACM Trans. Graphics},
	year = 2007,
	volume = 26,
	number = 3,
	XXXnote = "Proc. SIGGRAPH",
	XXXurl = "http://www.geometrie.tugraz.at/wallner/parallel.pdf",
	doi = "10.1145/1275808.1276458",
	pages = "65:1-11",
}

@article{KoenigsNetsBobenkoSuris,
	author = {Bobenko, Alexander I. and Suris, Yuri B.},
	title = {Discrete Koenigs Nets and Discrete Isothermic Surfaces},
	journal = {International Mathematics Research Notices},
	volume = {2009},
	number = {11},
	pages = {1976-2012},
	year = {2009},
	XXXmonth = {02},
	XXXissn = {1073-7928},
	doi = {10.1093/imrn/rnp008},
	XXXurl = {https://doi.org/10.1093/imrn/rnp008},
	XXXeprint = {https://academic.oup.com/imrn/article-pdf/2009/11/1976/2070263/rnp008.pdf},
}

@article{KoenigsOriginalPaper,
	author = {Koenigs, Gabriel},
	title = {Sur les r\'{e}seaux plans à invariants \'{e}gaux et les lignes asymptotiques},
	journal = {Comptes Rendus de l\'{}Academie des Sciences, S\'{e}rie 1: Math\'{e}matique},
	volume = {114},
	year = {1892},
	pages = {55-57},	
}

@book{BraunerGPR,
	author = "Brauner, Heinrich",
	year = "1976",
	publisher = "Bibliographisches Institut",
	title = "Geometrie projektiver Räume I",
}

@incollection{weingarten-aag-2020,
	title = "Architectural freeform surfaces designed for cost-effective
	paneling mold re-use",
	author = "Davide Pellis and Martin Kilian and Hui Wang and Caigui Jiang
	and Christian M{\"u}ller and Helmut Pottmann",
	year = 2020,
	booktitle = "Advances in Architectural Geometry 2020",
	address = "Paris",
	publisher = "{\'E}cole des Ponts",
}

@article{UsingIsometries,
	author = "Caigui Jiang and Hui Wang and Ceballos Inza, Victor
	and Felix Dellinger and Florian Rist
	and Johannes Wallner and Helmut Pottmann",
	year = 2021,
	title = "Using isometries for computational design and fabrication",
	journal = "ACM Trans. Graph.",
	volume = 40,
	number = 4,
	pages = "42:1-12",
	doi = "10.1145/3450626.3459839",
	XXurl = "http://www.geometrie.tugraz.at/wallner/isopanel.pdf",
	XXnote = "Proc. SIGGRAPH",
}

@article{Jiang2020Apl1,
	author = {Jiang, Caigui and Wang, Cheng and Rist, Florian and Wallner, Johannes and Pottmann, Helmut},
	title = {Quad-mesh based isometric mappings and developable surfaces},
	year = {2020},
	XXpublisher = {Association for Computing Machinery},
	XXaddress = {New York, USA},
	volume = {39},
	number = {4},
	XXissn = {0730-0301},
	XXXurl = {https://doi.org/10.1145/3386569.3392430},
	doi = {10.1145/3386569.3392430},
	journal = {ACM Trans. Graph.},
	XXmonth = jul,
	XXarticleno = {128},
	XXnumpages = {13},
	pages = {128:1-13}
}


\appendix

\section{Proof of Theorem \ref{Theorem:ClosedOneForm}}


\begin{bew}
	First we show that $q$ is closed around every second order face, i.e., multiplying the contribution of every edges of a second order face in counter clockwise order yields $1$. Let $c,c_1,c_{12}$ and $c_2$ be the vertices of a second order face, $p = \suppl{c}{c_1} \cap \suppl{c_2}{c_{12}}$ and $q = \suppl{c}{c_2}\cap \suppl{c_1}{c_{12}}$, see Figure \ref{Fig:ExactnessWhite}.
	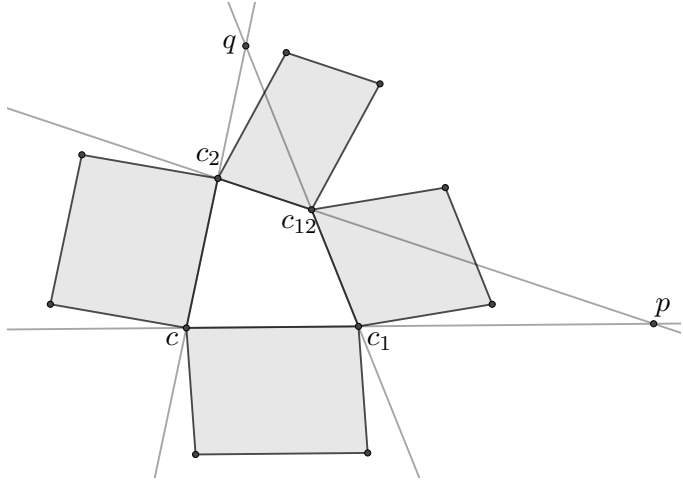
\begin{figure}[h]
		\begin{overpic}[width=0.68\linewidth]{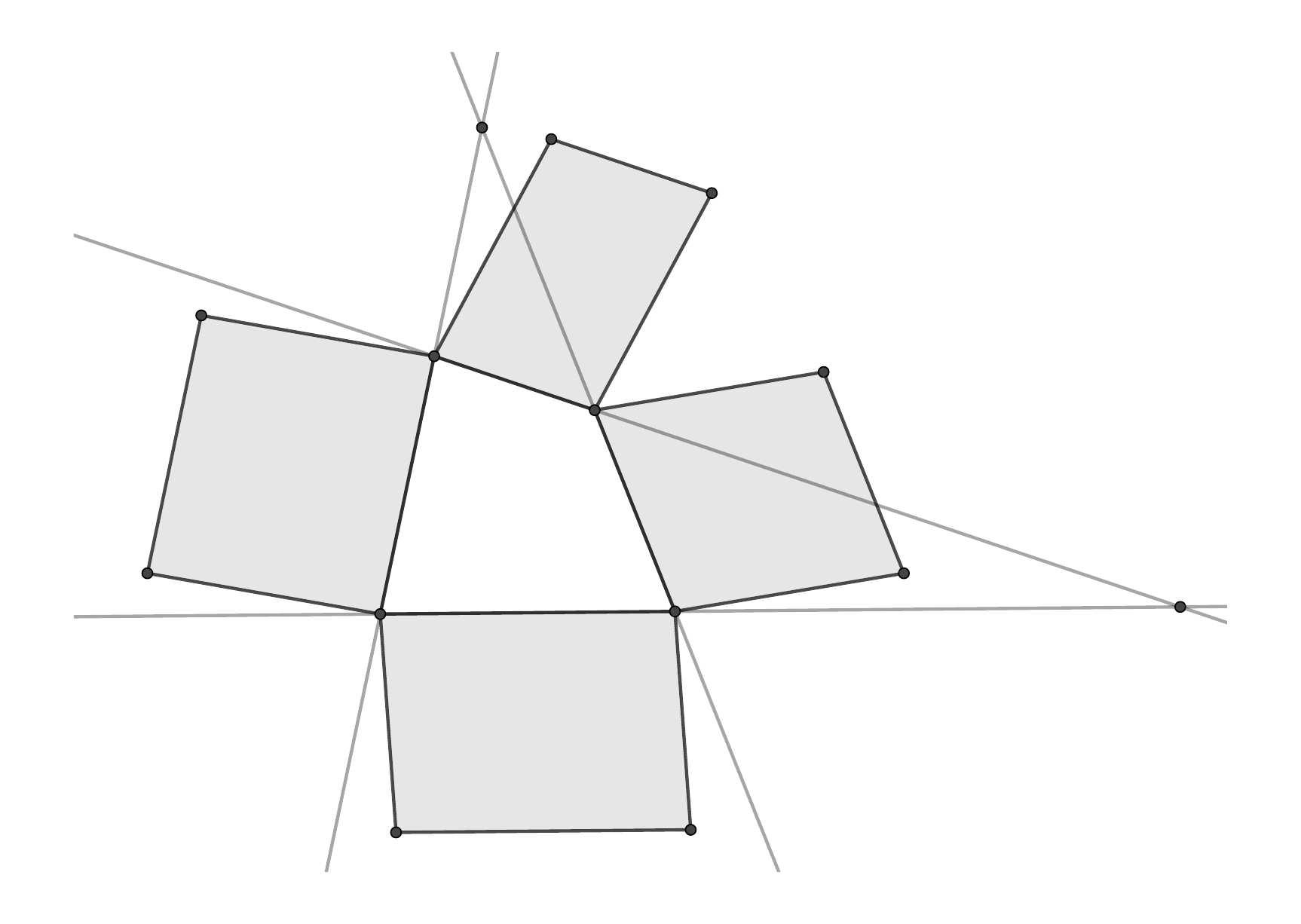}
			\put(26.5,21.5){$c$}
			\put(53,21.5){$c_1$}
			\put(41.7,37){$c_{12}$}
			\put(30.5,46){$c_2$}
			\put(91,26){$p$}
			\put(34,61){$q$}
		\end{overpic}
		\hfill
		\begin{minipage}[b]{0.3\linewidth}
			\caption{The multiplicative one form is automatically closed around every second order face as a consequence of Menelaus' Theorem.\vspace*{4em}}
			\label{Fig:ExactnessWhite}
		\end{minipage}
	\end{figure}
	Using Menelaus' Theorem \ref{Theorem:Menelaus} for the triangle $(c,c_1,q)$ and the triangle $(c_2,c_{12},q)$ we find that
	\begin{align*}
	&q(c,c_1) q(c_1,c_{12}) q (c_{12},c_2) q (c_2,c) = \\
	= &\frac{c - p}{c_1 - p} \frac{c_1 - q}{c_{12} - q} \frac{c_{12} - p}{c_{2} - p} \frac{c_2 - q}{c - q} =\\
	= &\underset{= -1}{\underbrace{\frac{c - p}{p - c_1} \frac{c_{12} - c_1}{q - c_{12}}
	\frac{c_2 - q}{c - c_2}}}
	\underset{=-1}{\underbrace{\frac{c - c_2}{q - c} 
			\frac{c_1 - q}{c_{12} - c_1} 
			\frac{p - c_{12}}{c_{2} - p}}} 
	= 1.
	\end{align*}
	Next we show that $q$ is closed on the edges of every first order face $\blackface = (c,c_1,c_{12},c_2)$ if and only if the checkerboard pattern is a Koenigs net. As the multiplicative one-form $q$ is projectively invariant, we choose a projective coordinate system such that $c = (0,0,1)$, $c_1 = (1,0,1)$, $c_2=(0,1,1)$ and $c_{12} = (1,1,1)$. The intersection points then have the following coordinates
	\begin{align*}
	p_1 &= \suppl{c}{c_1} \cap \suppl{c_2}{c_{12}} = (1,0,0)\\
	p_2 &= \suppl{c}{c_2} \cap \suppl{c_1}{ c_{12}} = (0,1,0)\\
	p_3 &= \suppl{c}{c_1} \cap \suppl{c_{-2}}{c_{1-2}} = (1,0,t)\\
	p_4 &= \suppl{c_2}{c_{12}} \cap \suppl{c_{22}}{c_{122}} = (s,1,1) \\
	p_5 &= \suppl{c}{c_2} \cap \suppl{c_{-1}}{c_{-12}} = (0,1,u)\\
	p_6 &= \suppl{c_1}{c_{12}} \cap \suppl{c_{11}}{c_{112}} = (1,v,1)
	\end{align*}
	for suitable $s,v \in \R$. Those six points lie on a common conic section if the system of equations $Ax_i^2 + Bx_iy_i + Cy_i^2 + D x_iz_i + E y_iz_i + F z_i^2 = 0$ has a nontrivial solution. Here $x_i,y_i$ and $z_i$ stand for the three homogeneous coordinates of $p_i$. We compute the determinant of the matrix of this system of equations:
	\begin{align*}
	&\det
	\begin{pmatrix}
	1 & 0 & 0 & 0 & 0 & 0\\
	0 & 0 & 1 & 0 & 0 & 0\\
	1 & 0 & 0 & t & 0 & t^2\\
	s^2 & s & 1 & s & 1 & 1\\
	0 & 0 & 1 & 0 & u & u^2\\
	1 & v &v^2& 1 & v & 1 
	\end{pmatrix} =
	t \left( s (u (1-t) - u^2 v) + v(u^2 - u(1-st)) \right).
	\end{align*}	
	A nontrivial solution exists if and only if the determinant is zero. We can exclude the cases $t=0$ and $u=0$ since no $p_i$ are the same. We find that the determinant is zero if and only if
	\begin{align*}
	s - st + vst &= v - vu + vus.
	\end{align*}
	Now we compute the multiplicative one-form along the edges of the quadrilateral. We find $p_3 = c_1 + (t-1) c$ and $p_1 = c_1 - c$. So if we use $c$ and $c_1$ as the bases for the line $\suppl{c}{c_1}$, we obtain
	\begin{align*}
	\DV(c,c_1,p_1,p_3) &= \frac{|c,p_1|}{|c_1, p_1|} \frac{|c_1,p_3|}{|c,p_3|}  =   \frac{\det \begin{bmatrix}
		1 & -1\\ 0 & 1
		\end{bmatrix}\det \begin{bmatrix}
		0 & t-1\\ 1 & 1
		\end{bmatrix} }{
		\det \begin{bmatrix}
		0 & -1\\
		1 & 1
		\end{bmatrix} \det \begin{bmatrix}
		1 & t-1\\
		0 & 1
		\end{bmatrix} } = 1 - t.
	\end{align*}
	For the next cross ratio we express the points $p_2$ and $p_6$ via $c_1$ and $c_{12}$, obtaining $p_2 = -c_1 + c_{12}$ and $p_6 =  (1-v) c_1 +  v c_{12}$. So the cross ratio is
	\begin{align*}
	\DV(c_1, c_{12},p_2, p_6) &= \frac{|c_1,p_2|}{|c_{12}, p_2|} \frac{|c_{12}, p_6|}{|c_1,p_6|} = \frac{
		\det \begin{bmatrix}
		1 & -1\\
		0 & 1
		\end{bmatrix}
		\det \begin{bmatrix}
		0 & 1-v\\
		1 & v
		\end{bmatrix}}{
		\det \begin{bmatrix}
		0 & -1\\
		1 & 1
		\end{bmatrix}
		\det \begin{bmatrix}
		1 & 1-v\\
		0 & v
		\end{bmatrix}} = \frac{v-1}{v}.
	\end{align*}
	Next, equation $p_1 = c_{12} -c_2$ and $p_4 = s c_{12} + (1-s) c_2$ yield
	\begin{align*}
	\DV(c_{12},c_2,p_1,p_4) &= \frac{|c_{12},p_1|}{|c_{2}, p_1|} \frac{|c_{2}, p_4|}{|c_{12},p_4|} = \frac{
		\det \begin{bmatrix}
		1 & 1\\
		0 & -1
		\end{bmatrix}
		\det \begin{bmatrix}
		0 & s\\
		1 & 1-s
		\end{bmatrix}}{
		\det \begin{bmatrix}
		0 & 1\\
		1 & -1
		\end{bmatrix}
		\det \begin{bmatrix}
		1 & s\\
		0 & 1-s
		\end{bmatrix}} = \frac{s}{s-1}.
	\end{align*}
	and $p_2 = c_2 - c$ and $p_5 = c_2 + (u-1)c$ yield
	\begin{align*}
	\DV(c_2,c,p_2,p_5) &= \frac{|c_2,p_2|}{|c,p_2|} \frac{|c,p_5|}{|c_2,p_5|} =
	\frac{
		\det \begin{bmatrix}
		1 & 1\\
		0 & -1
		\end{bmatrix}
		\det \begin{bmatrix}
		0 & 1\\
		1 & u-1
		\end{bmatrix}}{
		\det \begin{bmatrix}
		0 & 1\\
		1 & -1
		\end{bmatrix}
		\det \begin{bmatrix}
		1 & 1\\
		0 & u-1
		\end{bmatrix}} = \frac{1}{1-u}.
	\end{align*}
	Now $q$ is closed if and only if
	\begin{align*}
	1 &= (1-t)\frac{v-1}{v} \frac{s}{s - 1} \frac{1}{1 - u} \\
\iff \quad	s - st + stv &= v - uv + vsu.
	\end{align*}
	Thus the existence of the conic of Koenigs is equivalent to $q$ being closed.
\end{bew}

\section{Some Theorems}

The following lemma is known as trace polarity of a quadric. For a detailed description of trace polarity in German language see \cite{BraunerGPR}. However, since a proof in English of the following lemma is hard to find, we give a version of a proof suitable to our setting.

\begin{lem}\label{Lemma:Conjugacy}
	Let $s_1$ and $s_2$ be two orthogonally intersecting spheres on $\S^3$, i.e. intersections of $\S^3$ with conjugate hyperplanes $h_1$ and $h_2$. If $\psi$ is the central projection from point $Z\in\Pro\R^{4}$ onto a hyperplane $\zeta \cong \Pro\R^3$, then in this hyperplane $\psi(s_1)$ and $\psi(s_2)$ intersect in conjugate tangent planes with respect to the contour quadric $\psi(\S^3)^*$ of $\psi(\S^3)$. This means the corresponding tangent planes at the intersection points of $\psi(s_1)$ and $\psi(s_2)$ are orthogonal with respect to the inner product induced by $\psi(\S^3)^*$.
\end{lem}
\begin{bew}
	For the proof we use homogeneous coordinates of $\Pro\R^4$. We choose a basis $(b_0,\dots, b_4)$ such that the center of projection $Z = b_0$. Let $Q$ be the matrix such that the homogeneous coordinates of all points in $\S^3$ are given by $\{x\in \R^5 \backslash \{0\}\  : x^T Q x = 0\}$. We can assume without loss of generality that $Q$ is a diagonal matrix. Since projections onto different planes are projectively equivalent, we can further assume that $\zeta$ is the polar hyperplane of $Z$. Thus $\zeta$ is given by the equation $x_0 = 0$. We introduce the block notation $Q = \operatorname{diag}(q_0,\bar{Q})$.
	
	Since $\zeta$ is the polar hyperplane of $Z$, the contour quadric $\psi(\S^3)^*$ is the intersection of $\S^3$ with $\zeta$. In homogeneous coordinates it is given by $\{(0,\bar{x}) \in \R^5 : \bar{x}^T \bar{Q} \bar{x} = 0\}$.
	
	Let $P \in s_1 \cap s_2$ be a point in the intersection of $s_1$ and $s_2$ and let $\tau$ be the tangent hyperplane to $\S^3$ in $P$. Then the tangent planes to $s_1$ and $s_2$ are given by $\tau_1 = h_1 \cap \tau$ and $\tau_2 = h_2 \cap \tau$. Note that two hyperplanes are conjugate with respect to a quadric, if and only if each contains the polar point of the other. Let $H_1 \in h_2 \cap \tau$ be the polar point of $h_1$ and let $H_2 \in h_1 \cap \tau$ be the polar point of $h_2$. The line $g_1 := \suppl{P}{H_2}$ lies in $\tau_1$ and intersects $\zeta$. We denote the intersection point by $G_1 := g_1 \cap \zeta$. Analogously we define $G_2 := \suppl{P}{H_1} \cap \zeta$. The line $g:= \tau_1 \cap \tau_2$ also intersects $\zeta$ and we denote the intersection point by $T$.
	
	Writing $A \lor B \lor C$ for the plane spanned by $A,B,C$, we have $\tau_1 = P \lor T \lor G_1$ and $\tau_2 = P \lor T \lor G_2$. Since $T, G_1, G_2 \in \zeta$, we find that $\psi(\tau_1) = \psi(P) \lor T \lor G_1$ and $\psi(\tau_2) = \psi(P) \lor T \lor G_2$.
	
	We now show the orthogonality of $\psi(\tau_1)$ and $\psi(\tau_2)$ with respect to $\bar{Q}$. In order to facilitate the notation we identify points with their projective coordinates. As the projection $\psi$ just sets the first coordinate of $P$ to zero and due to the form of $Q$, we find that
	\begin{align*}
		G_1^T \bar{Q} \psi(P) &= G_1^T Q P = (\lambda_1 H_2 + \mu_1 P)^T Q P = 0,\\
		G_2^T \bar{Q} \psi(P) &= G_2^T Q P = (\lambda_2 H_1 + \mu_2 P)^T Q P = 0,\\
		G_1^T \bar{Q} T &= G_1^T Q T = (\lambda_1 H_2 + \mu_1 P)^T Q T = 0,\\
		G_2^T \bar{Q} T &= G_2^T Q T = (\lambda_2 H_1 + \mu_2 P)^T Q T = 0,\\
		G_1^T \bar{Q} G_2 &= G_1 Q G_2 = 0.
	\end{align*}
	Hence, the point $G_1$ is the polar point of $\psi(\tau_2)$ with respect to $\psi(\S^3)^*$ in $\zeta$ and vice versa. This shows the conjugacy of the tangent planes $\psi(\tau_1)$ and $\psi(\tau_2)$ of the projected spheres $s_1$ and $s_2$.
\end{bew}

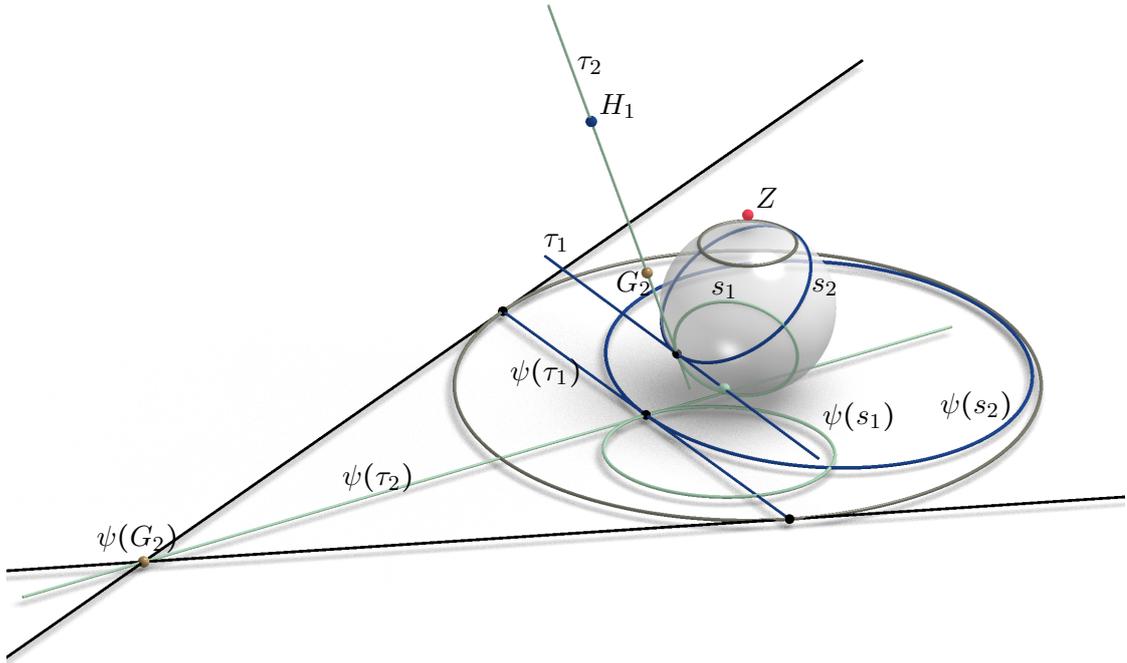
\begin{figure}[h]
	\begin{overpic}[width=\linewidth]{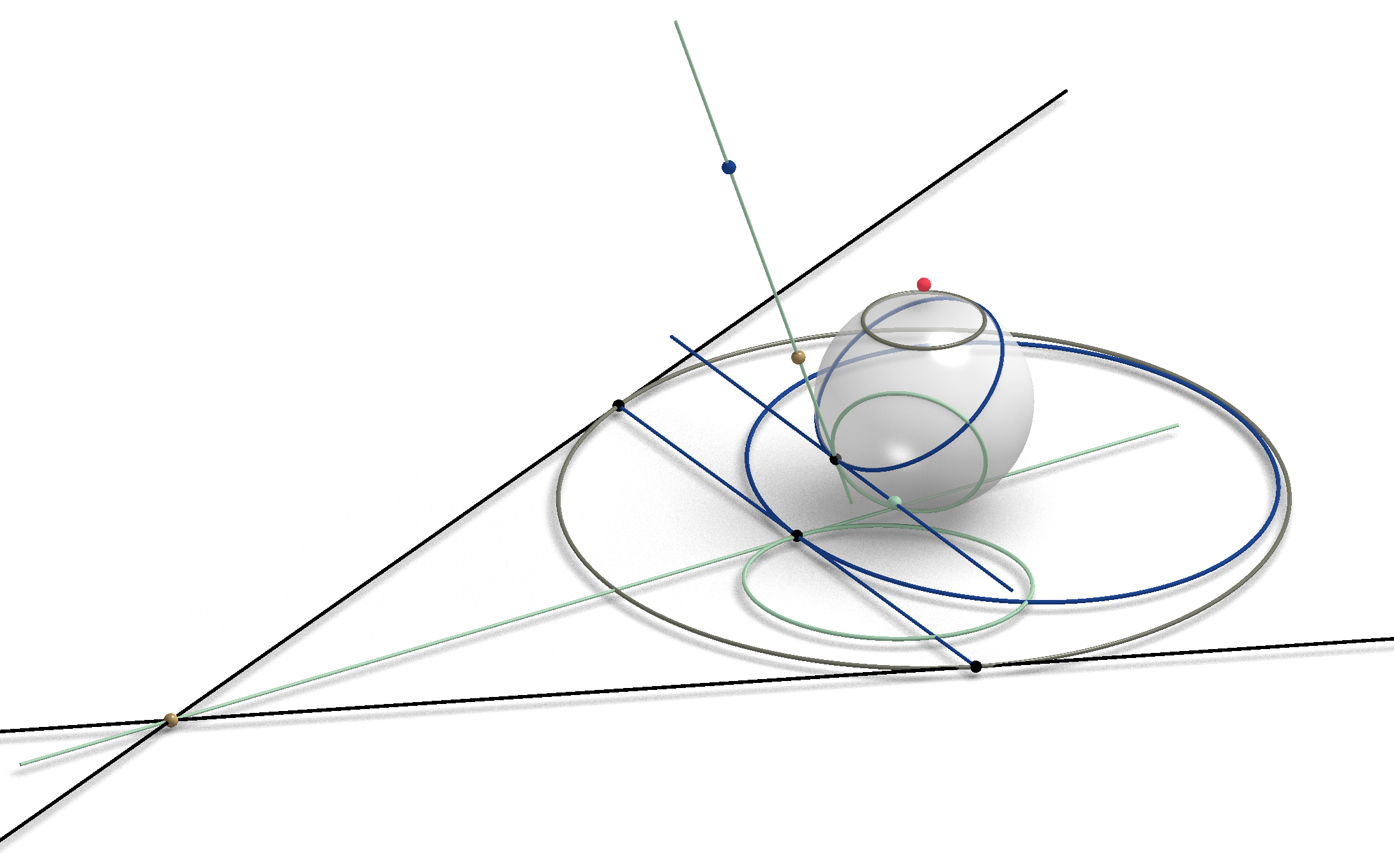}
		\put(8,12){$\psi(G_2)$}
		\put(54.5,35){$G_2$}
		\put(53,51){$H_1$}
		\put(67,42.5){$Z$}
		\put(48,39){$\tau_1$}
		\put(45,27){$\psi(\tau_1)$}
		\put(51,55){$\tau_2$}
		\put(30,17.5){$\psi(\tau_2)$}
		\put(83.5,24){$\psi(s_2)$}
		\put(72,35){$s_2$}
		\put(63,34.7){$s_1$}
		\put(73,23){$\psi(s_1)$}
	\end{overpic}
	\caption{The two-dimensional case of Lemma \ref{Lemma:Conjugacy}. The circles $s_1$ and $s_2$ on the unit sphere intersect orthogonally. The ellipses are their projections through the point $Z$. The gray ellipse $\psi(\S^2)^*$ is the contour of the unit sphere under the same projection.
	We see that the ellipses $\psi(s_1)$ and $\psi(s_2)$ intersect in conjugate lines with respect to $\psi(\S^2)^*$ as the polar point $\psi(G_2)$ of the tangent $\psi(\tau_2)$ is contained in $\psi(\tau_2)$. The preimage of this polar point drawn in beige is the intersection of the corresponding tangent line to the unit sphere with the polar plane of the center of projection. The gray circle on the unit sphere is the preimage of the contour quadric, i.e., the intersection of the unit sphere with the polar plane of $Z$.}
\end{figure}

\begin{satz}[Inscribed Angle Theorem for Hyperbolas]\label{InscribedAngleTheoremHyp}
	Consider $\R^2$ and coordinates $(x,y)$ with respect to a basis. The slope of a vector $(a,b)$ is defined as $\frac{b}{a}$ for $a\neq 0$. Four points $p_i = (x_i,y_i) \in \R^2$ with $x_j \neq x_k$ and $y_j \neq y_k$ lie on a hyperbola with equation $y = \frac{c}{x}$, if and only if the quotient of slopes of $\suppl{p_4}{ p_2}$ and $\suppl{p_4}{p_1}$ equals the quotient of slopes of $\suppl{p_3}{p_2}$ and $\suppl{p_3}{ p_1}$, compare Figure \ref{Fig:Peripheriewinkelsatz}. Computing this condition yields
	\begin{align}\label{Eq:AngleCondition}
	\frac{(y_4 - y_1)(x_4 - x_2)}{(x_4 - x_1)(y_4 - y_2)} = \frac{(y_3 - y_1)(x_3 - x_2)}{(x_3 - x_1)(y_3 - y_2)}.
	\end{align}
	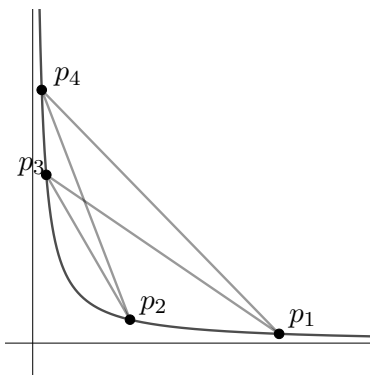
\begin{figure}
	\begin{overpic}[width=0.4\linewidth]
		{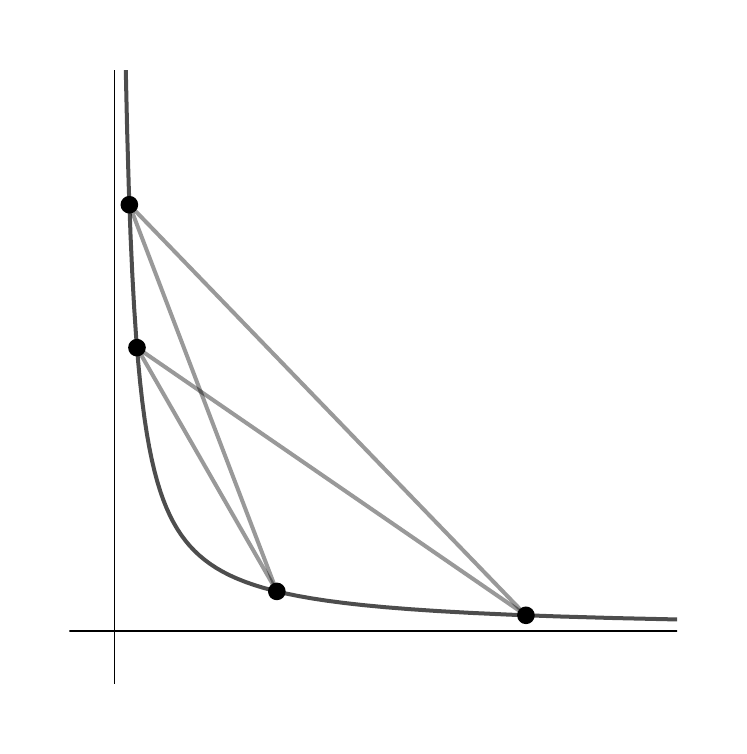}
		\put(72,21){$p_1$}
		\put(39,24){$p_2$}
		\put(12,55){$p_3$}
		\put(20,75){$p_4$}
	\end{overpic}
	\hfill
	\begin{minipage}[b]{0.4\linewidth}
		\caption{Inscribed angle theorem for hyperbolas. The four points $p_1,p_2,p_3$ and $p_4$ lie on a rectangular hyperbola if and only if the quotient of slopes of $\overline{p_4 p_2}$ and $\overline{p_4,p_1}$ equals the one of $\overline{p_3p_2}$ and $\overline{p_3p_1}$.\vspace*{3em}}
		\label{Fig:Peripheriewinkelsatz}
	\end{minipage}
	\end{figure}	
\end{satz}

\begin{satz}[Menelaus' Theorem]
	\label{Theorem:Menelaus}
	Let $A,B$ and $C$ be the vertices of a triangle and let $g$ be a straight line. For the three vertices $D = \suppl{A}{B} \cap g$, $E = \suppl{B}{C} \cap g$ and $F = \suppl{C}{A} \cap g$ the equation
	\begin{align*}
	\frac{(D-A)}{(B-D)} \cdot \frac{(E-B)}{(C - E)} \cdot \frac{(F-C)}{A - F} = -1
	\end{align*}
	holds. 
\end{satz}

\end{document}